\documentclass[a4paper,12pt]{article}
\usepackage{authblk,graphicx,amsthm,amsmath,amssymb,mathrsfs,amscd,bbm,float,booktabs,multirow,enumitem,hyperref,enumerate,amsmath,subfig,cite,sectsty,tikz,pgfplots} \usepackage{mathpazo}\usepackage[T1]{fontenc} \usepackage[utf8]{inputenc}
\usetikzlibrary{arrows.meta}\usetikzlibrary{decorations.pathmorphing}\usetikzlibrary{decorations.pathmorphing} \usepgfplotslibrary{fillbetween}
\usetikzlibrary{arrows} \pgfplotsset{compat=1.18}
\usepackage{times}
\newtheorem{prop}{Proposition}\newtheorem{theo}{Theorem}\newtheorem{Assu}{Assumption}\newtheorem{rhp}{RHP}\newtheorem{Lemma}{Lemma}\newtheorem{Remark}{Remark}\newtheorem{Drhp}{$\overline{\partial}$-RHP}\newtheorem{D}{$\overline{\partial}$-Problem}\newtheorem{Corollary}{Corollary}
\hypersetup{citecolor=blue}
\DeclareMathOperator*{\Res}{Res}\DeclareMathOperator*{\sech}{sech}\DeclareMathOperator*{\RRe}{Re}\DeclareMathOperator*{\IIm}{Im}
\makeatletter
\hypersetup{colorlinks=true,linkcolor=blue}
\date{}
\usepackage[a4paper, left=2.0cm, right=2cm, top=1cm, bottom=1.5cm]{geometry}
\sectionfont{\fontsize{14}{16}\selectfont}\subsectionfont{\fontsize{12}{14}\selectfont}
\numberwithin{equation}{section}
\title{ Long-time behavior of the reduced Maxwell-Bloch equations in the sharp-line limit}
\author{
	{\normalsize Kang Wu\textsuperscript{1}, ~~
	Yingcan Huang\textsuperscript{2}, ~~
	Jingsong He\textsuperscript{2}\thanks{Corresponding author: hejingsong@szu.edu.cn} \\
	\textsuperscript{1} School of Mathematical Sciences, Shenzhen University, Shenzhen 518060, PR China \\
	 \textsuperscript{2} Institute for Advanced Study, Shenzhen University, Shenzhen 518060, PR China}
}
 
\begin{document}
	{\normalsize  \maketitle}
	\begin{abstract}
		We study the Cauchy problem for the reduced Maxwell-Bloch equations with initial data for the electric field in weighted Sobolev spaces, assuming that all atoms initially reside in their ground state. Using the $\overline \partial$-steepest descent method, we derive long-time asymptotic expansions of the solutions, including both the electric field and the components of the Bloch vector, within any fixed cone. In particular, we formulate the inverse scattering transform as a properly posed Riemann-Hilbert problem, avoiding singularities in the scattering data by modifying the time evolution of the reflection coefficient. Under assumptions that allow only soliton generation, the leading-order asymptotics are determined by solitons inside the cone, while soliton-radiation interactions appear in lower-order terms. These results extend the applicability of the nonlinear steepest descent method to integrable systems with singularities in the associated Lax pair.
	\end{abstract}

	\tableofcontents
	\section{ Introduction}
	%直接给出这个方程以及本文的目的
	In this paper, we study the long-time asymptotic behavior of solutions to the Cauchy problem for the reduced Maxwell-Bloch (RMB) equations in the sharp-line limit  
	\begin{align}\label{eq01}
		\begin{aligned}
			&E_{t}( x,t) =-s( x,t)
			&s_{x}( x,t) =E( x,t) u( x,t) +\mu r( x,t)\\
			&u_{x}( x,t) =-E( x,t) s( x,t)
			&r_x( x ,t) =-\mu s( x,t), 
		\end{aligned}
	\end{align}
	with initial condition $E(x,0) =E_0(x).$ Here,  $E(x,t)\in \mathbb R $ represents the electric field intensity. The components of the Bloch vector $(r(x,t),s(x,t),u(x,t))\in \mathbb R ^3$, satisfying $r^2+s^2+u^2=1$, are linear combinations of the elements from the atomic density matrix \cite{Eilbeck1972b,J.C.Eilbeck_1973}.
	
	%RMB方程的历史，以及在数学和物理中的重要性
Originally derived by Eilbeck \cite{J.C.Eilbeck_1973}, the RMB equations describe the interaction between a classical electromagnetic field and a dielectric medium composed of quantized two-level atoms. More precisely, under the assumption that backscattering effects are negligible, the RMB equations \eqref{eq01} can be derived from the inhomogeneously broadened Maxwell–Bloch system  via a suitable change of variables and by taking the sharp-line limit, in which all atoms are assumed to share the same resonant frequency $\mu$ \cite{Bullough1979,Eilbeck1972a, J.C.Eilbeck_1973,Gibbon1973}. This system has received widespread attention from both the mathematics and physics communities, owing to its strong connections with the self-induced transparency (SIT) and sine-Gordon (SG) equations, as well as the valuable insights they provide into the SIT phenomenon \cite{Bullough1979,Dodd,Maimistov}.  SIT, a central topic in nonlinear optics, refers to the nearly lossless propagation of an ultrashort, high-energy, coherent light pulse resonant with a two-level medium, provided its intensity exceeds a certain critical threshold \cite{PhysRev.183.457}.

From a mathematical perspective, the inhomogeneously broadened RMB equations, when analyzed under the slowly varying envelope approximation and transformed into a rotating frame, reduce to the classical SIT equations considered in \cite{MR475310}. Nevertheless, one should note that the RMB equations describe the electric field, whereas the SIT equations describe the envelope of a resonant carrier wave. 
Owing to this fundamental distinction, the soliton and breather solutions of the RMB equations provide a more accurate and general description of the $2\pi$ pulse observed in SIT phenomenon, particularly beyond the range of validity of the slowly varying envelope approximation \cite{J.C.Eilbeck_1973,Maimistov2001}. In practical applications, the atomic frequency distribution in inhomogeneously broadened systems is often complex and not easily calculative. In such cases, it is a prudent strategy to first analyze the sharp-line limit as a tractable approximation, and then appropriately modify the results to account for a wide variety of broadening effects \cite{Bullough1979, J.C.Eilbeck_1973, Caudrey1973}. On the other hand, by setting the common resonant frequency $ \mu = 0$,  an equivalent form of the RMB equations \eqref{eq01}, as studied in \cite{J.C.Eilbeck_1973,Caudrey1974} reduces precisely to the sharp-line SIT equations, thereby establishing a direct connection with the SG equation.  This connection provides a clear mathematical pathway to understanding the integrability of the RMB equations \eqref{eq01}.
	
	%本文方程的研究成果罗列
	Summarizing known results, the RMB equations \eqref{eq01} represent a completely integrable model in nonlinear optics, with their integrable structure studied in \cite{Gibbon1973,Aiyer1983,Grauel1984,Grauel1986}.  A variety of exact solutions to \eqref{eq01}, such as multi-solitons, periodic and rogue waves, and their interactions, have been obtained using methods like the Darboux transformation and consistent Riccati expansion \cite{Wei2018, Huang2019,Xu2013}.  Under the assumption that the electric field vanishes and all atoms reside in their ground state as $x\to \pm \infty$,  i.e., $E(x,t)\to 0$ and $(r,s,u)\to (0,0,-1)$, the Cauchy problem for \eqref{eq01} was solved by Gibbon using the inverse scattering transform (IST) \cite{Gibbon1973}. Around the same time, Caudrey, applying the analytical method presented by Wadati \cite{Wadati1972}, demonstrated that the $N$-soliton solution $E(x,t)$ is stable and behave like particles, decomposing into a sum of $N$ single-solitons after collisions. As $t\to \pm\infty$, explicit calculations of the phase shifts confirm that the total phase shift is conserved,  similar to the case of the Korteweg-de Vries (KdV) equation. 
	
	However, the soliton resolution described above  corresponds to the reflectionless case within the IST framework, where the reflection coefficient vanishes. In the general case with nonzero reflection, the long-time asymptotics of solutions to \eqref{eq01} remain open. 
	
	The purpose of this paper is to analyze the long-time behavior of solutions to the Cauchy problem \eqref{eq01} with initial data $E_0(x)\in H^{1,1}(\mathbb{R})$, under the assumption that all atoms initially reside in their ground state. Here, for nonnegative integer $j$ and $k$, the weighted Sobolev space is defined as $$H^{j,k}( \mathbb{R} ) :=\left\{ f\in L^{2}( \mathbb{R} ) :\partial _{x}^{j}f,| x| ^{k}f\in L^{2}( \mathbb{R} ) \right\}.$$ A key development is the nonlinear steepest descent method introduced by Deift and Zhou, which provides a rigorous framework for the asymptotic analysis of oscillatory Riemann–Hilbert problems (RHPs)\cite{MR1207209, MR1989226}. This method has seen considerable refinements and extensions over time, allowing it to address more and more complex Cauchy problems \cite{Boutet2025,Charlier2023,Biondini2021}. In particular, the introduction of $\bar{\partial}$-techniques \cite{McLaughlin2006,McLaughlin2008} by McLaughlin and Miller \cite{McLaughlin2006,McLaughlin2008} has led to a series of important advances in the study of long-time asymptotics for integrable nonlinear wave equations \cite{MR3795020,Cuccagna2016,Dieng2019,Jenkins2018,Yang2022,Chen2021}.
	
	By employing the $\bar{\partial}$ generalization of the nonlinear steepest descent method, we derive the main asymptotic result stated in Theorem \ref{theo1}, which provides the precise long-time expansions of solutions to the RMB equations \eqref{eq01} within the conical region $ C(x_1, x_2, v_1, v_2) \subset \mathbb{R}^2 $. At leading order, the asymptotics are governed by the superposition of single-solitons and  single-kinks. Specifically, the electric field $ E(x,t) $ and the Bloch components $ u(x,t) $ and $ r(x,t) $ are asymptotically separate into sums of $ N(\mathcal{J}) $ single-solitons, while $ s(x,t) $ is described by a sum of $ N(\mathcal{J}) $ single-kinks. The velocities of the fundamental solitons and kinks are determined by the discrete spectrum lying within the spectral interval associated with the cone. Moreover, in the expansion of $E(x,t)$, the subleading $t^{-1/2}$ term characterizes the interaction between solitons and dispersive radiation, reflecting the influence of the continuous spectrum. Consequently, these results confirm that the solutions to \eqref{eq01} obey the soliton resolution conjecture.

	Although the asymptotic expansions derived in  Theorem \ref{theo1} are  the expectable results of the nonlinear steepest descent method, our analysis is highly nontrivial. Since the set of equations \eqref{eq01} corresponds to a negative flow in the Zakharov–Shabat hierarchy,   the $t$-part \eqref{eq03} of the associated Lax pair exhibits two real singularities at $z=\pm \mu/2$. Here, $\mu$ is a finite parameter that, from a physical perspective, plays the role of an effective resonance frequency. Without loss of generality, we assume $0<\mu\leq1$. These singularities at $z=\pm \mu/2$ induce rapid oscillations in the time-evolved  reflection coefficient near the singular points. Therefore, it is inappropriate to directly formulate the inverse problem within the IST framework using such oscillatory scattering data. To obtain a rigorous inverse problem, these singular behaviors should be carefully excluded. Such issues are typical for negative flows of integrable systems but have received little rigorous attention. In the case of \eqref{eq01}, the singularities arise from the sharp-line limit approximation.

    A brief discussion of such singularities was shown by Li and Miller in their analysis of the asymptotic behavior of causal solutions to the Maxwell–Bloch system \cite{Li2024}, where the RHP was formulated with a jump contour avoiding the singular points. In our setting, we similarly construct a properly posed RHP with a jump contour that avoids the singularities. Unlike \cite{Li2024}, this is achieved by modifying the time-evolution formula of the reflection coefficient (see \eqref{eq7.2}) so that it is defined on a $\kappa$-dependent contour $\widetilde\Sigma$, which avoids the singularities and tends to the real line as $\kappa \to 0$. For all $ z \in \mathbb{R} \setminus \{ \pm \mu/2 \} $, the modified time-evolved reflection coefficient converges to the original one as $ \kappa \to 0 $. Moreover, for any fixed $ t > 0 $, it lies in the Sobolev space $ H^{1,1}(\mathbb{R}) $. Building upon this, Theorem \ref{theo2} shows that the properly posed RHP indeed reconstructs the solutions to \eqref{eq01}.

     As discussed above, the electric field $E(x,t)$ may, in general, contain finitely many solitons and breathers. In the present analysis, we restrict our attention to the case where only solitons are present. To this end, we assume that the initial data $E_0(x)$ gives rise, under the direct scattering transform, to scattering data whose discrete spectrum is located entirely on the imaginary line and consists of simple poles, each corresponding to a soliton with a distinct velocity. Given a cone determined by  a fixed spatial interval and a prescribed velocity range, the solitons that appear within the cone correspond precisely to eigenvalues lying within a specified subinterval of the imaginary line.
     
     Moreover, it is well known from the nonlinear steepest descent method \cite{MR3795020,MR1207209} that contributions from jump matrices supported away from the real line are exponentially small in the long-time limit. This explains why the modified time-evolved reflection coefficient does not contribute significantly to the asymptotic expansions.
   
     %在应用Dbar 方法的时候，我们发现它有四个稳态相位点，但是只有两个会产生贡献,因此使得我们的分析类似mkdv
     %同样用PC模型来近似局部的跳跃条件，贡献了描述孤子与辐射相互作用的次领头项。  
	The paper is set out as follows: 
	In Section \ref{S2}, we formulate the inverse problem for the Cauchy problem of the RMB equations \eqref{eq01} as a properly posed RHP \ref{RHP1} within the IST framework. The reconstruction formulae are given in Theorem \ref{theo2}. As preliminary work, we review the direct scattering transform for the Zakharov–Shabat problem and examine the influence of singularities. The section concludes with the statement of the main asymptotic result, Theorem \ref{theo1}. In Section \ref{S3}, the RHP \ref{RHP1} is transformed into a  $\overline{\partial}$-RHP \ref{RHP3} via contour deformations and continuous extensions of jump factorization. Concretely, by introducing suitable transformations inspired by the $ \overline{\partial} $ steepest descent method, the jump contour of RHP \ref{RHP1} is deformed so that  the oscillatory factors $e^{\pm2it\theta(z)}$ decay and the residue conditions remain well-controlled. In Section \ref{S4}, we construct the asymptotic approximations to the solution of the $\overline{\partial}$-RHP \ref{RHP3} by decomposing it into four model problems, each responsible for different contributions to the long-time behavior. The out model $M_{\mathrm {out}}$ yields the leading-order term in the asymptotic expansion \eqref{THE1}; the local models $M_{\pm \zeta_0}$ are constructed to resolve the non-uniformly decaying jump conditions in neighborhoods of the stationary points in the $\overline{\partial}$-RHP \ref{RHP3}; and the error matrix $ E_\circ $ determines subleading contributions arising from soliton–radiation interactions. These components are assembled via formula \eqref{eq73} into the meromorphic part $M^{(2)}_{\mathrm{RHP}}$ of $M^{(2)}$. The non-analytic part is encoded in the solution $M^{(3)}$ to the $\overline{\partial}$-Problem \ref{RHP5}, contributing the final error term in \eqref{THE1}. In Section \ref{S5},  the proof of Theorem \ref{theo1} is completed by summarizing the model solutions.
  
	\section{Riemann-Hilbert problem construction for the reduced Maxwell-Bloch equations} \label{S2}
	This section provides a review of the direct scattering transform for the Zakharov-Shabat scattering problem. We then analyze the issues arising from singularities emerging in the $t$-part of the Lax pair \eqref{eq03}. To avoid the complexities, we formulate a properly posed RHP with a carefully designed jump contour that bypasses the singularities, effectively mitigating their adverse effects.  
	\subsection{Direct scattering transform}
	The  RMB equations \eqref{eq01} serve as the compatibility conditions, stating that for the spectral parameter $z\in \mathbb{C}\backslash \{\pm \mu/2\}$, there exists a $2\times 2$ matrix-valued function $\Phi(x,t;z)$ that simultaneously satisfies the following linear equations 
	\begin{align}
		&\begin{aligned}\label{eq02}
			\Phi _{x}=( -iz\sigma _{3}+Q_{1}) \Phi
		\end{aligned}\\ 
		&\begin{aligned}\label{eq03}
			\Phi _{t}=\frac{1}{4z^{2}-\mu ^{2}}( iz( \sigma _{3}+Q_{2}) +Q_{3}) \Phi,
		\end{aligned}
	\end{align}
	where
	\begin{align*}%\label{eq04}
		\begin{aligned}
			Q_{1}( x,t) :=\frac{1}{2}
			\begin{pmatrix}
				0 & -E \\
				E & 0
			\end{pmatrix},
			Q_{2}( x,t) :=
			\begin{pmatrix}
				-1-u & -s \\
				-s & 1+u
			\end{pmatrix},
			Q_{3}(x,t) :=\frac{\mu }{2}
			\begin{pmatrix}
				o & r \\
				-r  & 0
			\end{pmatrix},
		\end{aligned}
	\end{align*}
	and the standard Pauli matrices are denoted by
		\begin{align*}
			\begin{aligned}
				\sigma _{1}=\begin{pmatrix}
					0 & 1 \\
					1 & 0
				\end{pmatrix},
				\sigma _{2}=\begin{pmatrix}
					0 & -i \\
					i & 0
				\end{pmatrix},
				\sigma _{3}=\begin{pmatrix}
					1 & 0 \\
					0 & -1
				\end{pmatrix}.
			\end{aligned}
	\end{align*}
	The equation \eqref{eq02}, as the classical Zakharov-Shabat scattering problem (also known as the AKNS system), was first introduced in the study of the nonlinear Schr\"odinger equation (NLS) and rigorously analyzed by Beals and Coifman and Zhou \cite{Zakharov,beals,zhou}. Within the framework of IST, to facilitate a clearer understanding of the inverse problem established in Subsection \ref{s2.2}, we revisit the well-known direct scattering transform of the Zakharov-Shabat scattering problem. Readers already familiar with standard IST may proceed directly to Assumption \ref{Ass1}.
	
	The direct scattering transform refers to a Lipschitz continuous mapping from the initial data $E_0(x)$ to the associated scattering data. For $z\in \mathbb{R}$, let $\phi_\pm(x,z)$ be the two Jost solutions of equation \eqref{eq02} with the asymptotic behavior $e^{-izx\sigma_{3}}$ as $x\to \pm\infty$.  Then, define the normalized functions $m_\pm(x,z):=\phi_\pm(x,z)e^{izx\sigma_{3}}$, which satisfy the following Volterra integral equation 
	\begin{align}\label{eq05}
		\begin{aligned}
			m_{\pm }( x,z) =I+\int _{\pm \infty }^{x}e^{iz( y-x) \widehat{\sigma }_{3}}\left(  Q_{1}(y)m_\pm ( y,z) \right) dy,
		\end{aligned}
	\end{align}
	where the initial data $E_0(x)$ serves as the entry of $Q_1(x)$. Clearly, for  $z\in \mathbb{R}$, there exist uniquely continuous solutions $m_\pm(\cdot,z)\in  L^{\infty}( \mathbb{R} )$ satisfying the Volterra integral equations \eqref{eq05} when $E_{0}(x) \in L^{1}( \mathbb{R} )$. Thus, the Jost solutions $\phi_\pm(x,z)$ are bounded. Furthermore, $\phi_\pm(x,z)$ are linearly dependent, as they are solutions to a first-order differential equation with two components \eqref{eq02}. There exists a scattering matrix $S(z)$, independent of $x$, such that
	\begin{align}\label{eq1}
		\begin{aligned}
			\phi_-(x,z) =\phi _{+}( x,z) S( z) ,~~S(z) :=\begin{pmatrix}
				s_{11}( z)  & s_{12}(z) \\
				s_{21}( z)  & s_{22}( z) 
			\end{pmatrix},z\in \mathbb{R} .
		\end{aligned}
	\end{align}
	Since $\phi_\pm(x,z)$ have the unit determinant by Abel's theorem, the scattering coefficients $s_{11}(z)$ and $s_{21}(z)$ can be expressed as the Wronskian of the normalized functions
	\begin{align}\label{eq2}
		\begin{aligned}
			s_{11}(z) = Wr ( m_{-,1}(x,z) ,m_{+,2}( x,z) ) ,~~s_{21}(z) =-e^{-2izx} Wr( m_{-,1}(x,z) ,m_{+,1}( x,z) ) ,
		\end{aligned}
	\end{align}
	where $$m_\pm(x,z):=\left(m_{\pm,1}(x,z),m_{\pm,2}(x,z)\right).$$ 
	These coefficients are continuous for $z\in \mathbb{R}$.
	Additionally, the relationships
	\begin{align*}%\label{eq3}
		\begin{aligned}
			s_{22}(z) =\overline{s_{11}}(z) ,~~s_{12}(z)=-\overline{s_{21}}(z) ,
		\end{aligned}
	\end{align*}
	follow easily from the symmetries of the Jost solutions $\phi_\pm(x,z)=\sigma_2\overline{\phi_\pm}(x,z)\sigma_2$.
	Considering the first column of equation \eqref{eq1}, the reflection and transmission coefficients are defined by
	\begin{align*}%\label{eq4}
		\begin{aligned}
			r(z):=\frac{s_{21}(z)}{s_{11}(z)}, ~~\tau(z):= \frac{1}{s_{11}(z)},
		\end{aligned}
	\end{align*}
	which are crucial for our analysis.
	
	Furthermore, the Neumann series converges absolutely and uniformly to the respective columns of the normalized functions $m_\pm(x,z)$, which can thus be analytically extended to the half-planes $\mathbb{C}^\pm := \{ z \in \mathbb{C} : \pm \operatorname{Im} z > 0 \}$. 
	\begin{prop}\label{prop1}
		For $x\in \mathbb{R}$, the normalized functions $m_{-,1}(x,z)$ and $ m_{+,2}(x,z)$ can be analytically extended  to $\mathbb{C^+}$, while $m_{+,1}(x,z)$ and $ m_{-,2}(x,z)$ can be analytically extended to $\mathbb{C^-}$.
	\end{prop}
	As a result, the coefficient $s_{11}(z)$ admits analytic extension in $\mathbb{C^+}$ and is continuous for $\IIm z \geq 0$, due to its Wronskian representation in \eqref{eq2}. Using the symmetries $\phi_\pm(x,z)=\sigma_2\overline{\phi_\pm}(x,\overline z)\sigma_2$, we can obtain
	$s_{22}(z) =\overline{s_{11}}(\overline z) $, which implies that $s_{22}(z)$ admits analytic extension in $\mathbb{C^-}$  and is continuous for $\IIm z\leq 0$. However, the coefficients $s_{21}(z)$ remains continuous only for $z\in \mathbb{R}$, so the reflection coefficient $r(z)$ is continuous solely on $z\in \mathbb{R}$.
	
	\begin{prop}\label{prop2}
		For $x\in \mathbb{R}$, each column of the normalized functions $ m_\pm(x,z)$ exhibits asymptotics
		\begin{align*}%\label{eq5}
			\begin{aligned}
				\left(m_{-,1}(x,z),m_{+,2}(x,z)\right)\to I, ~~\left(m_{+,1}(x,z),m_{-,2}(x,z)\right)\to I,
			\end{aligned}
		\end{align*}	
		as $z\to \infty$ along a contour in the domains of their analyticity, extended with $\IIm z\neq 0$.
	\end{prop}
	This result is established via a direct computation of the limits of the corresponding elements in equation \eqref{eq05}. Then, taking the limit $z\to \infty$ with $\IIm z> 0$ in both sides of the first equation in \eqref{eq2}, we derive the asymptotics $s_{11}(z)\to 1$. Similarly, by the symmetry, $s_{22}(z)\to 1$ as $z\to \infty$ with $\IIm z< 0$.
	Before continuing, we make the following assumption on the initial date for $E_0(x)$ and the scattering coefficient $s_{11}(z)$ as the basis for subsequent discussion.
	\begin{Assu}\label{Ass1}
		The initial data $E_0(x):\mathbb{R}\mapsto \mathbb{R}$ belongs to the weighted Sobolev space $H^{1,1}( \mathbb{R} )$.
		Moreover, the zeros of the coefficients $s_{11}(z)$, arising from the direct scattering transform associated with $E_0(x)$, satisfy the following statements.
		\begin{itemize}
			\item[1.]  For $z\in \mathbb{R}$, there are no zeros of $s_{11}(z)$ implying the absence of spectral singularities.
			\item[2.] For $\IIm z> 0$, there are finitely many simple zeros of $s_{11}(z)$. 
		\end{itemize}
		Specifically, let $\mathcal{Z}:=\{z_k\in i\mathbb{R^+}\}_{k=1}^N$, with $z_j\neq z_k$ for $j\neq k$, be the set of zeros of $s_{11}(z)$. Then, there exists a corresponding set $C:=\{c_k\in i\mathbb R\backslash\{0\}\}_{k=1}^N$ such that the identity holds
		\begin{align*}%\label{eq6}
			\begin{aligned}
				\phi_{-,1}(x,z_k)=c_k\phi_{+,2}(x,z_k),
			\end{aligned}
		\end{align*}
		where $\phi_\pm(x,z):=\left(\phi_{\pm,1}(x,z),\phi_{\pm,2}(x,z)\right)$.
	\end{Assu}
	\begin{Remark}
		The symmetries of Jost solutions govern the distribution of the zeros of $s_{11}(z)$ and  $s_{22}(z)$.
		\begin{itemize}
			\item From $\phi_\pm(x,z)=\sigma_2\overline{\phi_\pm}(x,\overline{z})\sigma_2$, each purely imaginary zero $z_k \in i\mathbb{R}^+$ of $s_{11}(z)$ corresponds to an soliton in $E(x,t)$, with $\overline{z}_k$ being zero of $s_{22}(z)$.
			\item The additional symmetry $\phi_\pm(x,z)=\overline{\phi_\pm}(x,-\overline{z})$ implies that any non-imaginary zero $z_k \in \mathbb{C}^+$ must be accompanied by a symmetric zero at $- \overline z_k$ in $s_{11}(z)$. Each pair of  anti-hermitian conjugate eigenvalues corresponds to a breather in $E(x,t)$, analogous to the case of the real-valued mKdV equation \cite{Wadati19731}.
		\end{itemize}
		For simplicity, we restrict our analysis to zeros of $s_{11}(z)$ on $i\mathbb{R^+}$.
	\end{Remark}
	
	Thus, the direct scattering transform maps the initial data $E_0(x)$ to the scattering data, consisting of the reflection coefficient $r(z)$, the discrete spectrum $\mathcal{Z}$, and the corresponding norming constants $C$. According to results in \cite{zhou},  given initial data $E_0(x)\in H^{j,k}(\mathbb R)$, it is known that $s_{11}(z)-1\in H^{k,1}(\mathbb R)$ and $s_{21}(z)\in H^{k,j}(\mathbb R)$. Furthermore, from the first part of Assumption \ref{Ass1}, the reflection coefficient $r(z)$ belongs to the Sobolev space $H^{1,1}(\mathbb R)$.
	
	\subsection{Inverse problem formulated as a Riemann-Hilbert problem}\label{s2.2}
	The  inverse problem in IST refers to reconstruct the solutions of Cauchy problem \eqref{eq01} from the scattering data evolved over time $t$. Decades ago, the inverse problem was commonly formulated as an algebraic integral equation, known as the Gel'fand-Levitan-Marchenko (GLM) equation, via a Cauchy integral on the real line. A more modern and streamlined approach reformulates the inverse problem as a RHP, which seeks a sectionally meromorphic function in $\mathbb{C}$ satisfying specified jump conditions along a contour (typically the real line) and residue conditions at discrete poles $z_k\in \mathcal{Z}\cup \overline {\mathcal{Z}}$. 
	
	The temporal evolution of the scattering data is dictated by the $t$-part of the Lax pair \eqref{eq03}. Under the boundary condition $(r, s, u) \to (0, 0, -1)$ as $x \to \pm\infty$, the time-evolved reflection coefficient and norming constants are explicitly expressed as
	\begin{align*}%\label{eq7}
		\begin{aligned}
			r (z,t) := r(z) e^{-\frac{2iz}{4z^{2}-\mu ^{2}}t},~~c_{k}(t) =c_{k}e^{-\frac{2iz_k}{4z_k^{2}-\mu ^{2}}t}, ~~k=1,\ldots ,N.
		\end{aligned}
	\end{align*}
	Moreover, the coefficient $s_{11}(z,t)=s_{11}(z)$ remains time-invariant during time evolution, implying that the discrete spectrum $\mathcal Z$ is time-independent. 
	
	Note that, for $0<\mu\leq1$, the coefficient in equation \eqref{eq03} exhibits two first-order singularities at $z = \pm \mu/2$, which induce rapid oscillations in the evolution factor $\exp(-\frac{2iz}{4z^{2}-\mu^2}t)$ near these points. This gives rise to several analytical difficulties in the formulation of the inverse problem.
	\begin{itemize}
		\item 
		The real line is partitioned into three disjoint intervals by the singular points $z=\pm \mu/2$, which causes a loss of information at  $r(\pm\mu/2)$ when formulating the jump conditions for the RHP and the Cauchy integral representation associated with the GLM equation. 
		\item Direct computation of the derivative of 
		$$\frac{\partial }{\partial z}r ( z,t) =\left( r '(z) +\frac{2it( 4z^{2}+\mu ^{2}) }{( 4z^{2}-\mu^{2}) ^{2}}r( z) \right) e^{-\frac{2iz}{4z^{2}-\mu^{2}}t},$$ 
		demonstrates  $\partial_z r(z, t) \notin L^2 (\mathbb{R})$ for any $t>0$.
	\end{itemize}
	As a consequence, to rigorously formulate the inverse problem for $t \neq 0$, the singular behaviors should be carefully excluded. To this end, we define a modified time-evolved reflection coefficient by
	\begin{align}\label{eq7.2}
		\begin{aligned}
			\widetilde{r}( z,t) =r( \RRe z) e^{-\frac{2iz}{4z^{2}-\mu^2}t},~~z\in \widetilde\Sigma,
		\end{aligned}
	\end{align}
	where for fixed $\kappa>0$
		\begin{align*}%\label{eq7}
		\begin{aligned}
			\widetilde\Sigma:=\Sigma_{\mathrm {re}} \cup \Sigma _{1}^{\kappa}\cup \Sigma _{4}^{\kappa},&~~~\Sigma _{\mathrm {re}}:=\{ z\in \mathbb{R}:|z\pm {\mu }/{2}| \geq\kappa \} \\
				\Sigma _{1}^{\kappa}=\{ z\in \mathbb{C} :| z-{\mu }/{2}| =\kappa,\IIm z >0\} ,&~~~\Sigma _{4}^{\kappa}=\{ z\in \mathbb{C} :| z+{\mu }/{2}| =\kappa,\IIm z >0\} .\end{aligned}
	\end{align*}
	Moreover, for all $z\in \mathbb{R}\backslash \{\pm \mu/2\}$, we have
	\begin{align*}%\label{eq7}
		\lim _{\kappa\to 0}\widetilde{r}( z,t)={r}( z,t).
	\end{align*}
	That is, although $\widetilde{r}( z,t)$ is defined on a $\kappa$-dependent contour $ \widetilde\Sigma$ off the real line for each fixed $\kappa>0$, it converges to ${r}( z,t)$ as $\kappa \to 0$.
	
	\begin{prop}\label{prop0}
	Given the reflection coefficient $r(z)\in H^{1,1}(\mathbb R)$ and parameters $0<\kappa_0\leq\kappa <\mu \leq 1$, 
	there exists a constant $c(T, \kappa_0) > 0$ such that for any $t \in (0, T]$ with $T > 0$,
	\begin{align}\label{eq7.0}
		\|\widetilde r(z,t) \|_{H^{1,1}(\widetilde \Sigma)}\leq c(T,\kappa_0)\|r(z)\|_{H^{1,1}(\mathbb R)}.
	\end{align}
	\end{prop}
	\begin{proof}
		A direct computation gives
		 \begin{align}\label{eq7.1} 
		 	\RRe [{2iz}/(4z^{2}-\mu^2)]=\IIm z \widetilde G(\RRe z, \IIm z),
		 	\end{align}
		where  
		\begin{align*}
			\widetilde G(a, b):=\frac{2( \mu ^{2}+4a^{2}+4b^{2}) }{( ( 2a+\mu ) ^{2}+4b^{2}) ( ( 2a-\mu) ^{2}+4b^{2}) }.
		\end{align*}
		For $z \in \Sigma_1^\kappa$ and $\varepsilon > 0$, let $z = \mu/2 + \kappa e^{i\alpha}$ with $\varepsilon < \alpha < \pi - \varepsilon$. Then from \eqref{eq7.1}, we obtain
		\begin{align*}
			\RRe [{2iz}/(4z^{2}-\mu^2)]=\frac{\sin \alpha }{4\kappa}\left( 1+\frac{\kappa^{2}}{\kappa^{2}+2\mu \kappa\cos \alpha +\mu^{2}}\right) \geq \frac{\sin\varepsilon}{4\kappa}.
		\end{align*}
		A same bound can be derived for $z \in \Sigma_4^\kappa$. Using this and the assumption $r(z) \in H^{1,1}(\mathbb{R})$,  we have
		\begin{align*}
			\| \widetilde r( z,t) \|^2 _{L_{z}^{2}( \widetilde\Sigma) }&= \int _{\Sigma_{ \mathrm {re}}}| r( s) | ^{2}ds+\int _{\Sigma^\kappa_{ 1}+\Sigma^\kappa_{ 4}} |r ( \RRe s)|^2 e^{-t\RRe[4is/({4s^{2}-\mu ^2})]}|ds|\\
			&\leq \|r(z)\|^2_{L^2(\Sigma_{ \mathrm {re}})}+2\kappa e^{-t{\sin \varepsilon}/({2\kappa})}\int_\varepsilon^{\pi-\varepsilon}|r(\RRe s)|^2d\alpha \leq \infty.
		\end{align*}
		Similarly, $\| z\widetilde{r}(z, t) \|^2_{L^2_z(\widetilde{\Sigma})} < \infty$. 
		
		Now consider the derivative with respect to $z$
		\begin{align*}%\label{eq7}
			\begin{aligned}\frac{\partial }{\partial z}\widetilde{r}( z,t)=
				\left( \frac{1}{2}r '(\RRe z) +\frac{2it( 4z^{2}+\mu ^{2}) }{( 4z^{2}-\mu ^{2}) ^{2}}r( \RRe z) \right) e^{-\frac{2iz}{4z^{2}-\mu ^{2}}t}.
				\end{aligned} 
		\end{align*}
		Since $r(z) \in H^{1,1}(\mathbb{R})$, the first term is square-integrable on $\widetilde{\Sigma}$. For the second term, take $s = \pm \mu/2 + \kappa e^{i\alpha} \in \Sigma_1^\kappa \cup \Sigma_4^\kappa$ with $\kappa < \mu$. Then
	\begin{align*}%\label{eq7}
		\begin{aligned}
			\left| \frac{2it( 4s^{2}+\mu ^{2}) }{( 4s^{2}-\mu ^{2}) ^{2}}\right| \leq \frac{ 2t( 2\mu ^{2}+ 4\mu\kappa  +4\kappa^{2}) }{16\kappa^{2}| \pm \mu +\kappa e^{i\alpha }| ^{2}}\leq \frac{t( ( \mu + \kappa) ^{2}+\kappa^{2}) }{4\kappa^{2}| \mu -\kappa| ^{2}}.
		\end{aligned} 
	\end{align*}
	Thus, for some constant $c> 0$, we obtain
	\begin{align*}%\label{eq7}
		\begin{aligned}
		&	\int_{\widetilde{\Sigma}}\left|\frac{2it( 4s^{2}+\mu ^{2}) }{( 4s^{2}-\mu ^{2}) ^{2}}r( \RRe s) \right|^2e^{-t\RRe[4is/({4s^{2}-\mu ^2})]}|ds|\\
		&	\leq \frac{ct^2}{\kappa^4}  \|r(z)\|^2_{L^2(\Sigma_{ \mathrm {re}})}+	\int_{\Sigma^\kappa_{ 1}+\Sigma^\kappa_{ 4}}\left|\frac{2it( 4s^{2}+\mu ^{2}) }{( 4s^{2}-\mu ^{2}) ^{2}}r( \RRe s) \right|^2e^{-t\RRe[4is/({4s^{2}-\mu ^2})]}|ds|\\
		&\leq   \frac{ct^2}{\kappa^3} \left(\frac{1}{\kappa} \|r(z)\|^2_{L^2(\Sigma_{ \mathrm {re}})}+e^{-t{\sin \varepsilon}/({2\kappa})}\int_\varepsilon^{\pi-\varepsilon}|r(\RRe s)|^2d\alpha\right).
		\end{aligned} 
	\end{align*}
	Thus, $ \partial_z \widetilde{r}(z,t) \in L^2(\widetilde{\Sigma}) $. Combining this with previous estimates yields $ \widetilde{r}(z,t) \in H^{1,1}(\widetilde{\Sigma}) $. 
 	\end{proof}
 	
 	\begin{Remark}
 		Bullough \cite{Bullough1979}, using the Sokhotski–Plemelj formula in the distributional sense,
 		$$(z \pm \mu/2)^{-1} = \mathcal{P}(1/(z \pm \mu/2)) \mp i\pi \delta(z\pm \mu/2) ,$$
 		showed that $r(\pm \mu/2,t)=0$, where $\mathcal{P}$ denotes the Cauchy principal value and $\delta$ the Dirac delta function.
 		
 		A similar vanishing condition holds for the modified reflection coefficient $\widetilde{r}( z,t)$ defined in this work. Specifically, for $z\in \Sigma_1^\kappa\cup\Sigma_4^\kappa $ approaching $ \pm \mu/2$ as $\kappa\to 0$,  we have
 			\begin{align*}%\label{eq7}
 			|\widetilde{r}( z,t) |=|r( \RRe z)|e^{	-\RRe [{2iz}/(4z^{2}-\mu^2)] t}=|r( \RRe z)|  \exp \left(\frac{\sin \alpha }{4\kappa}\left( 1+\frac{\kappa^{2}}{\kappa^{2}\pm2\mu \kappa\cos \alpha +\mu^{2}}\right)\right) \to 0,
 		\end{align*}
 		Thus, the vanishing of \( \widetilde{r}(z, t) \) at \( z = \pm \mu/2 \) follows directly from the exponential decay, without requiring a distributional interpretation.
 		
 		Moreover, since the constant $c(T,\kappa_0)$ in estimate \eqref{eq7.0} grows at most polynomially in $t$, the time-evolved of $\widetilde{r}(z,t)$ remains controlled by the initial reflection coefficient $r(z)$. 
 		This control of time-evolved scattering data mirrors a key step in Pelinovsky's derivation \cite{Pelinovsky2017} of an a priori Lipschitz estimate for the derivative nonlinear Schrödinger equation, where bounding the solution $u(t,x)$ by its initial data $u_0(x)$ over $t \in (0,T]$ relied on such control.
 		Therefore, formulating the inverse problem using $\widetilde{r}(z,t)$ effectively ensures the global existence of solutions to the RMB equations.
 	%	In the analysis of the global existence of solutions to the derivative nonlinear Schr\"odinger equation, Pelinovsky \cite{Pelinovsky2017} derived an a priori Lipschitz estimate that bounds the solution $u(t,x)$ in terms of the initial data $u_0(x)$ for all $t \in (0,T]$, based on a similar control of the time-evolved scattering data. 
 	\end{Remark}
	
	From Proposition \ref{prop0},  the inverse problem can be formulated as a properly posed RHP whose jump contour avoids the singular points $ z = \pm \mu/2 $, thereby effectively circumventing the complications induced by these singularities.
	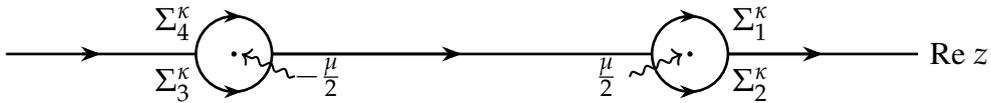
\begin{figure}[H] 
		 \centering
		 \begin{tikzpicture}[scale=1, transform shape]
			\draw[black,  line width=1pt] (-6,0) -- (-3.5,0);
			\draw[black, line width=1pt, -{Stealth}] (-5,0) -- (-4.75,0);
			\draw[black, line width=1pt] (-2.5,0) -- (2.5,0);
			\draw[black, line width=1pt, -{Stealth}] (-2.5,0) -- (0,0);
			\draw[black, line width=1pt] (3.5,0) -- (6,0) node[right] {\text{Re} $z$};
			\draw[black, line width=1pt, -{Stealth}] (3.5,0) -- (4.75,0);
			
			%left circle
			\draw[line width=1pt] (-3.5,0) arc (180:0:0.5);
			\draw[black, thick, -{Stealth}] (-3,0.5) -- (-2.9,0.5);
			\node[black] at (-3.8,0.4) {$\Sigma_4^\kappa$};
			\draw[line width=1pt] (-3.5,0) arc (-180:0:0.5);
			\draw[black, thick, -{Stealth}] (-3,-0.5) -- (-2.9,-0.5);
			\node[black] at (-3.8,-0.4) {$\Sigma_3^\kappa$};
			\fill (-3,0) circle (1pt);
			\node[black] at (-1.9,-0.3) {\(-\frac{\mu}{2}\)};
			\draw[->, thick, decorate, decoration={snake, amplitude=.4mm, segment length=2mm, post length=1mm}] (-2.2,-0.3) -- (-2.9,0);
			%right circle
			\draw[line width=1pt] (2.5,0) arc (180:0:0.5);
			\draw[black, thick, -{Stealth}] (3,0.5) -- (3.1,0.5);
			\node[black] at (3.8,0.4) {$\Sigma_1^\kappa$};
			\draw[line width=1pt] (2.5,0) arc (-180:0:0.5);
			\draw[black, thick, -{Stealth}] (3,-0.5) -- (3.1,-0.5);
			\node[black] at (3.8,-0.4) {$\Sigma_2^\kappa$};
			\draw[->, thick, decorate, decoration={snake, amplitude=.4mm, segment length=2mm, post length=1mm}] (2.2,-0.3) -- (2.9,0);
			\fill (3,0) circle (1pt);
			\node[black] at (1.9,-0.3) {\(\frac{\mu}{2}\)};
		\end{tikzpicture}
		\caption{\small {The jump contour $\Sigma_{1}$ of RHP \ref{RHP1} consists of four semicircular arcs $\Sigma_k^\kappa$ ($k=1,\dots,4$) and three intervals: $(-\infty, -\mu/2 - \kappa]$, $[- \mu/2 + \kappa, \mu/2 - \kappa]$, and $[\mu/2 + \kappa, \infty)$. Each $\Sigma_k^\kappa$ corresponds to the upper or lower semicircle of the circles $|z \pm \mu/2| = \kappa$.} }\label{F1}
	\end{figure}
	\begin{rhp}\label{RHP1} 
		Let $\kappa>0$ be sufficiently small and consider the contour $\Sigma_1$ shown in Figure \ref{F1}. For given scattering data denoted by $\sigma:=\{r(z),\{(z_k,c_k)\}_{k=1}^N\}$, seek for a meromorphic matrix-valued function $M(z)=M(z;x,t)$ defined on  $z\in \mathbb{C}\backslash\{\Sigma_1\}$ that satisfies the following conditions. 
		\begin{itemize}
			\item $M(z)$ has continuous boundary values $M_\pm(z):=\lim\limits_{\varepsilon\to 0^+}M(z\pm i\varepsilon)$  on $z\in \Sigma_1$, which satisfy the jump conditions
			\begin{align}\label{eq8}
				\begin{aligned}
					M_+(z)=M_-(z)V(z;x,t),
				\end{aligned}
			\end{align}
			where 
			\begin{align}\label{eq9}
				\begin{aligned} 
					V(z;x,t) :=\begin{cases}\begin{pmatrix}
							1+| r(z) | ^{2} & \overline{r}(z) e^{2i\varphi(z)} \\
							r( z) e^{-2i\varphi(z)}  & 1
						\end{pmatrix}&z\in \Sigma_{\mathrm {re}}\\
						\begin{pmatrix}
							1 & 0 \\
							r( \RRe z) e^{-2i\varphi ( z) } & 1
						\end{pmatrix}&z\in\Sigma^\kappa_1\cup\Sigma^\kappa_4\\
						\begin{pmatrix}
							1 & \overline{r }( \RRe z)e^{2i\varphi ( z) } \\
							0 & 1
						\end{pmatrix}&z\in\Sigma^\kappa_2\cup\Sigma^\kappa_3.
					\end{cases}	
				\end{aligned}		
			\end{align}
			and $\varphi(z)=\varphi(z;x,t):=\frac{z}{4z^2-\mu^2}t-zx$.
			
			Moreover, $M(z)$ satisfies the following residue conditions  at each simple poles $z_k$ and $\overline z_k$
			\begin{align}\label{eq10}
				\begin{aligned}
					\Res\limits_{z=z_{k}} M(z) =\lim_{z\rightarrow z_{k}} M( z) 
					\begin{pmatrix}
						0 &   0 \\
						c_{k}e^{-2i\varphi({z}_{k}) }  & 0
					\end{pmatrix},
					\Res\limits_{z=\overline{z}_{k}} M(z) =\lim_{z\rightarrow \overline{z}_{k}} M( z)
					\begin{pmatrix}
						0 & -\overline{c}_{k}e^{2i\varphi (\overline{z}_{k}) } \\
						0 & 0
					\end{pmatrix}.
				\end{aligned}
			\end{align}
			
			\item The matrix  $M(z)$ possesses the asymptotic behavior $M(z)\to I$, as $|z|\to \infty$,  where $I$ denotes the identity matrix.
		\end{itemize}
	\end{rhp}
	
	\begin{theo}\label{theo2}
		Let $E_0$ be an initial data satisfying Assumption \ref{Ass1}. Then, the unique solutions for the Cauchy problem \eqref{eq01} can be reconstructed from the solution of RHP \ref{RHP1} using the formulae
		\begin{align}\label{eq11}
			\begin{aligned} 
				&E(x,t) =-4i\lim _{z\rightarrow \infty }\left( zM ( z;x,t) \right)_{12},\quad&s(x,t) =-\frac{1}{2}\left(\rho( \pm\frac{\mu }{2};x,t) _{12}+\rho( \pm\frac{\mu }{2};x,t) _{21}\right),\\
				&u(x,t) =-\rho( \pm\frac{\mu }{2};x,t) _{11},\quad&r(x,t) =\mp \frac{1}{2i}\left( \rho( \pm\frac{\mu }{2};x,t) _{12}-\rho( \pm\frac{\mu }{2};x,t)_{21}\right), 
			\end{aligned}
		\end{align}
		whare $\rho( z;x,t) =M( z;x,t) \sigma _{3}M^{-1}( z;x,t) $.
	\end{theo}
	\begin{proof}
		Let $M (z;x,t)$ be the solution to RHP \ref{RHP1}. We define a matrix function 
		\begin{align}\label{eq12}
			\begin{aligned} 
				\widehat{\phi }( x,t;z) =M( z;x,t) e^{i\varphi ( z) \sigma _{3}},
			\end{aligned}
		\end{align}
		that satisfies the following properties:	
		\begin{itemize}
			
			\item The function $\widehat{\phi }( x,t;z) $ is meromorphic for $z\in \mathbb{C}\backslash\{\Sigma_1\cup\{\pm\mu/2\}\}$. From \eqref{eq8}, its boundary values satisfy the jump conditions 
			\begin{align*}%\label{eq13}
				\begin{aligned} 
					\widehat{\phi }_+( x,t;z) =\widehat{\phi }_-( x,t;z)  \widehat{V }( z) ,~~z\in \Sigma_1,
				\end{aligned}
			\end{align*}
			where $\widehat{V }( z)=e^{-i\varphi ( z) \widehat{\sigma _{3}}}V( x,t;z)$ is independent of $x$ and $ t$. The points $z_k\in \mathcal{Z}\cup\overline{\mathcal{Z}}$ are also the simple poles of the function $\widehat{\phi }( x,t;z)$, which satisfies the residue conditions
			\begin{align*}%\label{eq14}
				\begin{aligned}
					\Res\limits_{z=z_{k}} \widehat{\phi }(z) =\lim_{z\rightarrow z_{k}} \widehat{\phi }(z) 
					\begin{pmatrix}
						0 &   0 \\
						c_{k}  & 0
					\end{pmatrix},~~
					\Res\limits_{z=\overline{z}_{k}} \widehat{\phi }(z) =\lim_{z\rightarrow \overline{z}_{k}} \widehat{\phi }(z)
					\begin{pmatrix}
						0 & -\overline{c}_{k} \\
						0 & 0
					\end{pmatrix}.
				\end{aligned}
			\end{align*}
			Meanwhile, the factor $e^{i\varphi ( z) \sigma _{3}}$ introduces two essential singularities at $z=\pm\mu/2$ into the function $\widehat{\phi }( x,t;z) $.
			
			\item The function $\widehat{\phi }_x( x,t;z) \widehat{\phi }^{-1}( x,t;z)$ is analytic across the entire complex $z$-plane for the following reasons:
			
			The essential singularities of $\widehat{\phi }( x,t;z)$ at $z=\pm\mu/2$ are removable in $\widehat{\phi }( x,t;z)e^{-i\varphi(z;0,t)\sigma_{3}}$. Consequently, as shown by the relation
			\begin{align*}%\label{eq15}
				\begin{aligned}
					\widehat{\phi }_x( x,t;z) \widehat{\phi }^{-1}( x,t;z)=\left(\widehat{\phi }( x,t;z)e^{-i\varphi(z;0,t)\sigma_{3}}\right)_x\left(\widehat{\phi }( x,t;z)e^{-i\varphi(z;0,t)\sigma_{3}}\right)^{-1},
				\end{aligned}
			\end{align*}
			the function $\widehat{\phi}_x(x, t; z) \widehat{\phi}^{-1}(x, t; z)$ is analytic at $z = \pm \mu/2$. Each point in the set $ \mathcal{Z}\cup\overline{\mathcal{Z}}$ is also a removable pole of the function $ \widehat{\psi }(x,t;z) :=\widehat{\phi }( x,t;z) N( z)$, where $ N( z)$ is defined piecewise as
			\begin{align}\label{eq16}
				\begin{aligned}
					N( z) :=\begin{cases}\begin{pmatrix}
							1 & 0 \\
							-\frac{c_{k}}{z-z_{k}} & 1
						\end{pmatrix}&z\in  D( z_{k},\kappa ) ,~~z_{k}\in \mathcal{Z}\\
						\begin{pmatrix}
							1 & \frac{\overline{c}_{k}}{z-\overline{z}_{k}} \\
							0 & 1
						\end{pmatrix}&z\in D( \overline{z}_{k},\kappa ) ,~~ \overline{z}_k\in \overline{\mathcal{Z}}\\
						I&otherwise,\end{cases}
				\end{aligned}
			\end{align}
			where $D( a,\kappa ):=\{z\in\mathbb{C}:|z-a|<\kappa\}$ denotes an open disk centered at $z=a$ with radius $\kappa>0$, and $\partial D(a,\kappa )$ represents its boundary, oriented clockwise. The radius $\kappa$ should be small enough to prevent the disks' boundaries from intersecting with each other or with $\Sigma_1$. As a consequence, the function $\widehat{\psi }(x,t;z)$ is subjected to the following jump conditions on the boundary of each disk
			\begin{align*}%\label{eq17}
				\begin{aligned}
					\widehat{\psi }_+(x,t;z)=\widehat{\psi }_-(x,t;z)N^{-1}(z), ~~z\in \partial D( z_{k},\kappa ),~~z_k\in\mathcal{Z}\cup\overline{\mathcal{Z}},
				\end{aligned}
			\end{align*}
			in which the jump matrix $N^{-1}(z)$ is also independent of $x$ and $ t$. Meanwile, on $z\in  \Sigma_1$, $\widehat{\psi }(x,t;z)$ satisfies the same jump conditions as $\widehat{\phi }(x,t;z)$. It can be shown that the function $\widehat{\psi }_x( x,t;z) \widehat{\psi }^{-1}( x,t;z)$ is continuous across the boundaries $z\in \partial D( z_{j},\kappa )$ with $z_k\in \mathcal{Z}\cup\overline{\mathcal{Z}}$, and $z\in \Sigma_1$ (i.e.,  jump-free). Therefore, using the equivalence
			\begin{align}\label{eq18}
				\begin{aligned}
					\widehat{\phi }_x( x,t;z) \widehat{\phi }^{-1}( x,t;z)=\widehat{\psi }_x( x,t;z) \widehat{\psi }^{-1}( x,t;z),
				\end{aligned}
			\end{align}
			we conclude that $\widehat{\phi }_x( x,t;z) \widehat{\phi }^{-1}( x,t;z)$ is jump-free on $z\in \Sigma_1$, and that each pole $z_k\in \mathcal{Z}\cup\overline{\mathcal{Z}}$ is removable.
			
			Furthermore, by combining the definition \eqref{eq12} with the asymptotics $M(z)\to I+M_1(x,t)/z+\mathcal{O}(z^{-1})$, as $|z|\to \infty$, we calculate the asymptotics
			\begin{align*}%\label{eq19}
				\begin{aligned}
					\widehat{\phi }_x( x,t;z) \widehat{\phi }^{-1}( x,t;z)=-iz\sigma _{3}-i \left[ M_1,\sigma_{3} \right]+\mathcal{O}(z^{-1}) .
				\end{aligned}
			\end{align*}
			By Liouville's theorem, it can be shown that the function $\widehat{\phi}_x(x, t; z) \widehat{\phi}^{-1}(x, t; z)$ is a linear function of the form $-iz\sigma_3 - i[M_1, \sigma_3]$. Consequently, $\widehat{\phi}(x, t; z)$ can be considered a formal solution to \eqref{eq02}.
			
			\item The function $\widehat{\phi }_t( x,t;z) \widehat{\phi }^{-1}( x,t;z)$ is analytic for $z\in\mathbb{C}\backslash\{\pm\mu/2\}$ for the following reasons:
			
			Analogous to the previous argument, we replace the equivalence \eqref{eq18} with
			\begin{align*}%\label{eq20}
				\begin{aligned}
					\widehat{\phi }_t( x,t;z) \widehat{\phi }^{-1}( x,t;z)=\widehat{\psi }_t( x,t;z) \widehat{\psi }^{-1}( x,t;z),
				\end{aligned}
			\end{align*}
			which demonstrates that $\widehat{\phi }_t( x,t;z) \widehat{\phi }^{-1}( x,t;z)$ is jump-free on $z\in \Sigma_1$ and that each pole $z_k\in \mathcal{Z}\cup\overline{\mathcal{Z}}$ is removable. Furthermore, as $|z|\to\infty$, we obtain the asymptotic behavior
			\begin{align*}%\label{eq21}
				\begin{aligned}
					\widehat{\phi }_t( x,t;z) \widehat{\phi }^{-1}( x,t;z)=\mathcal{O}(z^{-1}), 
				\end{aligned}
			\end{align*}
			since the function $\widehat{\phi }( x,t;z)e^{-i\varphi(z;x,0)\sigma_{3}}$ shares the same asymptotics as $M(z;x,t)$ as $|z|\to\infty$.
			From \eqref{eq12}, it follows that $\widehat{\phi }_t( x,t;z) \widehat{\phi }^{-1}( x,t;z)$, which has simple poles at $z=\pm\mu/2$,  behaves formally as, by applying Liouville's theorem again
			\begin{align*}%\label{eq22}
				\begin{aligned}
					\widehat{\phi }_t( x,t;z) \widehat{\phi }^{-1}( x,t;z) =\frac{iz}{4z^{2}-\mu^{2}}M(z;x,t) \sigma _{3}M^{-1}( z;x,t) .
				\end{aligned}
			\end{align*}
			Therefore, $\widehat{\phi }( x,t;z)$ can also be regarded as a formal solution to \eqref{eq03}. Moreover, the solutions to the RMB equations \eqref{eq01}  can be reconstructed using the fromulae in \eqref{eq11}, which is made possible by compatibility.
			
			\item Finally, We verify that the initial data $E_0(x)$ can be reconstructed from the first formula of in \eqref{eq11}. When $t=0$, the essential singularities $z=\pm\mu/2$ of the factor $e^{\pm2i\varphi(z)}$ vanish from the jump matrices \eqref{eq9}. Now, via the transformation $\widehat{M}( z;x,0):={M}( z;x,0)L(z;x)$ with
			\begin{align*}%\label{eq23}
				\begin{aligned}
					L(z;x) :=\begin{cases} \begin{pmatrix}
							1 & 0 \\
							r( \RRe z) e^{-2i\varphi(z;x,0)}& 1
						\end{pmatrix}&\left| z\pm \frac{\mu}{2}\right|<\kappa,~~\IIm z >0\\
						\begin{pmatrix}
							1 & \overline{r}( \RRe z) e^{2i\varphi(z;x,0)}\\
							0 & 1
						\end{pmatrix}^{-1}&\left| z\pm \frac{\mu}{2}\right|<\kappa ,~~\IIm z <0\\
						I&otherwise, \end{cases}
				\end{aligned}
			\end{align*}
			the jump contours $\cup_{k=1}^4\Sigma_k^\kappa$ of $M(z;x,0)$ are deformed into the contours $$\left\{z\in\mathbb{R}:  0\leq|z\pm \mu/2|<\kappa \right\}$$ corresponding to $\widehat{M}( z;x,0)$. The function $\widehat{M}( z;x,0)$ only has jump on $z\in \mathbb R$ with
			\begin{align*}%\label{eq24}
				\begin{aligned}
					\widehat{M}_+( z;x,0)=\widehat{M}_-( z;x,0)
					\begin{pmatrix}
						1+| r(z) | ^{2} & \overline{r}(z) e^{2i\varphi(z;x,0)} \\
						r( z) e^{-2i\varphi(z;x,0)}  & 1
					\end{pmatrix},
				\end{aligned}
			\end{align*}
			and it shares both the residue conditions from \eqref{eq10} with $t=0$ and the asymptotic behavior of ${M}( z;x,0)$ at $z=\infty$. Therefore, $\widehat{M}( z;x,0)$ satisfies a RHP associated with the Zakharov-Shabat scattering problem with solitons as $t=0$, from which the initial data can be obtained by
			\begin{align*}%\label{eq25}
				\begin{aligned}
					E_0(x) =-4i\lim _{z\rightarrow \infty }\left( z \widehat M ( z;x,0) \right)_{12}.
				\end{aligned}
			\end{align*}
		\end{itemize}
	\end{proof}
	
	\begin{Remark}
		The existence and uniqueness of solutions to RHP \ref{RHP1} for all $(x, t) \in \mathbb{R}^2$ are guaranteed by Liouville’s theorem and Zhou’s vanishing lemma \cite{MR1000732}.  
	\end{Remark}
	
	\begin{theo}\label{theo1}
		Given the initial electric field $E_0(x) = E(x, 0)$ satisfying Assumption \ref{Ass1}, let $\sigma:=\{r(z),\{(z_k,c_k)\}_{k=1}^N\}$ be the scattering data corresponding to $E_0$, generated from the direct scattering transform of \eqref{eq02}. Let $E(x,t)$, $s(x,t)$, $u(x,t)$, and $r(x,t)$ be the solutions of Cauchy problem \eqref{eq01} recovered by $M(z;x,t)$ through the reconstructed formulae provided in Theorem \ref{theo2}. 
		
		For fixed $\mu, x_1, x_2, v_1, v_2 \in \mathbb{R}$ satisfying $\mu \in (0,1]$, $x_1 \leq x_2$, and $-1/\mu^2 < v_1 \leq v_2 < 0$, define the cone 
		$$C( x_{1},x_{2},v_{1},v_{2}) :=\{ ( x,t) :x = x_{0}+vt,x_{0}\in [ x_{1},x_{2}] ,v\in[ v_{1},v_{2}] \},$$
		and the interval
		$$\mathcal 	J := i[(-\mu^{-2}-v_1^{-1})^{1/2}/2,  (-\mu^{-2}-v_2^{-1})^{1/2}/2].$$
		Let 
		\begin{align*}%\label{eq31}
			\begin{aligned}
				\zeta_0 = ( \mu^2/4-t/(8x) ( 1 + \sqrt{1 - 8 \mu^2 x/t} ) )^{1/2},~~
				\beta =( 4\zeta_0^{2}-\mu ^{2}) ^{3}/( 4\zeta_{0}^{3}+3\mu ^{2}\zeta _{0}).
			\end{aligned}
		\end{align*}
		
		Then, for  $(x,t)\in C(\mu, x_1, x_2, v_1, v_2)$ as $t \to \infty$, the following asymptotic behaviors hold:
		\begin{align}\label{THE1}
			\begin{aligned}
			E(x,t) &= E(x,t; \sigma_\delta(\mathcal J)) + t^{-1/2} \sqrt{\beta} ( f_1(x,t) + f_2(x,t) ) +\mathcal O(t^{-\gamma_\circ}), \\
			s(x,t) &= s(x,t; \sigma_\delta(\mathcal J)) + \mathcal O(t^{-1/2}), \\
			u(x,t)& = u(x,t; \sigma_\delta(\mathcal J)) + \mathcal O(t^{-1/2}), \\
			r(x,t) &= r(x,t; \sigma_\delta(\mathcal J)) + \mathcal O(t^{-1/2}),
		\end{aligned}
		\end{align}
		where the constant $\gamma_\circ>1/2$, and $E(x,t;\sigma_\delta(\mathcal J)), u(x,t;\sigma_\delta(\mathcal J)), r(x,t; \sigma_\delta(\mathcal J))$ represent the superposition of  $N(\mathcal J)$ single-soliton, and $s(x,t; \sigma_\delta(\mathcal J))$ is the superposition of  $N(\mathcal J)$ single-kink, corresponding to the modified discrete scattering data
		$$\sigma_\delta(\mathcal J)=\{\{(z_k,c_{\delta,k}= c_k\delta^{-2}(z_k))\}_{k=1}^N:z_k\in \mathcal {Z\cap J} \},$$
		with
		$$\delta( z) =\exp{\left(i\int _{\mathbb{R} \backslash[ -\zeta_{0},\zeta _{0}] }\frac{\nu(s)}{s-z}ds\right)},~~ \nu(z) = -\frac{1}{{2\pi}}\ln (1+| r( z) | ^{2}).$$
		Moreover,
		\begin{align*}
			f_1(x,t)&=	\left( M_{\mathrm{out}}(\zeta_0) \begin{pmatrix}
				0 &\beta_{12}\delta_0^A \\
				\overline{\beta_{12}\delta_0^A}& 0
			\end{pmatrix} M^{-1}_{\mathrm{out}}(\zeta_0) \right) _{12}\\
			f_2(x,t)&=\left(M_{\mathrm{out}}(-\zeta_0) \begin{pmatrix}
				0 &\overline{\beta_{12}\delta_0^A} \\
				\beta_{12}\delta_0^A& 0
			\end{pmatrix} M^{-1}_{\mathrm{out}}(-\zeta_0)\right) _{12},
		\end{align*}
		where $M_{\mathrm {out}}(z)=M^{\bigtriangleup _{\mathcal J}}( z;\sigma _{\delta }(\mathcal J))$ solves RHP \ref{RHP10} with index set $\bigtriangleup=\bigtriangleup _{\mathcal J}$ and the scattering data $\sigma_\delta(\mathcal J)$ (see Proposition \ref{prop5.1}). The modulus and argument of $\beta_{12}\delta_0^A$ are given by 
		\begin{align*}%\label{eq109.1}
			|\beta_{12}\delta_0^A|^2=|\nu(\zeta_0)|,
		\end{align*}
		and 
		\begin{align*}%\label{eq109.2}
			\begin{aligned}\arg ( \beta _{12}\delta _{0}^{A}) &=-\pi/4-\arg r(\zeta_0)-\arg  \Gamma ( i\nu(\zeta_0)) -2\nu(\zeta_0) \ln ( \sqrt{\beta }/( 8\zeta_0\sqrt{t}) ) \\
				&+4\sum _{k\in \bigtriangleup }\arg ( \zeta_0-z_{k}) +2\int _{\zeta_0}^{\infty }\ln \left( \frac{s+\zeta_0}{s-\zeta_0}\right) d_{s}\nu (s) +\frac{16t\zeta_0^{3}}{( 4\zeta_0^{2}-\mu ^{2}) ^{2}},\end{aligned}
		\end{align*}
		where $\Gamma$ is the gamma function.
	\end{theo}
	
	\section{Deformation of the Riemann-Hilbert problem}\label{S3}
	In this section, we try to deform RHP \ref{RHP1} to consider the long-time asymptotic behaviors of the solutions to the Cauchy problem \eqref{eq01} along the characteristic $x=vt$ in soliton region. 
	
	The $z$-complex plane is divided based on the decay behavior of the oscillatory factors $e^{\pm2it\theta(z)}$ as $t\to\infty$. Inspired by the $\overline \partial$ generalization of the nonlinear steepest descent method presented in \cite{MR3795020}, we introduce two key transformation matrices, $F^{-\sigma_{3}}(z)$ and $R^{(2)}(z)$. Through the transformations \eqref{eq35} and \eqref{eq63}, the original jump contour $\Sigma_{1}$ is successfully deformed into $\Sigma_2$, where the factors $e^{\pm2it\theta(z)}$ exhibit exponential decay along the redefined path. Meanwhile, the residue conditions are effectively managed with proper asymptotic control.
	
	\subsection{Partitioning of the $z$--complex plane}
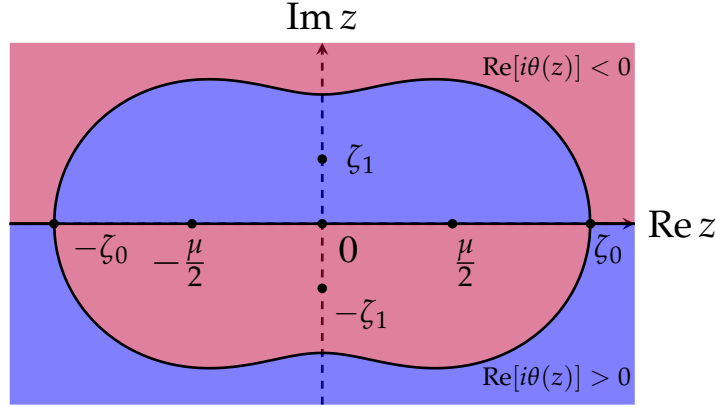
\begin{figure}[H] 
		\centering 
		       \begin{tikzpicture}    [scale=1.2, transform shape]       
			\begin{axis}[axis lines=middle,
				unit vector ratio=1 1,
				ymax=0.7,ymin=-0.7,
				xmax=1.2,xmin=-1.2,
				xlabel={$\RRe z$},ylabel={$\IIm z$},
				xlabel style={right},ylabel style={above},
				axis line style={dashed, line width=1pt},
				ytick={1},xtick={-2,2},
				]  
				\addplot+[name path=A, black, samples=1000,mark=none,domain=-1.2:1.2,line width=1pt] {sqrt(max(0,-4*x^2 + sqrt(16*x^2 + 1)))/2}; 
				\addplot+[name path=B, black, samples=1000,mark=none,domain=-1.2:1.2,line width=1pt] {-sqrt(max(0,-4*x^2 + sqrt(16*x^2 + 1)))/2};
				
				\addplot [name path=C, draw=none] coordinates {(-1.03,0) (1.03,0)};
				% 填充曲线与 x 轴之间的区域
				\addplot [
				purple,
				opacity=0.5
				] fill between[
				of=C and B,
				soft clip={domain=-1.03:1.03}
				];
				\addplot [
				blue,
				opacity=0.5
				] fill between[
				of=C and A,
				soft clip={domain=-1.03:1.03}
				];

				% 用极大的矩形覆盖整个上半平面
				\addplot [
				domain=-1.2:1.2,
				samples=2,
				name path=upperhalf,
				draw=none
				] {0.7}; % 上边界
				
				\addplot [
				domain=-1.2:1.2,
				samples=2,
				name path=lowerhalf,
				draw=none
				] {-0.7}; % 下边界

				\addplot [
				purple,
				opacity=0.5
				] fill between[
				of=A and upperhalf,
				];
				
				\addplot [
				blue, 
				opacity=0.5
				] fill between[
				of= B and lowerhalf,
				];
				
				\draw[line width=1pt] (-1.2,0) -- (1.17,0);
				
				\fill (0, 0.25) circle (1.5pt);
				\node[black, font=\footnotesize ] at (0.15, 0.25) {\(\zeta_1\)};
				\fill (0, -0.25) circle (1.5pt);
				\node[black, font=\footnotesize] at (0.15, -0.35) {\(-\zeta_1\)};
				
				\fill (1.03, 0) circle (1.5pt);
				\node[black, font=\footnotesize] at (1.1, -0.1) {\(\zeta_0\)};
				\fill (-1.03, 0) circle (1.5pt);
				\node[black, font=\footnotesize] at (-0.85, -0.1) {\(-\zeta_0\)};
				
				\fill (0.5, 0) circle (1.5pt);
				\node[black] at (0.55, -0.15) {\(\frac{\mu}{2}\)};
				\fill (-0.5, 0) circle (1.5pt);
				\node[black] at (-0.55, -0.15) {\(-\frac{\mu}{2}\)};
				
				\fill (0, 0) circle (1.5pt);
				\node[black] at (0.1, -0.1) {\(\text{0}\)};
				
				\node[black, font=\scriptsize] at (0.9, 0.6) {\(\RRe [i\theta(z)]<0\)};
				\node[black, font=\scriptsize ] at (0.9, -0.6) {\(\RRe [i\theta(z)]>0\)};
			\end{axis}
		\end{tikzpicture}
		\caption{ \small{The signature of $\RRe [i\theta(z)]$, the contour $G (\RRe z,\IIm z)=0$, and the stationary phase points $\pm \zeta_0$ and $\pm \zeta_1$ in the complex plane for  ${x}/{t}\in( -{1}/\mu ^{2},0) $ with $\mu \in ( 0,1]$.} }\label{F2}
	\end{figure}

	Recall the reflectionless 1-soliton solution for the electric field
	\begin{align}\label{eq31}
		\begin{aligned}
			E( x,t) =   4\eta \mathrm{sign} (c)\sech\left( 2\eta \left( \frac{1}{4\eta^{2}+\mu^{2}}t+x\right)+\ln 2\eta -\ln|c|\right) ,
		\end{aligned}
	\end{align}
	where the discrete spectrum is assumed to be $ \mathcal{Z} := \{i\eta\}$ with $\eta>0$, and the corresponding norming constant is $C:=\{ic\}$. From \eqref{eq31}, the 1-soliton has velocity $-{1}/(4\eta^{2}+\mu^{2})$. Here, we restrict our attention to the soliton region, i.e, $-{1}/{\mu^{2}} <{x}/{t} <0$. Denoting $\varphi(z;x,t):=t\theta(z;x,t)$, a direct calculation yields
	\begin{align}\label{eq32}
		\begin{aligned}
			\RRe\left[ i\theta ( z) \right]=\IIm zG(\RRe z,\IIm z) ,
		\end{aligned}
	\end{align}
	where
	\begin{align*}%\label{eq33}
		\begin{aligned}
			G ( a,b) :=\frac{4( a^{2}+b^{2}) +\mu^{2}}{\left(( 2a+\mu ) ^{2}+4b^{2}\right) \left( ( 2a-\mu ) ^{2}+4b^{2}\right) }+\frac{x}{t}.
		\end{aligned}
	\end{align*}
	For $-{1}/{\mu^{2}} <{x}/{t} <0$,  the complex $z$-plane can be  partitioned into four disjoint regions based on the signature of the function $\RRe\left[i\theta (z) \right]$, as shown in Figure \ref{F2}. The points $\pm \zeta_0$ and $\pm \zeta_1$ represent the four stationary phase points of the function $e^{i\theta(z)}$, situated on the real and imaginary axes, respectively.
	\begin{Remark}\label{Remark1}
		For  ${x}/{t}\in( -{1}/\mu ^{2},0) $ with $\mu \in ( 0,1] $, the positions of stationary phase points $\pm \zeta_0$ and $\pm \zeta_1$ relative to $\pm \mu/2$ in the complex $z$-plane are described as follows: 
		\begin{itemize}
			\item ${\sqrt{3}}\mu/{2}<\zeta_0 $, and $G ( \pm \zeta_0,0) =0$, i.e., $\RRe\left[ i\theta ( \pm \zeta_0) \right]=0$.
			\item $\zeta_1$ remains purely imaginary with $0<|\zeta_1|<{\mu}/{2}$, and $\RRe\left[ i\theta ( \pm \zeta_1) \right]\neq0$.
		\end{itemize}
	\end{Remark}
	Indeed,  for $(\mu, \lambda)\in  ( 0,1]\times( -{1}/{\mu ^{2}},0)$, we define 
	\begin{align*}%\label{eq34}
		\begin{aligned}
			f_{\pm }( \mu ,\lambda ) :=\frac{1}{4}\mu ^{2}+\frac{1}{8\lambda}\left( -1\pm \sqrt{1-8\mu ^{2}\lambda}\right) .
		\end{aligned}
	\end{align*}
	For fixed $\mu$, the function $f_{+}(\mu,\lambda)$ is monotonically decreasing with respect to $ \lambda$ and satisfies $-\mu^2/4<f_{+}(\mu,\lambda)<0.$ Conversely, $f_{-}(\mu,\lambda)$ increases monotonically with respect to $\lambda$ and satisfies $3\mu^2/4<f_{-}(\mu,\lambda)$. Using the identities $\zeta_1^{2}= f_{+ }( \mu ,{x}/{t})$ and $\zeta_0^{2}= f_{- }( \mu ,{x}/{t})$ to make available the above statements in Remark \ref {Remark1}. It follows that for varying parameters $\mu $ and ${x}/{t}$, the sign pattern of $\RRe[ i\theta ( z)]$ undergoes only a scaling transformation, while the points $z=\pm{\mu}/{2}$ always lie within the contour defined by $G (\RRe z,\IIm z)=0$. 
	
	\subsection{Separation of oscillatory factors in the jump matrix}
	To separate the  oscillatory factors $ e^{\pm 2i\varphi} $ in the jump matrix $ V(z) $  into distinct triangular matrices via upper/lower triangular factorization, and to control the residue conditions as $t \to \infty $, we introduce the transformation
	\begin{align}\label{eq35}
		\begin{aligned}
			M ^{(1)}( z;x,t) =M( z;x,t) F ( z)^{-\sigma _{3}},
		\end{aligned}
	\end{align}
	where
	\begin{align*}%\label{eq36}
		\begin{aligned}
			F( z) :=\prod _{k\in \bigtriangleup }\frac{z-\overline{z}_{k}}{z-z_{k}}\delta ( z),
		\end{aligned}
	\end{align*}
	with the index set $\bigtriangleup$ and the uniformly bounded scalar function $\delta(z)$ are defined as
	\begin{align}\label{eq37}
		\begin{aligned} \bigtriangleup:=\left\{k\in\{1,2,\dots,N\}:z_k\in\mathcal{Z},G(0,\IIm z_k)<0\right\},
		\end{aligned}
	\end{align}	
	and 
	\begin{align}\label{eq38}
		\begin{aligned} 
			\delta( z) :=e^{i\int _{\mathbb{R} \backslash[ -\zeta_{0},\zeta _{0}] }\frac{\nu(s)}{s-z}ds}
			,\text{with~} \nu(z) = -\frac{1}{2\pi}\ln (1+| r( z) | ^{2}).
		\end{aligned}
	\end{align}	
	It can be readily verified that the function $F(z)$ exhibits a jump discontinuity on $ z\in \mathbb{R} \backslash[ -\zeta_{0},\zeta _{0}]$, with $F_+(z)=F_-(z)(1+|r(z)|^2)$, is meromorphic in the complement of its jump contours, satisfies the symmetry $\overline {F}(\overline z)=1/F(z)$, and tends to 1 as $|z|\to \infty$ with $\IIm z\neq 0$. As a result, $M ^{(1)}(z; x,t)$ satisfies the following RHP that it is well-condition as $t\to \infty$. 
	
	\begin{rhp}\label{RHP2}  For given scattering data $\sigma$, seek for a meromorphic matrix-valued function $M^{(1)}(z)=M^{(1)}(z;x,t)$ defined on  $z\in \mathbb{C}\backslash\{\Sigma_1\}$ that satisfies the following conditions.
		\begin{itemize}
			\item $M^{(1)}(z)$ has continuous boundary values $M^{(1)}_\pm(z)$ on $z\in \Sigma_1$, which satisfy the jump conditions
			\begin{align*}%\label{eq39}
				\begin{aligned}
					M^{(1)}_+(z)=M^{(1)}_-(z)V^{(1)}(z;x,t),
				\end{aligned}
			\end{align*}
			where
			\begin{align}\label{eq40}
				\begin{aligned} 
					V^{(1)}( z) :=\begin{cases} \begin{pmatrix}
							1 & 0 \\
							\frac{r( z) }{1+| r(z)|^2}F_-^{-2}( z) e^{ -2it\theta(z) } & 1
						\end{pmatrix}
						\begin{pmatrix}
							1 & \frac{\overline{r}(z) }{1+|r( z)| ^2} F_{+}^{2}( z) e^{ 2it\theta(z) } \\
							0 & 1
						\end{pmatrix}
						~~	z\in\mathbb{R} \backslash[ -\zeta_{0},\zeta _{0}]\\
						\begin{pmatrix}
							1 & \overline{r}(z) F^{2}( z) e^{ 2it\theta(z) } \\
							0 & 1
						\end{pmatrix}
						\begin{pmatrix}
							1 & 0 \\
							r( z) F^{-2}( z) e^{ -2it\theta(z) } & 1
						\end{pmatrix}
						~~	z\in( -\zeta_{0},\zeta_0) \backslash \overline D( \pm \frac{\mu}{2},\kappa )\\
						\begin{pmatrix}
							1 & 0 \\
							r( \RRe z) F^{-2}(z) e^{ -2it\theta(z) } & 1
						\end{pmatrix}\quad \quad \quad \quad \quad \quad \quad \quad \quad ~ z\in\Sigma^\kappa_1\cup\Sigma^\kappa_4\\
						\begin{pmatrix}
							1 & \overline{r}( \RRe z) F^{2}( z) e^{ 2it\theta(z) } \\
							0 & 1
						\end{pmatrix}\quad \quad \quad \quad \quad \quad \quad \quad \quad \quad~~ z\in\Sigma^\kappa_2\cup\Sigma^\kappa_3. \end{cases}
				\end{aligned}		
			\end{align}
			Moreover, $M^{(1)}(z)$ satisfies the following residue conditions at each simple poles $z_k$ and $\overline z_k$
			\begin{align*}%\label{eq41}
				\begin{aligned}
					\Res \limits _{z=z_{k}}M^{(1)}( z) =
					\begin{cases}
						\lim\limits_{z\rightarrow z_{k}}M^{(1)}( z) 
						\begin{pmatrix}
							0 & c_{k}^{-1}[ ( F^{-1}) '( z_{k})] ^{-2}e^{2it\theta( z_{k}) } \\
							0 & 0
						\end{pmatrix}
						&k\in \bigtriangleup\\
						\lim\limits_{z\rightarrow z_{k}}M^{(1)}( z) 
						\begin{pmatrix}
							0 & 0 \\
							c_{k} F^{-2}( z_{k})e^{-2it\theta( z_{k}) } & 0
						\end{pmatrix}
						&k\in \bigtriangledown\end{cases}\\
					\Res \limits_{z= \overline z_{k}}M^{(1)}(z)= \begin{cases} 
						\lim\limits_{z\rightarrow \overline z_{k}}M ^{(1)}( z) 
						\begin{pmatrix}
							0 & 0 \\
							-\overline c_{k}^{-1}( F'(\overline{z}_{k})) ^{-2}e^{ -2it\theta (\overline{z}_{k}) } & 0
						\end{pmatrix}
						&k\in \bigtriangleup\\
						\lim\limits_{z\rightarrow \overline z_{k}}M ^{(1)}( z) 
						\begin{pmatrix}
							0 & -\overline c_{k}  F^{2}( \overline{z}_{k})e^{ 2it\theta (\overline{z}_{k}) } \\
							0 & 0 
						\end{pmatrix}
						&k\in \bigtriangledown,\end{cases}
				\end{aligned}
			\end{align*}
			where the index set $\bigtriangledown:=\{1,2\dots,N\}\backslash\bigtriangleup.$		
			\item  The matrix  $M^{(1)}(z)$ possesses the asymptotic behavior $M^{(1)}(z)\to I$, as $|z|\to \infty$.
		\end{itemize}
	\end{rhp}
	
	Furthermore, we replace the function $\delta(z)$ in \eqref{eq38} with
	\begin{align}\label{eq42}
		\begin{aligned}
			\delta (z) :=\left( \frac{z+\zeta_0}{\zeta_0-z}\right) ^{i\nu( \zeta_0)}e^{i\beta( z,\zeta _{0}) }
		\end{aligned},
	\end{align}
	where 
	$$\beta( z,\zeta _{0}):=\int 
	_{\zeta_0}^{\infty }\left( \nu( s) -\chi_{(\zeta_0,\zeta _{0}+1)}(s)\nu(\zeta_0)\right) \left( \frac{1}{s-z}-
	\frac{1}{s+z}\right) ds+\nu( \zeta_0) \ln \frac{\zeta_0 +1-z}{\zeta_0 +1+z},
	$$ 
	with   $ \chi_{(\zeta_0,\zeta _{0}+1)}(s)$ denoting the characteristic function of the interval $(\zeta_0,\zeta _{0}+1)$, and the logarithm being principally branched along $(-\infty,-\zeta _{0}-1)\cup (\zeta_0+1,\infty)$. 
	
	To clarify, we introduce the notation
	$$\widehat{\delta}(z):=e^{i\int ^{\infty }_{\zeta_0}\frac{\nu(s)}{s-z}ds}= ( \zeta _{0}-z) ^{-i\nu(\zeta_0)}e^{i\int _{\zeta_0}^{\infty }\frac{\nu(s) -\chi_{(\zeta_0,\zeta_0+1)}(s)\nu(\zeta_0)}{s-z}ds+i\nu(\zeta_0)\ln(\zeta _{0}+1-z)},
	$$
	and note that the reflection coefficient satisfies the symmetry $r(z)=\overline {r}(-z)$, which arises from the symmetries of the Jost solutions $\phi_\pm(x,z)=\overline {\phi_\pm}(x,-z)$. Consequently, we obtain $\nu(z)=\nu(-z)$ and write  $\delta(z)=\widehat{\delta}(-z)^{-1}\widehat{\delta}(z)$, which corresponds to \eqref{eq42}.
	\begin{prop}
		For a fixed $0<|\alpha|<\pi $, consider $z\to \pm \zeta_0$ along the rays $z=\pm \zeta_0+\gamma e^{i\alpha}$ with $\gamma\in [0,\infty)$. Then,  there exists constants $c$ such that
		\begin{itemize}
			\item As $z\to \zeta_0$, 
			\begin{align*}%\label{eq43}
				\begin{aligned}
					|\left( (z+\zeta_0)/(\zeta_0-z)\right) ^{-i\nu( \zeta_0)}(F(z)-F_A(z))|\leq c|z-\zeta_0|^\frac{1}{2}, 
				\end{aligned}
			\end{align*}
			where 
			\begin{align*}%\label{eq44}
				\begin{aligned}
					F_A(z) :=\left( \frac{z+\zeta_0}{\zeta_0-z}\right) ^{i\nu( \zeta_0)}\prod _{k\in \bigtriangleup }\frac{\zeta_0-\overline{z}_{k}}{\zeta_0-z_{k}}e^{i\beta( \zeta_0,\zeta _{0}) }.
				\end{aligned}
			\end{align*}
			\item As $z\to -\zeta_0$,
			\begin{align}\label{eq45}
				\begin{aligned}
					|\left( (z+\zeta_0)/(\zeta_0-z)\right) ^{-i\nu( \zeta_0)}(F(z)-F_B(z))|\leq c|z+\zeta_0|^\frac{1}{2},
				\end{aligned}
			\end{align}
		\end{itemize}
		where 
		\begin{align*}%\label{eq46}
			\begin{aligned}
				F_B(z) :=\left( \frac{z+\zeta_0}{\zeta_0-z}\right) ^{i\nu( \zeta_0)}\prod _{k\in \bigtriangleup }\frac{\zeta_0+\overline{z}_{k}}{\zeta_0+z_{k}}e^{-i\beta( \zeta_0,\zeta _{0}) }.
			\end{aligned}
		\end{align*}
	\end{prop}
	\begin{proof}
		Let $y(z):=\nu(z)-\chi_{(\zeta_0,\zeta_0+1)}(z)\nu(\zeta_0)$. Then, $y(z)\in H^{1}(\mathbb{R})$, which follows from the fact that $r(z)\in H^{1,1}(\mathbb{R})$. Thus, by applying the Lemma 23.3 from \cite{MR954382}, we yield the estimate
		$$|\beta(z,\zeta_0)-\beta(\zeta_0,\zeta_0)|\leq c_1\sqrt{2}\left\|y'\right\|_{L^2(\zeta_0,\infty)}|z-\zeta_0|^{\frac{1}{2}}.$$
	\end{proof}
	
	\subsection{The continuous extensions of jump factorization}
	
	To deform the segment $\Sigma_{\mathrm {re}} $ onto the steepest descent contours for $e^{\pm it \theta(z)}$, we construct a continuous (but non-analytic) extension $R^{(2)}(z;x,t)$ of the jump matrix factors in \eqref{eq40} off the real line, following the methodology of \cite{MR3795020}. 
	
\begin{figure}[h] 
		\centering 
		 \begin{tikzpicture}
			% RRe z 
			\draw[black, thick, dashed] (-8,0) -- (-3.5,0);
			\draw[black, thick, dashed] (-2.5,0) -- (2.5,0);
			\draw[black, thick, dashed] (3.5,0) -- (8,0) node[right] {\text{Re} $z$};
			%zeta1
			\fill (0,1.5) circle (1.5pt);
			\node[black] at (0, 2) {\(\zeta_1\)};
			%-zeta1
			\fill (0,-1.5) circle (1.5pt);
			\node[black] at (0, -2) {\(-\zeta_1\)};
			%left
			%sigma4c
			\draw[black, line width=1pt] (-2.5,0) -- (0,1.5);
			\draw[black, -{Stealth},line width=1pt] (-2.5,0) -- (-1,0.9);
			\node[black] at (-1.5,1.1) {\(\Sigma_4^C\)};
			%omega12
			\node[black] at (-1.3, 0.3) {\(\Omega_{12}\)};
			
			%sigma3c
			\draw[black, line width=1pt] (-2.5,0) -- (0,-1.5);
			\draw[black,-{Stealth}, line width=1pt] (-2.5,0) -- (-1, -0.9);
			\node[black] at (-1.5,-1.1) {\(\Sigma_3^C\)};
			%omega11
			\node[black] at (-1.3, -0.3) {\(\Omega_{11}\)};
			
			%sigma6D
			\draw[black, line width=1pt] (0,1.5) -- (0,0);
			\draw[black, -{Stealth},line width=1pt] (0,1.5) -- (0,0.5);
			\node[black] at (0.4, 0.8) {\(\Sigma_6^D\)};
			
			%sigma5D
			\draw[black, line width=1pt] (0,-1.5) -- (0,0);
			\draw[black, -{Stealth},line width=1pt] (0,-1.5) -- (0,-0.5);
			\node[black] at (0.4, -0.8) {\(\Sigma_5^D\)};
			
			%sigma1C
			\draw[black, line width=1pt] (2.5,0) -- (0,1.5);
			\draw[black, -{Stealth}, line width=1pt](0,1.5) -- (1.5,0.6);
			\node[black] at (1.5,1.1) {\(\Sigma_1^C\)};
			%omega9
			\node[black] at (1.2, 0.3) {\(\Omega_{9}\)};
			
			%sigma2C
			\draw[black, line width=1pt] (2.5,0) -- (0,-1.5);
			\draw[black,-{Stealth}, line width=1pt] (0,-1.5) -- (1.5, -0.6);
			\node[black] at (1.5,-1.1) {\(\Sigma_2^C\)};
			%omega10
			\node[black] at (1.3, -0.3) {\(\Omega_{10}\)};
			
			%sigma4D
			\draw[black, line width=1pt] (-3.5,2) -- (-3.5,0);
			\draw[black, -{Stealth},line width=1pt] (-3.5,2) -- (-3.5,1);
			\node[black] at (-3,1.8) {\(\Sigma_4^D\)};
			
			%sigma3D
			\draw[black, line width=1pt] (-3.5, -2) -- (-3.5,0);
			\draw[black, -{Stealth},line width=1pt] (-3.5, -2) -- (-3.5,-1);
			\node[black] at (-3, -1.8) {\(\Sigma_3^D\)};
			
			%sigmaB
			\draw[black, line width=1pt] (-8,-2.5) -- (-3.5,2);
			\draw[black, line width=1pt] (-8,2.5) -- (-3.5,-2);
			%sigma3B
			\draw[black, -{Stealth}, line width=1pt] (-8,-2.5) -- (-6.0,-0.5);
			\node[black] at (-6.5,-1.8) {\(\Sigma_3^B\)};
			%sigm1B
			\draw[black, -{Stealth}, line width=1pt] (-5.5,0) -- (-4.75,0.75);
			\node[black] at (-4.5,1.8) {\(\Sigma_1^B\)};
			%sigma4B
			\draw[black, -{Stealth}, line width=1pt] (-5.5, 0) -- (-4.75, -0.75);
			\node[black] at (-4.5,-1.8) {\(\Sigma_4^B\)};
			%sigm4B
			\draw[black, thick,black, -{Stealth}, line width=1pt] (-8,2.5) -- (-6.0,0.5);
			\node[black] at (-6.5,1.8) {\(\Sigma_2^B\)};

			%right
			%sigma1D
			\draw[black, line width=1pt] (3.5,2) -- (3.5,0);
			\draw[black, -{Stealth},line width=1pt] (3.5,0) -- (3.5,1);
			\node[black] at (3,1.8) {\(\Sigma_1^D\)};
			
			%sigma2D
			\draw[black, line width=1pt] (3.5, -2) -- (3.5,0);
			\draw[black, -{Stealth},line width=1pt] (3.5, 0) -- (3.5,-1);
			\node[black] at (3, -1.8) {\(\Sigma_2^D\)};
			
			%sigmaA
			\draw[black, line width=1pt] (8,-2.5) -- (3.5,2);
			\draw[black, line width=1pt] (8,2.5) -- (3.5,-2);
			%sigma4A
			\draw[black, -{Stealth}, line width=1pt] (5.5,0) -- (6.25,-0.75);
			\node[black] at (6.5,-1.8) {\(\Sigma_4^A\)};
			%sigm2A
			\draw[black, -{Stealth}, line width=1pt] (3.5,2) -- (5,0.5);
			\node[black] at (4.5,1.8) {\(\Sigma_2^A\)};
			%sigma1A
			\draw[black, -{Stealth}, line width=1pt] (5.5, 0) -- (6.25, 0.75);
			\node[black] at (6.5,1.8) {\(\Sigma_1^A\)};
			%sigm3A
			\draw[black, thick,black, -{Stealth}, line width=1pt] (3.5,-2) -- (5,-0.5);
			\node[black] at (4.5,-1.8) {\(\Sigma_3^A\)};
			
			% %left circle
			% \draw[thick] (-3.5,0) arc (180:0:0.5);
			% \draw[thick, -{Stealth}] (-3.5,0) arc (180:90:0.5);
			
			% \draw[thick] (-3.5,0) arc (-180:0:0.5);
			% \draw[thick, -{Stealth}] (-3.5,0) arc (-180:-90:0.5);
			
			% \fill (-3,0) circle (1pt);
			% \node[black] at (-3,-1) {\(-\frac{\mu}{2}\)};
			
			% %right circle
			% \draw[thick] (2.5,0) arc (180:0:0.5);
			% \draw[thick, -{Stealth}] (2.5,0) arc (180:90:0.5);
			
			% \draw[thick] (2.5,0) arc (-180:0:0.5);
			% \draw[thick, -{Stealth}] (2.5,0) arc (-180:-90:0.5);
			
			% \fill (3,0) circle (1pt);
			% \node[black] at (3,-1) {\(\frac{\mu}{2}\)};
			
			%left%zeta0
			\fill (-5.5,0) circle (1.5pt);
			\node[black] at (-5.5, -0.5) {\(-\zeta_0\)};
			%omega5
			\node[black] at (-4.8, 0.3) {\(\Omega_5\)};
			%omega8
			\node[black] at (-4.8, -0.3) {\(\Omega_8\)};
			
			%omega6
			\node[black] at (-6.2, 0.3) {\(\Omega_6\)};
			%omega7
			\node[black] at (-6.2, -0.3) {\(\Omega_7\)};

			%-mu-kappa
			\fill (-3.5,0) circle (1.5pt);
			\draw[->, thick, decorate, decoration={snake, amplitude=.4mm, segment length=2mm, post length=1mm}] (-2.7, -0.8) -- (-3.4,-0.1);
			\node[black] at (-2.8, -1) {\(-\frac{\mu}{2}-\kappa\)};
			%-mu+kappa
			\fill (-2.5,0) circle (1.5pt);
			\draw[->, thick, decorate, decoration={snake, amplitude=.4mm, segment length=2mm, post length=1mm}] (-2.7, 0.8) -- (-2.6,0.1);
			\node[black] at (-2.8, 1) {\(-\frac{\mu}{2}+\kappa\)};
			
			%right
			%zeta0
			\fill (5.5,0) circle (1.5pt);
			\node[black] at (5.5, -0.5) {\(\zeta_0\)};
			
			%omega2
			\node[black] at (4.8, 0.3) {\(\Omega_2\)};
			%omega3
			\node[black] at (4.8, -0.3) {\(\Omega_3\)};
			
			%omega1
			\node[black] at (6.2, 0.3) {\(\Omega_1\)};
			%omega4
			\node[black] at (6.2, -0.3) {\(\Omega_4\)};

			%mu-kappa
			\fill (3.5,0) circle (1.5pt);
			\draw[->, thick, decorate, decoration={snake, amplitude=.4mm, segment length=2mm, post length=1mm}] (2.9, -0.8) -- (3.4,-0.1);
			\node[black] at (2.8, -1) {\(\frac{\mu}{2}-\kappa\)};
			%mu+kappa
			\fill (2.5,0) circle (1.5pt);
			\draw[->, thick, decorate, decoration={snake, amplitude=.4mm, segment length=2mm, post length=1mm}] (2.7, 0.8) -- (2.6,0.1);
			\node[black] at (2.8, 1) {\(\frac{\mu}{2}+\kappa\)};

			%omega13
			\node[black] at (-1,2.3) {\(\Omega_{13}\)};
			
			% origin point O
			\fill (0, 0) circle (1.5pt);
			\node[black] at (0.3, -0.3) {\(\text{0}\)};			
		\end{tikzpicture}
		\caption{\small {The contours defined by \eqref{eq46.1} and regions $\Omega_k$ (for $k = 1, \ldots, 13$).}}\label{F3}
	\end{figure}
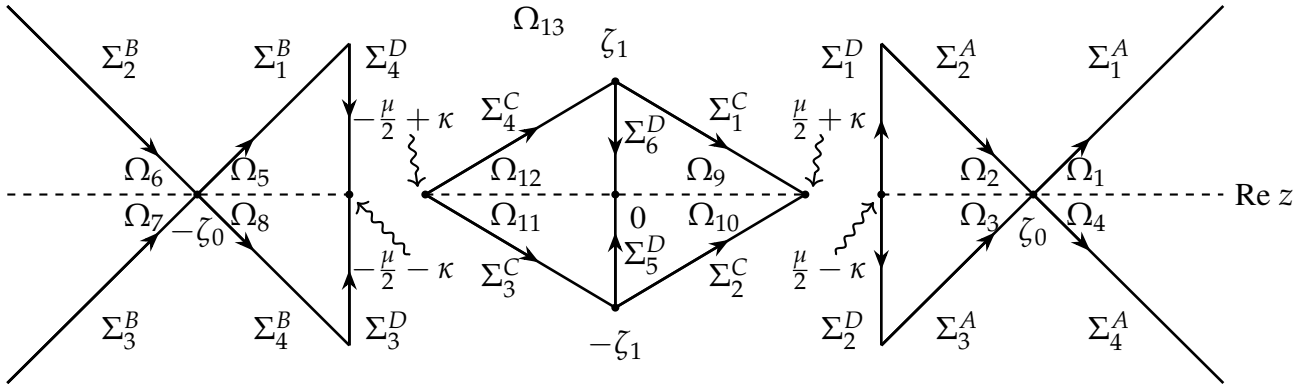

	To ensure that the deformed contour does not intersect the curve defined by $G(\RRe z, \IIm z)=0$, we define
	\begin{align}\label{eq46.1}
		\begin{aligned}
			&\Sigma_k^A:=\zeta_0+\mathbb{R^+} e^{i\frac{(2k-1)\pi}{4}}&&~~k=1,4\\
			&\Sigma_k^A:=\zeta_0+\gamma e^{i\frac{(2k-1)\pi}{4}}&\gamma\in (0,d_1),&~~k=2,3\\
			&\Sigma_k^C:=\mu/2 -\kappa+\gamma e^{i(-1)^{k-1}(\pi-\vartheta)}&\gamma \in ( 0, d_2),&~~k=1,2\\
			&\Sigma_k^C:=-\mu/2 +\kappa+\gamma e^{i(-1)^{k}\vartheta}&\gamma \in ( 0, d_2),&~~k=3,4\\
			&\Sigma_k^B:=\zeta_0+ \gamma e^{i\frac{(2k-1)\pi}{4}}&\gamma\in (0,d_1),&~~k=1,4\\
			&\Sigma_k^B:=\zeta_0+\mathbb{R^+} e^{i\frac{(2k-1)\pi}{4}}&&~~k=2,3\\
			&\Sigma_k^D:=\mu/2+\kappa+	(-1)^{k-1}i\gamma&\gamma\in (0,d_3),&~~k=1,2\\
			&\Sigma_k^D:=-\mu/2-\kappa+	(-1)^{k}i\gamma&\gamma\in (0,d_3),&~~k=3,4\\
			&\Sigma_k^D:=(-1)^{k}\gamma&\gamma\in (0,\zeta_1),&~~k=5,6,
		\end{aligned}
	\end{align}
	where 
	\begin{align*}
		\begin{aligned}
			d_1:=\sqrt{2}d_3,\quad
			d_2:=\sqrt{( {\mu }/{2}-\kappa) ^{2}-\zeta_{1}^{2}},\quad
			d_3:=\zeta_0-\mu/2-\kappa,\quad  \vartheta:=\arctan \frac{-i\zeta_{1}}{{\mu }/{2}-\kappa}.
		\end{aligned}
	\end{align*}
	The orientations of the paths are as shown in Figure \ref{F3}. We denote the regions separated by these paths and $\Sigma_{\mathrm {re}}$ as $\Omega_k$ (for $k = 1, \ldots, 13$).
	
	Additionally, let $\chi_\mathcal{Z}\in C_0^\infty$ be the cutoff function  supported on disks enclosing the discrete spectrum, defined as
	\begin{align*}
		\begin{aligned}
			\chi_\mathcal{Z}(z):=\begin{cases}
				1&z\in  \cup_{z_k \in \mathcal{Z} \cup \overline{\mathcal{Z}}} D( z_k,\kappa ),\\
				0& otherwise.
			\end{cases}
		\end{aligned}
	\end{align*}
	The extension $R^{(2)}(z)$ is constructed piecewise in the regions $\Omega_k$ for $k = 1, \ldots, 13$, with each component $R_k$ (for $k = 1, \ldots, 12$) satisfying the prescribed boundary conditions along the corresponding interfaces.
	\begin{align}
		R_{1}(z) &=\begin{cases}\label{eq47}
			\frac{\overline{r}(z)}{1+| r( z) | ^{2}}F_+^{2}(z) &  \hspace{5em} z \in(\zeta_0,\infty)\\
			\frac{\overline{r}(\zeta_0)}{1+|r(\zeta_0)|^{2}}F_{A}^{2}(z) &\hspace{5em} z\in \Sigma_1^A\end{cases}\\
		R_{4}(z) &=\begin{cases}\label{eq48}
			\frac{r( z) }{1+| r( z) | ^{2}}F_-^{-2}(z) & \hspace{4.5em} z \in(\zeta_0,\infty)\\
			\frac{r(\zeta_0)}{1+|r(\zeta_0)|^{2}}F_{A}^{-2}(z) & \hspace{4.5em} z\in \Sigma_4^A\end{cases}\\
		R_{2}(z) &=\begin{cases}\label{eq49}
			r(z) F^{-2}(z) & \hspace{6em} z\in({\mu }/{2}+\kappa,\zeta_0)\\
			r(\zeta_0)F_{A}^{-2}(z)  & \hspace{6em} z\in \Sigma _{2}^{A}\end{cases}\\
		R_{3}(z) &=\begin{cases}\label{eq50}
			\overline{r}(z) F^{2}(z) & \hspace{6.5em} z\in({\mu }/{2}+\kappa,\zeta_0)\\
			\overline{r}(\zeta_0)F_{A}^{2}(z)& \hspace{6.5em} z\in \Sigma _{3}^{A}\end{cases}
		\\
		R_{9}(z) &=\begin{cases}\label{eq51}
			r(z)F^{-2}(z) &\hspace{3.5em} z\in(0,{\mu }/{2}-\kappa)\\
			r({\mu }/{2}-\kappa) F^{-2}({\mu }/{2}-\kappa)(1-\chi_\mathcal{Z}(z)) &\hspace{3.5em} z\in  \Sigma _{1}^{C}\end{cases}
		\\
		R_{10}(z) &=\begin{cases} \label{eq52}
			\overline {r}(z)F^{2}(z) &\hspace{4.1em} z\in(0,{\mu }/{2}-\kappa)\\
			\overline {r}({\mu }/{2}-\kappa) F^{2}({\mu }/{2}-\kappa)(1-\chi_\mathcal{Z}(z)) &\hspace{4.1em} z\in  \Sigma _{2}^{C}\end{cases}\\
		R_{12}(z) &=\begin{cases}\label{eq53}
				r(z)F^{-2}(z) &\hspace{1.9em} z\in(-{\mu }/{2}+\kappa,0)\\
				r(-{\mu }/{2}+\kappa) F^{-2}(-{\mu }/{2}+\kappa)(1-\chi_\mathcal{Z}(z)) &\hspace{1.9em} z\in  \Sigma _{4}^{C}\end{cases}
			\\
		R_{11}(z) &=\begin{cases}\label{eq54}
				\overline {r}(z)F^{2}(z) &\hspace{2.5em} z\in(-{\mu }/{2}+\kappa,0)\\
				\overline {r}(-{\mu }/{2}+\kappa) F^{2}(-{\mu }/{2}+\kappa)(1-\chi_\mathcal{Z} (z)) &\hspace{2.5em} z\in  \Sigma _{3}^{C}\end{cases}\\
		R_{5}(z) &=\begin{cases}\label{eq55}r(z) F^{-2}(z) &\hspace{6em}z\in(-\zeta_0,-{\mu }/{2}-\kappa)\\
					\overline {r}(\zeta_0)F_{B}^{-2}(z) &\hspace{6em}z\in \Sigma _{1}^{B}\end{cases}
				\\
		R_{8}(z) &=\begin{cases}\label{eq56}
					\overline{r}(z) F^{2}(z) &\hspace{6.6em}z\in(-\zeta_0,-{\mu }/{2}-\kappa)\\
					r(\zeta_0)F_{B}^{2}(z) &\hspace{6.6em}z\in \Sigma _{4}^{B}\end{cases}
		\end{align}
		\begin{align}
		R_{6}(z) &=\begin{cases}\label{eq57}
			\frac{\overline{r}(z)}{1+| r( z) | ^{2}}F_+^{2}(z) &\hspace{4.9em} z \in(-\infty,-\zeta_0)\\
			\frac{r(\zeta_0)}{1+|r(\zeta_0)|^{2}}F_{B}^{2}(z)&\hspace{4.9em} z\in \Sigma_2^B\end{cases}
		\\
		R_{7}(z) &=\begin{cases}\label{eq58}
			\frac{r( z) }{1+| r( z) | ^{2}}F_-^{-2}(z) &\hspace{4.4em} z \in(-\infty,-\zeta_0)\\
			\frac{\overline{r}(\zeta_0)}{1+|r(\zeta_0)|^{2}}F_{B}^{-2}(z) &\hspace{4.4em} z\in \Sigma_3^B.
		\end{cases}
	\end{align}
	
	\begin{Lemma}\label{lemma1}
		Let the complex variable be expressed as $z=a+ib$. Given the reflection coefficient  $r(z)\in H^{1,1}(\mathbb{R})$, there exist uniformly bounded functions $R_{k}$  (for $k=1,\ldots,12$) that satisfy boundary conditions specified in \eqref{eq47}-\eqref{eq58}, with the property that their nonzero  $\overline \partial$ derivatives obey the estimations
		\begin{align}
			\begin{aligned}\label{eq59}
				&|\overline \partial R_{k}(z) | \lesssim |r'|+|a-\zeta_0|^{-1/2}\hspace{7.5em} \text{~~for~}k=1,\ldots ,4\\
				&|\overline \partial R_{k}(z) | \lesssim |r'|+|a+\zeta_0|^{-1/2}\hspace{7.5em} \text{~~for~}k=5,\ldots ,8\\
				&|\overline \partial R_{k}(z) | \lesssim |r'|+|a-\mu/2+\kappa|^{-1/2}+|\overline \partial  \chi_\mathcal{Z}(z)|\text{~~for~}k=9,10 \\
				&|\overline \partial R_{k}(z) | \lesssim |r'|+|a+\mu/2-\kappa|^{-1/2}+|\overline \partial  \chi_\mathcal{Z}(z)|\text{~~for~}k=11,12 .
			\end{aligned}
		\end{align}
	\end{Lemma}
	\begin{proof}
		Here, we provide the explicit constructions of the functions $R_{1}$, $R_{5}$, and $R_{9}$ in \eqref{eq60}-\eqref{eq62} below, while an analogous approach and the symmetry $\overline {F}(\overline z)=1/F(z)$ be applied to the remaining functions. Define the functions 
		\begin{align}
			&\begin{aligned}\label{eq60}
				R_{1}(z) :=  & \frac{\overline{r}(a)}{1 + |r(a)|^2} F^2(a + ib) \left ( 1 - \frac{b}{a - \zeta_0} \right) \\
				&+\frac{b}{a - \zeta_0} \frac{\overline{r}(\zeta_0)}{1 + |r(\zeta_0)|^2}  F_A^2(a + ib) 
			\end{aligned}& z\in \Omega_1\\
			&\begin{aligned}\label{eq61}
				R_{5}(z) := & r (a) F^{-2}( a+ib) \left( 1-\frac{b}{a+\zeta_0}\right)  
				+\frac{b}{a+\zeta_0}\overline{r}(\zeta_0) F_{B}^{-2}( a+ib)
			\end{aligned} &z\in \Omega_5\\
			&\begin{aligned}\label{eq62}
				R_{9}(z) := &\left\{r(a) F^{-2}( a+ib) \left( 1-\frac{b({\mu }/{2}-\kappa)}{{i\zeta_1}( a-{\mu }/{2}+\kappa) }\right) \right.\\
				&\left.+\frac{b({\mu }/{2}-\kappa)}{{i\zeta_1}( a-{\mu }/{2}+\kappa) }r( {\mu }/{2}-\kappa) F^{-2}( {\mu }/{2}-\kappa) \right\} (1 - \chi_\mathcal{Z}(a + ib))
			\end{aligned}&z\in \Omega_9.
		\end{align}
		We proceed to analyze in detail the properties of the function $R_{5} $. By Definition \eqref{eq61}, it can be directly verified that $ R_{5} $ satisfies the boundary condition \eqref{eq55}. Furthermore,  $ R_{5} $ is uniformly bounded, with the bound controlled by the norm $ \| r \|_{H^{1,1}(\mathbb{R})} $, since $b/(a+\zeta_0)$ remains bounded as $z \to -\zeta_0$ along any straight line within the region $\Omega_5$. Using the operator \( \bar{\partial} = \frac{1}{2}(\partial_a + i \partial_b) \), it follows that
		\begin{align*}
			\begin{aligned}
				\textstyle   \overline{\partial }R_{5}(z) &=\frac{1}{2}\left( -\frac{b}{a+\zeta_0 }+i\right) \frac{\overline{r}( \zeta _{0})F_{B}^{-2}( a+ib)  - r(a) F^{-2}( a+ib)}{a+\zeta_0}\\
				&	+ \frac{1}{2}r'(a)F^{-2}( a+ib) \left( 1-\frac{b}{a+\zeta_0}\right).
			\end{aligned}
		\end{align*}
		For the first term, applying the Cauchy–Schwarz inequality yields 
		\begin{align*}
			\begin{aligned}
				|\overline r(\zeta_0)- r(a)|=\left| \int^ {-\zeta_0}_a r'(s)ds\right|\leq \|r\|_{H^{1,1}(\mathbb{R})}|a+\zeta_0|^{{1}/{2}}.
			\end{aligned}
		\end{align*}
		Combining this estimate with \eqref{eq45}, and noting that for $z\in \Omega_5$, the inequality $|z+\zeta_0|\leq c|a+\zeta_0|$ holds  ($c$ is a constant), we conclude that there exists another constant $c > 0$ such that
		\begin{align*}
			\begin{aligned}
				|\overline{r}( \zeta _{0})F_{B}^{-2} - r(a) F^{-2}|= |\overline{r}( \zeta _{0})(F_{B}^{-2} -F^{-2})+(\overline r(\zeta_0)- r(a)) F^{-2}|\leq c|a+\zeta_0|^{{1}/{2}},\end{aligned}
		\end{align*}
		which directly establishes the validity of the inequality \eqref{eq59} for the function $R_{5}$.  The same method can be used to estimate $R_1$.
		\begin{align*}
			\begin{aligned}
				\overline{\partial }R_{9}(z) &=\left\{ \frac{\mu/2-\kappa}{2i\zeta_1}\left( -\frac{b}{a-\mu/2+\kappa}+i\right) \frac{{r}( \mu/2-\kappa )F^{-2}(  \mu/2-\kappa)  - r(a) F^{-2}( a+ib)}{a-\mu/2+\kappa}\right. \\
				&	\left.+ \frac{1}{2}r'(a)F^{-2}( a+ib) \left( 1-\frac{b({\mu }/{2}-\kappa)}{{i\zeta_1}( a-{\mu }/{2}+\kappa) }\right)\right\}(1 - \chi_\mathcal{Z}(a + ib))  \\
				&	-\left\{r(a) F^{-2}( a+ib) \left( 1-\frac{b({\mu }/{2}-\kappa)}{{i\zeta_1}( a-{\mu }/{2}+\kappa) }\right) \right.\\
				&\left.+\frac{b({\mu }/{2}-\kappa)}{{i\zeta_1}( a-{\mu }/{2}+\kappa) }r( {\mu }/{2}-\kappa) F^{-2}( {\mu }/{2}-\kappa)  \right\} \overline{\partial }\chi _{z}( a+ib). \end{aligned}
		\end{align*}
		For the function $ R_{9} $, it is essential to emphasize that $F(z)$ is continuous at the point $z ={\mu}/{2} - \kappa $. As a result, the estimations in \eqref{eq59} associated with $R_1$ and $R_9$ is established. The remainder estimations in \eqref{eq59} are treated by analogous methods.
	\end{proof}
	Employing the functions $R_k$ (for $k=1,\ldots,12$) defined in Lemma \ref{lemma1}, we construct the matrix function $R^{(2)}(z)$ as
	\begin{align*}%\label{eq63.1}
		\begin{aligned} 
			R^{(2)}(z) =\begin{cases}
				\begin{pmatrix}
					1 & R_{k}e^{2it\theta(z)} \\
					0 & 1
				\end{pmatrix}^{-1}  &z\in \Omega_k,\text{~~for ~}k=1,6\\
				\begin{pmatrix}
					1 & R_{k}e^{2it\theta(z) } \\
					0 & 1
				\end{pmatrix}  &z\in \Omega_k,\text{~~for ~}k=3,8,10,11\\
				\begin{pmatrix}
					1 & 0 \\
					R_{k}e^{-2it\theta(z)} & 1
				\end{pmatrix}^{-1}  &z\in \Omega_k,\text{~~for ~}k=2,5,9,12\\
				\begin{pmatrix}
					1 & 0 \\
					R_{k}e^{-2it\theta(z)} & 1
				\end{pmatrix}  &z\in \Omega_k,\text{~~for ~}k=4,7\\
				I   &z\in \Omega_k,\text{~~for ~}=13.\end{cases}
		\end{aligned}
	\end{align*}	
	Clearly, $R^{(2)}(z)$ is uniformly bounded and satisfies the asymptotic condition $R^{(2)}(z)\to I$ as $z\to \infty$.  
	
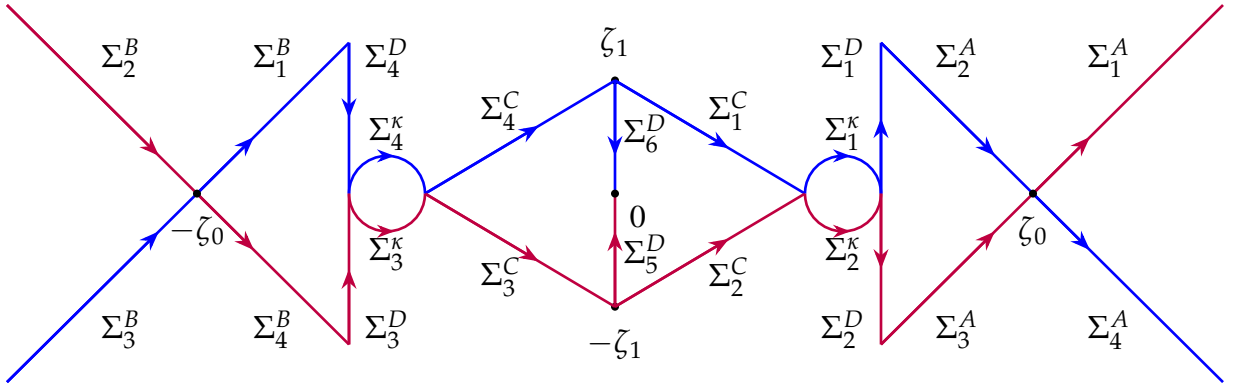
\begin{figure}[h] 
		\centering 
		\begin{tikzpicture}
			% RRe z 
			%\draw[black, thick, dashed] (-8,0) -- (-3.5,0);
			%\draw[black, thick, dashed] (-2.5,0) -- (2.5,0);
			%\draw[black, thick, dashed] (3.5,0) -- (8,0) node[right] {\text{Re} $z$};
			%zeta1
			\fill (0,1.5) circle (1.5pt);
			\node[black] at (0, 2) {\(\zeta_1\)};
			%-zeta1
			\fill (0,-1.5) circle (1.5pt);
			\node[black] at (0, -2) {\(-\zeta_1\)};
			%left
			%sigma4c
			\draw[blue, line width=1pt] (-2.5,0) -- (0,1.5);
			\draw[blue, -{Stealth},line width=1pt] (-2.5,0) -- (-1,0.9);
			\node[black] at (-1.5,1.1) {\(\Sigma_4^C\)};
			%omega12
			%\node[black] at (-1.3, 0.3) {\(\Omega_{12}\)};
			
			%sigma3c
			\draw[purple, line width=1pt] (-2.5,0) -- (0,-1.5);
			\draw[purple,-{Stealth}, line width=1pt] (-2.5,0) -- (-1, -0.9);
			\node[black] at (-1.5,-1.1) {\(\Sigma_3^C\)};
			%omega11
			%\node[black] at (-1.3, -0.3) {\(\Omega_{11}\)};
			
			%sigma6D
			\draw[blue, line width=1pt] (0,1.5) -- (0,0);
			\draw[blue, -{Stealth},line width=1pt] (0,1.5) -- (0,0.5);
			\node[black] at (0.4, 0.8) {\(\Sigma_6^D\)};
			
			%sigma5D
			\draw[purple, line width=1pt] (0,-1.5) -- (0,0);
			\draw[purple, -{Stealth},line width=1pt] (0,-1.5) -- (0,-0.5);
			\node[black] at (0.4, -0.8) {\(\Sigma_5^D\)};
			
			%sigma1C
			\draw[blue, line width=1pt] (2.5,0) -- (0,1.5);
			\draw[blue, -{Stealth}, line width=1pt](0,1.5) -- (1.5,0.6);
			\node[black] at (1.5,1.1) {\(\Sigma_1^C\)};
			%omega9
			%\node[black] at (1.2, 0.3) {\(\Omega_{9}\)};
			
			%sigma2C
			\draw[purple, line width=1pt] (2.5,0) -- (0,-1.5);
			\draw[purple,-{Stealth}, line width=1pt] (0,-1.5) -- (1.5, -0.6);
			\node[black] at (1.5,-1.1) {\(\Sigma_2^C\)};
			%omega10
			%\node[black] at (1.3, -0.3) {\(\Omega_{10}\)};
			
			%sigma4D
			\draw[blue, line width=1pt] (-3.5,2) -- (-3.5,0);
			\draw[blue, -{Stealth},line width=1pt] (-3.5,2) -- (-3.5,1);
			\node[black] at (-3,1.8) {\(\Sigma_4^D\)};
			
			%sigma3D
			\draw[purple, line width=1pt] (-3.5, -2) -- (-3.5,0);
			\draw[purple, -{Stealth},line width=1pt] (-3.5, -2) -- (-3.5,-1);
			\node[black] at (-3, -1.8) {\(\Sigma_3^D\)};
			
			%sigmaB
			\draw[blue, line width=1pt] (-8,-2.5) -- (-3.5,2);
			\draw[purple, line width=1pt] (-8,2.5) -- (-3.5,-2);
			%sigma3B
			\draw[blue, -{Stealth}, line width=1pt] (-8,-2.5) -- (-6.0,-0.5);
			\node[black] at (-6.5,-1.8) {\(\Sigma_3^B\)};
			%sigm1B
			\draw[blue, -{Stealth}, line width=1pt] (-5.5,0) -- (-4.75,0.75);
			\node[black] at (-4.5,1.8) {\(\Sigma_1^B\)};
			%sigma4B
			\draw[purple, -{Stealth}, line width=1pt] (-5.5, 0) -- (-4.75, -0.75);
			\node[black] at (-4.5,-1.8) {\(\Sigma_4^B\)};
			%sigm4B
			\draw[black, thick,purple, -{Stealth}, line width=1pt] (-8,2.5) -- (-6.0,0.5);
			\node[black] at (-6.5,1.8) {\(\Sigma_2^B\)};
			
			%right
			%sigma1D
			\draw[blue, line width=1pt] (3.5,2) -- (3.5,0);
			\draw[blue, -{Stealth},line width=1pt] (3.5,0) -- (3.5,1);
			\node[black] at (3,1.8) {\(\Sigma_1^D\)};
			
			%sigma2D
			\draw[purple, line width=1pt] (3.5, -2) -- (3.5,0);
			\draw[purple, -{Stealth},line width=1pt] (3.5, 0) -- (3.5,-1);
			\node[black] at (3, -1.8) {\(\Sigma_2^D\)};
			
			%sigmaA
			\draw[blue, line width=1pt] (8,-2.5) -- (3.5,2);
			\draw[purple, line width=1pt] (8,2.5) -- (3.5,-2);
			%sigma4A
			\draw[blue, -{Stealth}, line width=1pt] (5.5,0) -- (6.25,-0.75);
			\node[black] at (6.5,-1.8) {\(\Sigma_4^A\)};
			%sigm2A
			\draw[blue, -{Stealth}, line width=1pt] (3.5,2) -- (5,0.5);
			\node[black] at (4.5,1.8) {\(\Sigma_2^A\)};
			%sigma1A
			\draw[purple, -{Stealth}, line width=1pt] (5.5, 0) -- (6.25, 0.75);
			\node[black] at (6.5,1.8) {\(\Sigma_1^A\)};
			%sigm3A
			\draw[black, thick,purple, -{Stealth}, line width=1pt] (3.5,-2) -- (5,-0.5);
			\node[black] at (4.5,-1.8) {\(\Sigma_3^A\)};
			
			\draw[blue,line width=1pt] (2.5,0) arc (180:0:0.5);			
			\draw[blue, thick, -{Stealth}] (3,0.5) -- (3.1,0.5);			
			\node[black] at (3,0.8) {$\Sigma_1^\kappa$};			
			\draw[purple,line width=1pt] (2.5,0) arc (-180:0:0.5);			
			\draw[purple, thick, -{Stealth}] (3,-0.5) -- (3.1,-0.5);			
			\node[black] at (3,-0.8) {$\Sigma_2^\kappa$};
			\draw[blue,line width=1pt] (-3.5,0) arc (180:0:0.5);			
			\draw[blue, thick, -{Stealth}] (-3,0.5) -- (-2.9,0.5);		
			\node[black] at (-3,0.8) {$\Sigma_4^\kappa$};		
			\draw[purple, line width=1pt] (-3.5,0) arc (-180:0:0.5);		
			\draw[purple, thick, -{Stealth}] (-3,-0.5) -- (-2.9,-0.5);			
			\node[black] at (-3,-0.8) {$\Sigma_3^\kappa$};
			
			%left%zeta0
			\fill (-5.5,0) circle (1.5pt);
			\node[black] at (-5.5, -0.5) {\(-\zeta_0\)};
			%omega5
			%\node[black] at (-4.8, 0.3) {\(\Omega_5\)};
			%omega8
			%\node[black] at (-4.8, -0.3) {\(\Omega_8\)};
			
			%omega6
			%\node[black] at (-6.2, 0.3) {\(\Omega_6\)};
			%omega7
			%\node[black] at (-6.2, -0.3) {\(\Omega_7\)};
			
			%-mu-kappa
			%\fill (-3.5,0) circle (1pt);
			%\draw[->, thick, decorate, decoration={snake, amplitude=.4mm, segment length=2mm, post length=1mm}] (-2.7, -0.8) -- (-3.4,-0.1);
			%\node[black] at (-2.8, -1) {\(-\frac{\mu}{2}-\kappa\)};
			%-mu+kappa
			%\fill (-2.5,0) circle (1pt);
			%\draw[->, thick, decorate, decoration={snake, amplitude=.4mm, segment length=2mm, post length=1mm}] (-2.7, 0.8) -- (-2.6,0.1);
			%\node[black] at (-2.8, 1) {\(-\frac{\mu}{2}+\kappa\)};
			
			%right
			%zeta0
			\fill (5.5,0) circle (1.5pt);
			\node[black] at (5.5, -0.5) {\(\zeta_0\)};
			
			%omega2
			%\node[black] at (4.8, 0.3) {\(\Omega_2\)};
			%omega3
			%\node[black] at (4.8, -0.3) {\(\Omega_3\)};
			
			%omega1
			%\node[black] at (6.2, 0.3) {\(\Omega_1\)};
			%omega4
			%\node[black] at (6.2, -0.3) {\(\Omega_4\)};
			
			%mu-kappa
			%\fill (3.5,0) circle (1.5pt);
			%\draw[->, thick, decorate, decoration={snake, amplitude=.4mm, segment length=2mm, post length=1mm}] (2.9, -0.8) -- (3.4,-0.1);
			%\node[black] at (2.8, -1) {\(\frac{\mu}{2}-\kappa\)};
			%mu+kappa
			%\fill (2.5,0) circle (1.5pt);
			%\draw[->, thick, decorate, decoration={snake, amplitude=.4mm, segment length=2mm, post length=1mm}] (2.7, 0.8) -- (2.6,0.1);
			%\node[black] at (2.8, 1) {\(\frac{\mu}{2}+\kappa\)};
			
			%omega13
			%\node[black] at (-1,2.3) {\(\Omega_{13}\)};
			
			% origin point O
			\fill (0, 0) circle (1.5pt);
			\node[black] at (0.3, -0.3) {\(\text{0}\)};			
		\end{tikzpicture}
		\caption{The jump contour $\Sigma_2$ of the $\overline{\partial}$-RHP~\ref{RHP3} consists of the segments $\Sigma_k^\kappa$, $\Sigma_k^A$, $\Sigma_k^B$, and $\Sigma_k^C$ for $k = 1, \dots, 4$, together with $\Sigma_k^D$ for $k = 1, \dots, 6$. The exponential factor $e^{it\theta(z)}$ decays along the red contours, whereas $e^{-it\theta(z)}$ decays along the blue contours.}\label{F4}
	\end{figure}

	We now introduce a new unknown matrix 
	\begin{align}\label{eq63}
		\begin{aligned} M^{(2)}(z;x,t):=M^{(1)}(z)R^{(2)}(z)	.
		\end{aligned}
	\end{align}	
	Under the boundary conditions specified in \eqref{eq47}–\eqref{eq58}, the original jump contour $\Sigma_1$ of $M^{(1)}(z)$ is deformed into the modified contour $ \Sigma_2 $ for $ M^{(2)}(z) $, as shown in Figure \ref{F4}. By the definition of $ \chi_{\mathcal{Z}}(z) $, we have $ R^{(2)}(z) \equiv I $ within the disks $D(z_k,\kappa)$ for $ z_k\in \mathcal{Z}\cup \overline{\mathcal{Z}}$. Consequently, $ M^{(2)}(z) $ retains the same residue structure as $ M^{(1)}(z) $. 
	
	However, as a consequence of the contour deformation described above, the matrix $ M^{(2)}(z) $ is only continuous (and not analytic) on the domain  
	$$\Omega_c:=\mathbb{C} \backslash \left( \Sigma_2 \cup \Omega_{13}  \cup_{z_k \in \mathcal{Z} \cup \overline{\mathcal{Z}}} D(z_k, \kappa) \right).$$
	This lack of analyticity necessitates that \( M^{(2)}(z) \) satisfies a mixed \(\overline{\partial}\)-Riemann-Hilbert problem, where the classical jump conditions of a Riemann-Hilbert problem are coupled with a \(\overline{\partial}\)-equation governing the non-analytic contributions.  
	\setcounter{Drhp}{2}
	\begin{Drhp}\label{RHP3} For given scattering data $\sigma$, seek for a continuous matrix-valued function $M^{(2)}(z)$ defined on  $z\in \mathbb{C}\backslash\{\Sigma_2\}$ that satisfies the following conditions.
		\begin{itemize}
			\item $M^{(2)}(z)$ has continuous boundary values $M^{(2)}_\pm(z)$ on $z\in \Sigma_2$, which satisfy the jump conditions
			\begin{align}\label{eq64}
				\begin{aligned}
					M^{(2)}_+(z)=M^{(2)}_-(z)V^{(2)}(z;x,t),
				\end{aligned}
			\end{align}
			where
			\begin{align}\label{eq65}
				V^{(2)}(z) := 
				\begin{cases}
					\begin{aligned}
						& \begin{pmatrix} 1 & R_1 e^{2it \theta(z)} \\ 0 & 1 \end{pmatrix}
						&&\hspace{-5em}  z \in \Sigma_1^A
						&& \hspace{-5em} \begin{pmatrix} 1 & 0 \\ R_2 e^{-2it \theta(z)} & 1 \end{pmatrix}
						&& \hspace{-0em}z \in \Sigma_2^A \cup \Sigma_1^D\\[6pt]
						& \begin{pmatrix} 1 & R_3 e^{2it \theta(z)} \\ 0 & 1 \end{pmatrix}
						&& \hspace{-5em}z \in \Sigma_3^A \cup \Sigma_2^D
						&& \hspace{-5em}\begin{pmatrix} 1 & 0 \\ R_4 e^{-2it \theta(z)} & 1 \end{pmatrix}
						&& \hspace{-0em}z \in \Sigma_4^A\\[6pt]
						& \begin{pmatrix} 1 & 0 \\ R_5 e^{-2it \theta(z)} & 1 \end{pmatrix} 
						&& \hspace{-5em}z \in \Sigma_1^B \cup \Sigma_4^D
						&& \hspace{-5em}\begin{pmatrix} 1 & R_6 e^{2it \theta(z)} \\ 0 & 1 \end{pmatrix} 
						&& \hspace{-0em}z \in \Sigma_2^B\\[6pt]
						& \begin{pmatrix} 1 & 0 \\ R_7 e^{-2it \theta(z)} & 1 \end{pmatrix} 
						&& \hspace{-5em}z \in \Sigma_3^B
						&& \hspace{-5em}\begin{pmatrix} 1 & R_8 e^{2it \theta(z)} \\ 0 & 1 \end{pmatrix} 
						&&\hspace{-0em} z \in \Sigma_4^B \cup \Sigma_3^D\\[6pt]
						& \begin{pmatrix} 1 & 0 \\ R_9 e^{-2it \theta(z)} & 1 \end{pmatrix} 
						&& \hspace{-5em}z \in \Sigma_1^C
						&& \hspace{-5em}\begin{pmatrix} 1 & R_{10} e^{2it \theta(z)} \\ 0 & 1 \end{pmatrix} 
						&& \hspace{-0em}z \in \Sigma_2^C\\[6pt]
						& \begin{pmatrix} 1 & R_{11} e^{2it \theta(z)} \\ 0 & 1 \end{pmatrix} 
						&& \hspace{-5em}z \in \Sigma_3^C 
						&&\hspace{-5em} \begin{pmatrix} 1 & 0 \\ R_{12} e^{-2it \theta(z)} & 1 \end{pmatrix}
						&&\hspace{0em} z \in \Sigma_4^C \\[6pt]
						& \begin{pmatrix} 1 & 0 \\ (R_{12} - R_9) e^{-2it \theta(z)} & 1 \end{pmatrix} 
						&& z \in\Sigma_6^D \\[6pt]
						& \begin{pmatrix} 1 & (R_{11} - R_{10}) e^{2it \theta(z)} \\ 0 & 1 \end{pmatrix} 
						&& z \in\Sigma_5^D \\[6pt]
						&\begin{pmatrix}1 & 0 \\r( \RRe z) F^{-2}(z) e^{ -2it\theta(z) } & 1\end{pmatrix} 
						&&z\in\Sigma^\kappa_1\cup\Sigma^\kappa_4 \\[6pt]
						& \begin{pmatrix}1 & \overline{r}( \RRe z) F^{2}( z) e^{ 2it\theta(z) } \\0 & 1\end{pmatrix}
						&&z\in\Sigma^\kappa_2\cup\Sigma^\kappa_3.
					\end{aligned}
				\end{cases}
			\end{align}
			Moreover, $M^{(2)}(z)$ satisfies the following residue conditions at each simple poles $z_k$ and $\overline z_k$
			\begin{align}\label{eq66}
				\begin{aligned}
					\Res \limits _{z=z_{k}}M^{(2)}( z) =
					\begin{cases}
						\lim\limits_{z\rightarrow z_{k}}M^{(2)}( z) 
						\begin{pmatrix}
							0 & c_{k}^{-1}[ ( F^{-1}) '( z_{k})] ^{-2}e^{2it\theta( z_{k}) } \\
							0 & 0
						\end{pmatrix}
						&k\in \bigtriangleup,\\
						\lim\limits_{z\rightarrow z_{k}}M^{(2)}( z) 
						\begin{pmatrix}
							0 & 0 \\
							c_{k} F^{-2}( z_{k})e^{-2it\theta( z_{k}) } & 0
						\end{pmatrix}
						&k\in \bigtriangledown,\end{cases}\\
					\Res \limits_{z= \overline z_{k}}M^{(2)}(z)= \begin{cases} 
						\lim\limits_{z\rightarrow \overline z_{k}}M ^{(2)}( z) 
						\begin{pmatrix}
							0 & 0 \\
							-\overline c_{k}^{-1}( F'(\overline{z}_{k})) ^{-2}e^{ -2it\theta (\overline{z}_{k}) } & 0
						\end{pmatrix}
						&k\in \bigtriangleup,\\
						\lim\limits_{z\rightarrow \overline z_{k}}M ^{(2)}( z) 
						\begin{pmatrix}
							0 & -\overline c_{k}  F^{2}( \overline{z}_{k})e^{ 2it\theta (\overline{z}_{k}) } \\
							0 & 0 
						\end{pmatrix}
						&k\in \bigtriangledown,\end{cases}
				\end{aligned}
			\end{align}
			\item  For $z\in \Omega_c$, the continuous matrix function 
			$M^{(2)}(z)$ satisfies the $\overline \partial$-equations 
			\begin{align}\label{eq67}
				\overline \partial M^{(2)}(z)=M^{(2)}(z)\overline \partial R^{(2)}(z),
			\end{align}
			where
			\begin{align}\label{eq68}
				\overline \partial R^{(2)}(z):=
				\begin{aligned} 
					\begin{cases}
						\begin{pmatrix}
							0 &- \overline \partial R_{k}e^{2it\theta(z)} \\
							0 & 0
						\end{pmatrix} &z\in \Omega_k,{~~for ~}k=1,6\\
						\begin{pmatrix}
							0 & \overline \partial R_{k}e^{2it\theta(z) } \\
							0 & 0
						\end{pmatrix}  &z\in \Omega_k,\text{~~for ~}k=3,8,10,11\\
						\begin{pmatrix}
							0 & 0 \\
							-\overline \partial R_{k}e^{-2it\theta(z)} & 0
						\end{pmatrix}  &z\in \Omega_k,\text{~~for ~}k=2,5,9,12\\
						\begin{pmatrix}
							0 & 0 \\
							\overline \partial R_{k}e^{-2it\theta(z)} & 0
						\end{pmatrix}  &z\in \Omega_k,\text{~~for ~}k=4,7.
					\end{cases}
				\end{aligned}
			\end{align}
			\item  The matrix  $M^{(2)}(z)$ possesses the asymptotic behavior $M^{(2)}(z)\to I$, as $|z|\to \infty$.
		\end{itemize}
	\end{Drhp}
	\subsection{Decoupling of the $\overline{\partial}$-Riemann-Hilbert problem}
	Here, we express the solution to the $\overline{\partial}$-RHP \ref{RHP3} as  
	\begin{align}\label{eq69}  
		M^{(2)}(z) := M^{(3)}(z;x,t)M^{(2)}_{\mathrm{RHP}}(z;x,t),  
	\end{align}  
	which decomposes the complex $\overline{\partial}$-RHP into the following two subproblems that are ideally conditions for asymptotic analysis as $t \to \infty $. This decoupling allows $ M^{(3)}(z) $ to absorb the non-analytic effects, leaving $ M^{(2)}_{\mathrm{RHP}}(z) $ as a classical RHP with decaying jump matrix.  
	\begin{itemize}
		\item [1.] Meromorphic component: 
		\setcounter{rhp}{3}
		\begin{rhp}\label{RHP4}
			For given scattering data $\sigma$, the matrix $ M^{(2)}_{\mathrm{RHP}}(z) $ encodes the meromorphic structure of $ M^{(2)}(z) $, which satisfies the condition $\overline{\partial}M^{(2)}_{\mathrm{RHP}}(z) =0$ for $ z\in \mathbb{C}\backslash \Sigma_{2}$.
			\begin{itemize}
				\item [ $\bullet$]	$ M^{(2)}_{\mathrm{RHP}}(z) $,  inheriting the jump conditions from \eqref{eq64},  satisfies
				\begin{align}\label{eq70}  
					M^{(2)}_{\mathrm{RHP},+}(z)=M^{(2)}_{\mathrm{RHP},-}(z)V^{(2)}(z), z\in   \Sigma_{2}.
				\end{align}  
				At each simple pole $z_k$ and $\overline z_k$, the function $M_\mathrm{RHP}^{(2)}( z)$ inherits the residue relations specified in \eqref{eq66} by replacing  $M^{(2)}(z)$ in those conditions.
				\item [ $\bullet$]	 The matrix  $M^{(2)}_{\mathrm{RHP}}(z)$ possesses the asymptotic behavior $M^{(2)}_{\mathrm{RHP}}(z)\to I$, as $|z|\to \infty$.
			\end{itemize}
		\end{rhp}
		
		\item [2.] Non-analytic component:  \\
		If $ M^{(2)}_{\mathrm{RHP}}(z) $ satisfies the above description, then $ M^{(3)}(z)$ is determined by a pure $\overline{\partial}$-problem \ref{RHP5}, below, where the non-analytic contributions (arising from $ R^{(2)}(z) $) are encoded in the \(\overline{\partial}\)-equation \eqref{eq67}. 
		\setcounter{D}{4}\begin{D}\label{RHP5}
			The matrix $ M^{(3)}(z) $, encoding the non-analytic structure of $M^{(2)}(z)$, is globally continuous across the complex plane and $\overline \partial$-differentiable on $\mathbb{C}\backslash(\mathbb{R}\cup\Sigma_{2})$. Specifically,
			$M^{(3)}(z)$ satisfies the following conditions.
			\begin{itemize}
				\item [ $\bullet$] $\overline \partial $-equation: For $z\in \cup_{k=1}^{12}\Omega_k$,
				\begin{align}\label{eq72}
					\begin{aligned}
						\overline \partial M^{(3)}(z)=M^{(3)}(z)W(z;x,t), 
					\end{aligned}
				\end{align}
				where  $W(z):=M^{(2)}_{\mathrm{RHP}}(z)\overline \partial R^{(2)}(z)(M^{(2)}_{\mathrm{RHP}}(z))^{-1}$ and $\overline \partial R^{(2)}(z)$ is defined in \eqref{eq68}.
				\item [ $\bullet$]  The matrix  $M^{(3)}(z)$ possesses the asymptotic behavior $M^{(3)}(z)\to I$, as $|z|\to \infty$.
			\end{itemize}
		\end{D}
		
	\end{itemize}
	\begin{proof}
		By adapting the analytical framework from \cite{MR3795020} with slight modifications, we can readily prove that $ M^{(3)}(z) $  is jump-free on $ \Sigma_2 $ and the simple poles $ z_k \in \mathcal{Z} \cup \overline{\mathcal{Z}} $, originally present in $ M^{(2)}(z) $, are reduced to removable singularities in $ M^{(3)}(z) $. In addition, for $z\in \cup_{k=1}^{12}\Omega_k$, combining \eqref{eq67} and \eqref{eq69}, we derive
		\begin{align*}
			\overline \partial M^{(3)}(z)&=\overline \partial M^{(2)}(z)(M^{(2)}_{\mathrm{RHP}}(z)) ^{-1}\\
			&=M^{(2)}(z)\overline \partial R^{(2)}(z)(M^{(2)}_{\mathrm{RHP}}(z)) ^{-1} \\
			&=M^{(3)}(z)M^{(2)}_{\mathrm{RHP}}(z)\overline \partial R^{(2)}(z)(M^{(2)}_{\mathrm{RHP}}(z)) ^{-1} ,
		\end{align*}
		which is consistent with \eqref{eq72}.  
	\end{proof}
	
	\section{Resolution of the pure Riemann--Hilbert and pure $\overline{\partial}$ problems}\label{S4}
	%\section{RHP 4的分解及其近似}
	In this section, we further decompose RHP \ref{RHP4} based on the decay behavior of jump matrix $V^{(2)}(z)$, as detailed in Lemma \ref{Lemma4}, and analyze its asymptotic behavior as $ t \to \infty $.  The distinct contributions of the jump and residue conditions  to the asymptotic expansion of $ M^{(2)}_{\mathrm{RHP}}(z) $ are separated using the approximately solvable models within specific regions of the $ z $-complex plane.  
	
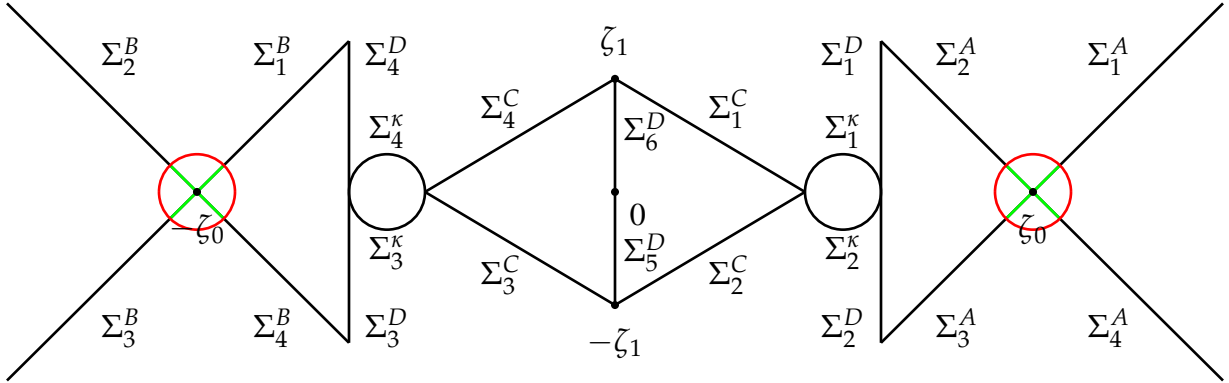
\begin{figure}[h] 
	\centering 
	\begin{tikzpicture}
		% RRe z 
		%\draw[black, thick, dashed] (-8,0) -- (-3.5,0);
		%\draw[black, thick, dashed] (-2.5,0) -- (2.5,0);
		%\draw[black, thick, dashed] (3.5,0) -- (8,0) node[right] {\text{Re} $z$};
		%zeta1
		\fill (0,1.5) circle (1.5pt);
		\node[black] at (0, 2) {\(\zeta_1\)};
		%-zeta1
		\fill (0,-1.5) circle (1.5pt);
		\node[black] at (0, -2) {\(-\zeta_1\)};
		%left
		%sigma4c
		\draw[black, line width=1pt] (-2.5,0) -- (0,1.5);
		%\draw[black, -{Stealth},line width=1pt] (-2.5,0) -- (-1,0.9);
		\node[black] at (-1.5,1.1) {\(\Sigma_4^C\)};
		%omega12
		%\node[black] at (-1.3, 0.3) {\(\Omega_{12}\)};
		
		%sigma3c
		\draw[black, line width=1pt] (-2.5,0) -- (0,-1.5);
		%\draw[black,-{Stealth}, line width=1pt] (-2.5,0) -- (-1, -0.9);
		\node[black] at (-1.5,-1.1) {\(\Sigma_3^C\)};
		%omega11
		%\node[black] at (-1.3, -0.3) {\(\Omega_{11}\)};
		
		%sigma6D
		\draw[black, line width=1pt] (0,1.5) -- (0,0);
		%\draw[black, -{Stealth},line width=1pt] (0,1.5) -- (0,0.5);
		\node[black] at (0.4, 0.8) {\(\Sigma_6^D\)};
		
		%sigma5D
		\draw[black, line width=1pt] (0,-1.5) -- (0,0);
		%\draw[black, -{Stealth},line width=1pt] (0,-1.5) -- (0,-0.5);
		\node[black] at (0.4, -0.8) {\(\Sigma_5^D\)};
		
		%sigma1C
		\draw[black, line width=1pt] (2.5,0) -- (0,1.5);
		%\draw[black, -{Stealth}, line width=1pt](0,1.5) -- (1.5,0.6);
		\node[black] at (1.5,1.1) {\(\Sigma_1^C\)};
		%omega9
		%\node[black] at (1.2, 0.3) {\(\Omega_{9}\)};
		
		%sigma2C
		\draw[black, line width=1pt] (2.5,0) -- (0,-1.5);
		%\draw[black,-{Stealth}, line width=1pt] (0,-1.5) -- (1.5, -0.6);
		\node[black] at (1.5,-1.1) {\(\Sigma_2^C\)};
		%omega10
		%\node[black] at (1.3, -0.3) {\(\Omega_{10}\)};
		
		%sigma4D
		\draw[black, line width=1pt] (-3.5,2) -- (-3.5,0);
		%\draw[black, -{Stealth},line width=1pt] (-3.5,2) -- (-3.5,1);
		\node[black] at (-3,1.8) {\(\Sigma_4^D\)};
		
		%sigma3D
		\draw[black, line width=1pt] (-3.5, -2) -- (-3.5,0);
		%\draw[black, -{Stealth},line width=1pt] (-3.5, -2) -- (-3.5,-1);
		\node[black] at (-3, -1.8) {\(\Sigma_3^D\)};
		
		%sigmaB
		\draw[black, line width=1pt] (-8,-2.5) -- (-3.5,2);
		\draw[black, line width=1pt] (-8,2.5) -- (-3.5,-2);
		%sigma3B
		%\draw[black, -{Stealth}, line width=1pt] (-8,-2.5) -- (-6.0,-0.5);
		\node[black] at (-6.5,-1.8) {\(\Sigma_3^B\)};
		%sigm1B
		%\draw[black, -{Stealth}, line width=1pt] (-5.5,0) -- (-4.75,0.75);
		\node[black] at (-4.5,1.8) {\(\Sigma_1^B\)};
		%sigma4B
		%\draw[black, -{Stealth}, line width=1pt] (-5.5, 0) -- (-4.75, -0.75);
		\node[black] at (-4.5,-1.8) {\(\Sigma_4^B\)};
		%sigm4B
		%\draw[black, thick,black, -{Stealth}, line width=1pt] (-8,2.5) -- (-6.0,0.5);
		\node[black] at (-6.5,1.8) {\(\Sigma_2^B\)};
		
		%right
		%sigma1D
		\draw[black, line width=1pt] (3.5,2) -- (3.5,0);
		%\draw[black, -{Stealth},line width=1pt] (3.5,0) -- (3.5,1);
		\node[black] at (3,1.8) {\(\Sigma_1^D\)};
		
		%sigma2D
		\draw[black, line width=1pt] (3.5, -2) -- (3.5,0);
		%\draw[black, -{Stealth},line width=1pt] (3.5, 0) -- (3.5,-1);
		\node[black] at (3, -1.8) {\(\Sigma_2^D\)};
		
		%sigmaA
		\draw[black, line width=1pt] (8,-2.5) -- (3.5,2);
		\draw[black, line width=1pt] (8,2.5) -- (3.5,-2);
		%sigma4A
		%\draw[black, -{Stealth}, line width=1pt] (5.5,0) -- (6.25,-0.75);
		\node[black] at (6.5,-1.8) {\(\Sigma_4^A\)};
		%sigm2A
		%\draw[black, -{Stealth}, line width=1pt] (3.5,2) -- (5,0.5);
		\node[black] at (4.5,1.8) {\(\Sigma_2^A\)};
		%sigma1A
		%\draw[black, -{Stealth}, line width=1pt] (5.5, 0) -- (6.25, 0.75);
		\node[black] at (6.5,1.8) {\(\Sigma_1^A\)};
		%sigm3A
		%\draw[black, thick,black, -{Stealth}, line width=1pt] (3.5,-2) -- (5,-0.5);
		\node[black] at (4.5,-1.8) {\(\Sigma_3^A\)};
		
		\draw[black,line width=1pt] (2.5,0) arc (180:0:0.5);			
		%\draw[black, thick, -{Stealth}] (3,0.5) -- (3.1,0.5);			
		\node[black] at (3,0.8) {$\Sigma_1^\kappa$};			
		\draw[black,line width=1pt] (2.5,0) arc (-180:0:0.5);			
		%\draw[black, thick, -{Stealth}] (3,-0.5) -- (3.1,-0.5);			
		\node[black] at (3,-0.8) {$\Sigma_2^\kappa$};
		\draw[black,line width=1pt] (-3.5,0) arc (180:0:0.5);			
		%\draw[black, thick, -{Stealth}] (-3,0.5) -- (-2.9,0.5);		
		\node[black] at (-3,0.8) {$\Sigma_4^\kappa$};		
		\draw[black, line width=1pt] (-3.5,0) arc (-180:0:0.5);		
		%\draw[black, thick, -{Stealth}] (-3,-0.5) -- (-2.9,-0.5);			
		\node[black] at (-3,-0.8) {$\Sigma_3^\kappa$};
		\draw[red, line width=1pt] (-6,0) arc (180:0:0.5);
		\draw[red,line width=1pt] (-6,0) arc (-180:0:0.5);
		\draw[red, line width=1pt] (5,0) arc (180:0:0.5);
		\draw[red,line width=1pt] (5,0) arc (-180:0:0.5);
		
		\draw[green, line width=1pt] (-5.84,0.34) -- (-5.16,-0.34);
		\draw[green, line width=1pt] (-5.84,-0.34) -- (-5.16,0.34);
		\draw[green, line width=1pt] (5.84,0.34) -- (5.16,-0.34);
		\draw[green, line width=1pt] (5.84,-0.34) -- (5.16,0.34);

		%left%zeta0
		\fill (-5.5,0) circle (1.5pt);
		\node[black] at (-5.5, -0.5) {\(-\zeta_0\)};
		%omega5
		%\node[black] at (-4.8, 0.3) {\(\Omega_5\)};
		%omega8
		%\node[black] at (-4.8, -0.3) {\(\Omega_8\)};
		
		%omega6
		%\node[black] at (-6.2, 0.3) {\(\Omega_6\)};
		%omega7
		%\node[black] at (-6.2, -0.3) {\(\Omega_7\)};
		
		%-mu-kappa
		%\fill (-3.5,0) circle (1pt);
		%\draw[->, thick, decorate, decoration={snake, amplitude=.4mm, segment length=2mm, post length=1mm}] (-2.7, -0.8) -- (-3.4,-0.1);
		%\node[black] at (-2.8, -1) {\(-\frac{\mu}{2}-\kappa\)};
		%-mu+kappa
		%\fill (-2.5,0) circle (1pt);
		%\draw[->, thick, decorate, decoration={snake, amplitude=.4mm, segment length=2mm, post length=1mm}] (-2.7, 0.8) -- (-2.6,0.1);
		%\node[black] at (-2.8, 1) {\(-\frac{\mu}{2}+\kappa\)};
		
		%right
		%zeta0
		\fill (5.5,0) circle (1.5pt);
		\node[black] at (5.5, -0.5) {\(\zeta_0\)};
		
		%omega2
		%\node[black] at (4.8, 0.3) {\(\Omega_2\)};
		%omega3
		%\node[black] at (4.8, -0.3) {\(\Omega_3\)};
		
		%omega1
		%\node[black] at (6.2, 0.3) {\(\Omega_1\)};
		%omega4
		%\node[black] at (6.2, -0.3) {\(\Omega_4\)};
		
		%mu-kappa
		%\fill (3.5,0) circle (1.5pt);
		%\draw[->, thick, decorate, decoration={snake, amplitude=.4mm, segment length=2mm, post length=1mm}] (2.9, -0.8) -- (3.4,-0.1);
		%\node[black] at (2.8, -1) {\(\frac{\mu}{2}-\kappa\)};
		%mu+kappa
		%\fill (2.5,0) circle (1.5pt);
		%\draw[->, thick, decorate, decoration={snake, amplitude=.4mm, segment length=2mm, post length=1mm}] (2.7, 0.8) -- (2.6,0.1);
		%\node[black] at (2.8, 1) {\(\frac{\mu}{2}+\kappa\)};
		
		%omega13
		%\node[black] at (-1,2.3) {\(\Omega_{13}\)};
		
		% origin point O
		\fill (0, 0) circle (1.5pt);
		\node[black] at (0.3, -0.3) {\(\text{0}\)};			
	\end{tikzpicture}
		\caption{As $t\to \infty$, uniform decay of $V^{(2)}(z) $ to the identity matrix on black contours, and point-wise convergence on green contours.  The red boundaries are $\partial D(\pm \zeta_0,\kappa).$}\label{F5}
	\end{figure}

	As shown in Figure \ref{F5}, outside the disks $D(\pm\zeta_0, \kappa) $, the jump matrix $ V^{(2)}(z) $ defined in \eqref{eq65} decays uniformly and exponentially to the identity matrix as $ t \to \infty $. In this region,  the residue conditions (encoding solitons), provide the dominant contribution, while the jump conditions become negligible. 
	
	Conversely, within the disks $ D(\pm\zeta_0, \kappa) $, the jump matrix $ V^{(2)}(z) $ approaches the identity point-wise. Choosing $ \kappa > 0 $ sufficiently small to ensure $ D(\pm\zeta_0, \kappa) $ is disjoint from $ \cup_{z_k \in \mathcal{Z} \cup \overline{\mathcal{Z}}} D(z_k, \kappa) $  allows the analysis to be focused entirely on the jump conditions in \eqref{eq70}, due to the absence of discrete spectrum (poles) in $ D(\pm\zeta_0, \kappa) $.  
	
	Based on the above analysis, we aim to construct an appropriate form for $M_{\mathrm{RHP}}^{(2)}(z) $ by introducing a piecewise definition
	\begin{align}\label{eq73}
		\begin{aligned}
			M_{\mathrm{RHP}}^{(2)}(z) =
			\begin{cases}
				E_{\circ}(z;x,t) M_{\mathrm{out}}(z;x,t), & z \in \mathbb{C} \backslash D(\pm \zeta_0, \kappa), \\
				E_{\circ}(z;x,t) M_{\pm \zeta_0}(z;x,t), & z \in D(\pm \zeta_0, \kappa),
			\end{cases}
		\end{aligned}
	\end{align}  
	where the matrix $ E_{\circ}(z;x,t) $ serves as an error correction term. This matrix accounts for the discrepancy between the actual solution $ M_{\mathrm{RHP}}^{(2)}(z) $ and the model matrices $ M_{\mathrm{out}}(z;x,t) $ and $ M_{\pm \zeta_0}(z;x,t) $.  
	
	\subsection{Out model}
	\begin{Lemma}\label{Lemma4}
		Let the complex variable be expressed as $z=a+ib$. For $x/t\in (-1/\mu^2,0)$ with a fixed $\mu\in (0,1]$, the jump matrix $V^{(2)}(z)$ exhibits exponential decay on the contour $\Sigma_{2}$ away from the stationary phase points $\pm\zeta_0$.
		\begin{align}\label{eq74}
			\begin{aligned}
				\left\| V^{(2)}-I\right\| _{L^{\infty }}=\begin{cases}\mathcal{O}( {e}^{c_1|a-\zeta_0|^2t}) & z\in \cup^{4}_{k=1}\Sigma _{k}^{A},\\
					\mathcal{O}(e^{c_1|a+\zeta_0|^2t})  & z\in \cup^{4}_{k=1}\Sigma _{k}^{B},\\
					\mathcal{O}( e^{-c_2t}) &\text{otherwise~on~} \Sigma_2,\end{cases}
			\end{aligned}
		\end{align}
		where $c_1<0$ is defined in \eqref{eq74.1} and $c_2>0.$
	\end{Lemma}
	\begin{proof}
		By restricting the reformulated expression of \eqref{eq32} 
		\begin{align*}
			\RRe [i\theta(z)] =b\left( \frac{ 16( a^{2}-\zeta_0^{2})( a^{2}-\zeta _{1}^{2}){x}/{t} +4b^{2}(  1+2( 4a^{2}+2b^{2}+\mu ^{2}){x}/{t} ) }{  ( ( 2a+\mu ) ^{2}+4b^{2})  ( (2a-\mu )^{2}+4b^{2}) }\right) 
		\end{align*}
		to the contour $ \Sigma_{1}^A := \{ z = a + \mathrm{i}b : b = a - \zeta_0, \, a \in (\zeta_0, \infty) \} $, we derive the asymptotic expansion
		\begin{align*}
			\RRe [i\theta(z)] =\frac{32\zeta_{0}( \zeta_0^{2}-\zeta_{1}^{2}) {x}/{t}}{( 2\zeta_{0}+\mu ) ^{2}( 2\zeta_{0}-\mu ) ^{2}}( a-\zeta_{0}) ^{2}+\mathcal{O}( a-\zeta_{0}) ^{3}
		\end{align*}
		which defines the decay constant  
		\begin{align}\label{eq74.1}
			c_1=\frac{64\zeta_{0}( \zeta_0^{2}-\zeta_{1}^{2}) {x}/{t}}{( 2\zeta_{0}+\mu ) ^{2}( 2\zeta_{0}-\mu ) ^{2}}.
		\end{align}
		
		For  $z\in \Sigma^\kappa_{k}$, setting $z:=\pm {\mu }/{2}+\kappa e^{i \alpha}$,  we evaluate the limits  
		\begin{align*}
			\begin{aligned}
				\lim_{\kappa \to 0}\RRe( i\theta ( \pm {\mu }/{2}+\kappa e^{i \alpha} ) ) &=\lim_{\kappa \to 0}\kappa \sin \alpha \left(\frac{4((\kappa\cos \alpha \pm \mu /2)^2+\kappa^2\sin^2 \alpha)+\mu^2}{16\kappa^2((\kappa \cos \alpha \pm \mu)^2+\kappa^2\sin^2\alpha)}+\frac{x}{t}\right)\\
				&=\lim_{\kappa \to 0}\frac{\sin \alpha}{4(\mu^2+\kappa^2\pm2\mu \kappa\cos \alpha)}\left( \kappa \pm\mu \cos \alpha+ \frac{\mu^2}{\kappa}\right)+\frac{x}{t}\kappa\sin \alpha.
			\end{aligned}
		\end{align*}
		For $k=1,4$, where $\alpha \in (0,\pi)$, we have 
		\begin{align*}
			\begin{aligned}
				\lim_{\kappa \to 0}\RRe(- i\theta ( \pm {\mu }/{2}+\kappa e^{i \alpha} ) ) 
				=\lim_{\kappa \to 0}\frac{-\sin \alpha}{4\mu^2}\left(  \pm\mu \cos \alpha+ \frac{\mu^2}{\kappa}\right)=-\infty.
			\end{aligned}
		\end{align*}
		Similarly, for $k=2,3$, with $\alpha \in (-\pi,0)$,
		\begin{align*}
			\begin{aligned}
				\lim_{\kappa \to 0}\RRe( i\theta ( \pm {\mu }/{2}+\kappa e^{i \alpha} ) ) =-\infty.
			\end{aligned}
		\end{align*}
		The continuity of $\RRe[i\theta(z)]$ implies the existence of a constant $c_2 > 0 $ ensuring the uniform validity of the exponential decay estimates in \eqref{eq74} for $z \in \cup_{k=1}^4 \Sigma_k^\kappa $. Analogously, for $ z \in \cup_{k=1}^4 \Sigma_k^C \cup_{j=1}^6 \Sigma_j^D $, another constant $c_2 > 0 $ guarantees the validity of the decay rates in \eqref{eq74}.
	\end{proof}
	
	The conclusion of lemma \ref{Lemma4} provides sufficient basis for decomposition \eqref{eq73}. We now formulate RHP \ref{RHP6} governing the solvable model $M_{\mathrm{out}}(z).$
	\setcounter{rhp}{5}
	\begin{rhp}\label{RHP6}
		For given scattering data $\sigma:=\{r(z),\{(z_k,c_k)\}_{k=1}^N\}$,  seek for a meromorphic matrix-valued function $M_{\mathrm{out}}(z;\sigma)$ defined on  $z\in \mathbb{C}$ that satisfies the following conditions.
		\begin{itemize}
			\item At each simple pole $z_k$ and $\overline z_k$, the function $M_\mathrm{out}( z;\sigma)$ inherits the residue relations specified in \eqref{eq66} by replacing  $M^{(2)}(z;\sigma)$ in these conditions.
			\item The matrix  $M_{\mathrm{out}}(z;\sigma)$ possesses the asymptotic behavior $M_{\mathrm{out}}(z;\sigma)\to I$, as $|z|\to \infty$.
		\end{itemize}
	\end{rhp}
	\begin{prop}\label{prop5}
		There exist a unique solution $M_{\mathrm{out}}(z;\sigma)$ of RHP \ref{RHP6}. Moreover, this solution generates the $N$-soliton or $N$-kink solutions through the reconstruction formulae
		\begin{align}\label{eq75}
			\begin{aligned} 
				&E(x,t;\sigma_{\delta}) =-4i\lim _{z\rightarrow \infty }\left( zM^{\bigtriangleup }(z;\sigma_\delta)\right)_{12},\quad&s(x,t;\sigma_{\delta}) =-\frac{1}{2}\left(\rho( \pm\frac{\mu }{2};x,t) _{12}+\rho( \pm\frac{\mu }{2};x,t) _{21}\right),\\
				&u(x,t;\sigma_{\delta}) =-\rho( \pm\frac{\mu }{2};x,t) _{11},\quad&r(x,t;\sigma_{\delta}) =\mp \frac{1}{2i}\left( \rho( \pm\frac{\mu }{2};x,t) _{12}-\rho( \pm\frac{\mu }{2};x,t)_{21}\right), 
			\end{aligned}
		\end{align}
		where $\rho( z;x,t) =M^{\bigtriangleup }(z;\sigma_\delta) \sigma _{3}M^{\bigtriangleup }(z;\sigma_\delta) $. Here, $M^{\bigtriangleup }(z;\sigma_\delta) $ solves RHP \ref{RHP10} with the  modified scattering data 
		\begin{align}\label{eq75.1}
			\sigma_\delta:=\{\{(z_k,c_{\delta,k}=c_k\delta^{-2}(z_k))\}_{k=1}^N\},
		\end{align} where   $\delta(z)$ is given in \eqref{eq38}.
	\end{prop}
	\begin{Remark}
		The solutions generated by $M_{\mathrm{out}}(z;\sigma)$ via \eqref{eq75} are all localized traveling waves. Specifically, $E(x,t;\sigma_\delta)$, $u(x,t;\sigma_\delta)$, and $r(x,t;\sigma_\delta)$ represent localized $N$-soliton solutions, while  $s(x,t;\sigma_\delta)$ corresponds to an $N$-kink solution.
	\end{Remark}
	\begin{proof}
		The proof of Proposition \ref{prop5} proceeds in two key stages.
		
		\textbf{Step 1.}  
		We first prove that the reflectionless case of RHP \ref{RHP1} admits a unique solution $M(z;\sigma_0)$, where $\sigma_0:=\{r(z)\equiv 0,\{(z_k,c_k)\}_{k=1}^N\}$ represents the discrete scattering data with vanishing reflection coefficient. 
		
		Under the condition $r(z)\equiv 0$, the jump matrix simplifies to $V(z)\equiv I$ in \eqref{eq9}, thereby reducing RHP \ref{RHP1} to a residue problem. Let the solution matrix be expressed element-wise as
		\begin{align*}
			M(z;\sigma _{0}) :=\begin{pmatrix}
				m_{11}( z;x,t)  & m_{12}( z;x,t)   \\
				m_{21}( z;x,t)  & m_{22}( z;x,t)  
			\end{pmatrix}.
		\end{align*}
		Combining the symmetry relation $M(z;\sigma_0)=\sigma_{2} \overline M(\overline z;\sigma_0) \sigma_{2}$ with  the residue condition \eqref{eq10}, induces the partial fraction decomposition
		\begin{align*}
			M(z;\sigma _{0}) =I+\sum ^{N}_{k=1}\frac{1}{z-z_{k}}\begin{pmatrix}
				m_{12}( z_{k}) c_k e^{-2i\varphi ( z_{k}) } & 0\\
				m_{22}( z_{k}) c_k e^{-2i\varphi ( z_{k}) } & 0
			\end{pmatrix}+\sum ^{N}_{k=1}\frac{1}{z-\overline{z}_{k}}\begin{pmatrix}
				0 & -\overline{m_{22}}( z_{k}) \overline{c}_k e^{2i\overline \varphi ( z_{k}) } \\
				0 & \overline{m_{12}}( z_{k}) \overline{c}_k e^{2i\overline \varphi ( z_{k})  }
			\end{pmatrix}.
		\end{align*}
		Substituting this expansion into \eqref{eq10}, we obtain the linear system for the unknowns $m_{12}(z_n)$ and $m_{22}(z_n)$ (for  $n=1,\ldots ,N$)
		\begin{align*}
			\begin{pmatrix}
				m_{12}( z_{n}) c_{n}e^{-2i\varphi ( z_{n}) } & 0 \\
				m_{22}( z_{n}) c_{n}e^{-2i\varphi ( z_{n}) } & 0
			\end{pmatrix}=\begin{pmatrix}
				0 & 0 \\
				c_k e^{-2i\varphi ( z_{k}) } & 0
			\end{pmatrix}+\sum ^{N}_{k=1}\frac{1}{z_{n}-\overline{z}_{k}}\begin{pmatrix}
				-\overline{m_{22}}( z_{k}) | c_k| ^{2}e^{-4i{\varphi}( z_{k}) } & 0 \\
				\overline{m_{12}}( z_{k}) | c_k| ^{2}e^{-4i{\varphi}( z_{k}) } & 0
			\end{pmatrix},
		\end{align*}
		where the relation $\overline{\varphi}( z_{k}) =-{\varphi}( z_{k}) $ arises from $z_k \in i \mathbb{R}_+$. 
		The uniqueness and existence of $M(z;\sigma_0)$ follow directly from unique solvability of  this linear system (see Proposition B.1. in \cite{MR3795020}) and Liouville’s theorem.
		
		Moreover, from Theorem \ref{theo2}, the matrix $M(z;\sigma_0)$ encodes the reflectionless solutions $E(x,t; \sigma_{0})$, $s(x,t; \sigma_{0})$, $u(x,t; \sigma_{0})$, and $r(x,t; \sigma_{0})$ by replacing $M(z;x,t)$ with $M(z;\sigma_0)$ in the reconstruction formulae.
		
		\textbf{Step 2.}  By using \eqref{eq76} and \eqref{eq77} below, we relate $M_{\mathrm{out}}( z;\sigma)$ to $M(z;\sigma_0)$, which guarantees the existence and uniqueness of $M_{\mathrm{out}}( z;\sigma) $.
		
		Revisiting the residue condition \eqref{eq10}, we note that $z_k$ and  $\overline z_k$ (for $k=1,\ldots,N$) are simple poles of the first and second columns of  $M(z;\sigma_0)$ , respectively. Let $ \bigtriangleup$ be as defined in \eqref{eq37} and $\bigtriangledown=\{1,\ldots,N\}\backslash \bigtriangleup$. Define
		\begin{align*}
			F_{\bigtriangleup }(z):=\prod _{k\in\bigtriangleup}\frac{z-\overline z_k}{z-z_k}.
		\end{align*}
		Then, the transformation 
		\begin{align}\label{eq76}
			M^{\bigtriangleup }( z;\sigma _{0}) =M( z;\sigma _{0}) ( F_{\bigtriangleup }(z)) ^{-\sigma _{3}},
		\end{align}
		transfers poles $z_k$ ($k\in\bigtriangleup$) to the second column and $\overline{z}_k$ ($k\in\bigtriangleup$) to the first column of $M^{\bigtriangleup }( z;\sigma _{0}) $, while leaving others unchanged. Consequently, $M^{\bigtriangleup }( z;\sigma _{0}) $ satisfies RHP \ref{RHP10}.
		
		\begin{rhp}\label{RHP10}
			For given scattering data $\sigma_0$, seek for a meromorphic matrix-valued function $M^{\bigtriangleup }( z;\sigma _{0}) $ defined on  $z\in \mathbb{C}$ that satisfies the following conditions.
			\begin{itemize}
				\item	$M^{\bigtriangleup }( z;\sigma _{0}) $ satisfies the residue conditions at each simple poles $z_k$ and $\overline z_k$
				\begin{align}\label{eq76.1}
					\begin{aligned}
						\Res \limits _{z=z_{k}}M^{\bigtriangleup }( z;\sigma _{0})  &=
						\begin{cases}
							\lim\limits_{z\rightarrow z_{k}}M^{\bigtriangleup }( z;\sigma _{0}) 
							\begin{pmatrix}
								0 & c_{k}^{-1}[ ( F_{\bigtriangleup }^{-1}) '( z_{k})] ^{-2}e^{2it\theta( z_{k}) } \\
								0 & 0
							\end{pmatrix}
							&k\in \bigtriangleup\\
							\lim\limits_{z\rightarrow z_{k}}M^{\bigtriangleup }( z;\sigma _{0}) 
							\begin{pmatrix}
								0 & 0 \\
								c_{k} F_{\bigtriangleup }^{-2}( z_{k})e^{-2it\theta( z_{k}) } & 0
							\end{pmatrix}
							&k\in \bigtriangledown\end{cases}\\
						\Res \limits_{z= \overline z_{k}}M^{\bigtriangleup }( z;\sigma _{0}) &= \begin{cases} 
							\lim\limits_{z\rightarrow \overline z_{k}}M^{\bigtriangleup }( z;\sigma _{0}) 
							\begin{pmatrix}
								0 & 0 \\
								-\overline c_{k}^{-1}( F_{\bigtriangleup }'(\overline{z}_{k})) ^{-2}e^{ -2it\theta (\overline{z}_{k}) } & 0
							\end{pmatrix}
							&k\in \bigtriangleup\\
							\lim\limits_{z\rightarrow \overline z_{k}}M^{\bigtriangleup }( z;\sigma _{0}) 
							\begin{pmatrix}
								0 & -\overline c_{k}  F_{\bigtriangleup }^{2}( \overline{z}_{k})e^{ 2it\theta (\overline{z}_{k}) } \\
								0 & 0 
							\end{pmatrix}
							&k\in \bigtriangledown,\end{cases}
					\end{aligned}
				\end{align}
				\item  The matrix  $M^{\bigtriangleup }( z;\sigma _{0})$ possesses the asymptotic behavior $M^{\bigtriangleup }( z;\sigma _{0})\to I$, as $|z|\to \infty$.
			\end{itemize}
		\end{rhp}
		
		Furthermore, as $|z|\to \infty$ the asymptotic expansion of $F_{\bigtriangleup }(z)$ is given by 
		\begin{align*}
			F_{\bigtriangleup }(z)=1+\frac{1}{z}\sum_{k\in \bigtriangleup}2 z_k +\mathcal{O}(z^{-1}),
		\end{align*}
		leading to the expansion of 
		\begin{align*}
			\begin{aligned}
				M^{\bigtriangleup }(z;\sigma _{0}) 	=I+\frac{1}{z}\left( M_1(x ,t;\sigma_{0}) -2\sigma_3\sum _{k\in \bigtriangleup } z_k\right) +\mathcal{O}( z^{-1}),
			\end{aligned}
		\end{align*}
		where $M(z;\sigma_{0})=I+M_1(x,t;\sigma_0)/z+\mathcal{O} (z^{-1})$. Form \eqref{eq76}, we have 
		\begin{align*}
			\begin{aligned}
				\rho ( z;x,t) &=M ( z;\sigma _{0}) \sigma _{3}M ^{-1}( z;\sigma _{0}) \\
				&=M^{\bigtriangleup }( z;\sigma _{0}) F_{\bigtriangleup }(z) ^{\sigma _{3}}\sigma _{3}F_{\bigtriangleup }(z) ^{-\sigma _{3}}\left(M^{\bigtriangleup }( z;\sigma _{0}) \right) ^{-1}\\
				&=M^{\bigtriangleup }( z;\sigma _{0}) \sigma _{3}\left(M^{\bigtriangleup }( z;\sigma _{0}) \right) ^{-1}.
			\end{aligned}
		\end{align*}
		Then, the reflectionless solutions are recoverable from $M^{\bigtriangleup }( z;\sigma _{0}) $ through the reconstruction formulae presented in Theorem \ref{theo2}, by replacing  $M(z;x,t)$ with $M^{\bigtriangleup }(z;\sigma_0)$ 
			\begin{align*}
				\begin{aligned} 
					&E(x,t;\sigma_{0}) =-4i\lim _{z\rightarrow \infty }\left( zM^{\bigtriangleup }(z;\sigma_0)\right)_{12},\quad&s(x,t;\sigma_{0}) =-\frac{1}{2}\left(\rho( \pm\frac{\mu }{2};x,t) _{12}+\rho( \pm\frac{\mu }{2};x,t) _{21}\right),\\
					&u(x,t;\sigma_{0}) =-\rho( \pm\frac{\mu }{2};x,t) _{11},\quad&r(x,t;\sigma_{0}) =\mp \frac{1}{2i}\left( \rho( \pm\frac{\mu }{2};x,t) _{12}-\rho( \pm\frac{\mu }{2};x,t)_{21}\right), 
				\end{aligned}
			\end{align*}
			whare $\rho( z;x,t) =M^{\bigtriangleup }(z;\sigma_0) \sigma _{3}M^{\bigtriangleup }(z;\sigma_0) $.
		
		Comparing RHP \ref{RHP6} and RHP \ref{RHP10}, we find the relation
		\begin{align}\label{eq77}
			M_{\mathrm {out}}( z;\sigma ) =M^{\bigtriangleup }( z;\sigma _{\delta }) ,
		\end{align}
		where $\sigma_\delta$ is given in \eqref{eq75.1}.
	\end{proof}
	
	For $N=1$, assuming the scattering data $\sigma_{0}:=\{r(z)\equiv0,(i\eta,ic)\}$ with $\eta>0$, the exact solutions to the RMB equations \eqref{eq01} encoded by $M(z;\sigma_0)$ are given by 
	\begin{align*}%\label{eq78}
		\begin{aligned}
			E(x,t;\sigma_0) &= 4\eta \,\mathrm{sign}(c) \sech\Theta  \quad\quad	s(x,t;\sigma_0) = \frac{8\eta^2\cosh\Theta - 4\eta c e^{-\Lambda}}{D}&\\
			u(x,t;\sigma_0) &= -1 + \frac{8\eta^2}{D} \quad\quad\quad \quad\quad r(x,t;\sigma_0) = \frac{4\eta\mu\cosh\Theta}{D},
		\end{aligned}
	\end{align*}
	where
	$$
	\begin{aligned}
		D := (4\eta^2 + \mu^2)\cosh^2\Theta ,~~	\Theta := \Lambda + \ln({2\eta}/{|c|}),~~\Lambda:=2\eta( {t}/{(4\eta^{2}+\mu^{2})}+x).
	\end{aligned}
	$$
	We choose the parameters  $c>0$.
	\begin{itemize}  
		\item  The electric field $ E(x,t;\sigma_0) $ is a single-soliton of maximum amplitude $4\eta $ propagating at speed  $  -(4\eta^2 + \mu^2)^{-1} $. 
		\item  $ s(x,t;\sigma_0) $, carrying phase information, is a kink  propagating at speed  $  -(4\eta^2 + \mu^2)^{-1} $ with maximum amplitude  $4\eta^2/(4\eta^2+\mu^2)$. 
		\item  $ u(x,t;\sigma_0) $, the atomic inversion, is a single-soliton propagating on the background $-1$  at speed  $  -(4\eta^2 + \mu^2)^{-1} $ with maximum amplitude $(4\eta^2 - \mu^2)/(4\eta^2 + \mu^2)$.  
		\item   $r(x,t;\sigma_0) $, representing the atomic dipole,  is a single-soliton propagating at speed  $  -(4\eta^2 + \mu^2)^{-1} $ with maximum amplitude $4\eta\mu/(4\eta^2 + \mu^2) $.  
	\end{itemize}

	For $N>1$, the solutions $E(x,t;\sigma_0)$, $u(x,t;\sigma_0)$ and $r(x,t;\sigma_0)$ are localized traveling $N$-soliton, while $s(x,t;\sigma_0)$ corresponds to a localized $N$-kink. None of these solutions blow up within finite time. After interactions, the solitons or kinks asymptotically separate as $t\to \infty $ into  $N$ single-soliton or $N$ single-kink solutions traveling at speeds $ -(4\eta_k^2 + \mu^2)^{-1} $, one for each point in the discrete spectrum $\{z_k=i\eta_k\}_{k=1}^N$, where $\eta_j\neq \eta_k$ for $j\neq k$. 
	
\begin{figure}[h]
		\centering 

		 \begin{tikzpicture}
			% % 坐标轴
			% \draw[->] (-2,0) -- (6,0) node[right] {$x$};
			% \draw[->] (0,-1) -- (0,5) node[above] {$t$};
			
			% 特征线 x = x1 + v1 t
			\draw[line width=1pt] (2,0) -- (-2,2.4);
			\node[below] at (2,0) {$x_1$};
			\node[rotate=0] at (-1.5,1) {$x = x_1 + v_1 t$};
			
			% 特征线 x = x2 + v2 t
			\draw[line width=1pt] (3,0) -- (2, 2.4);
			\node[below] at (3,0) {$x_2$};
			\node[rotate=0] at (4,1) {$x = x_2 + v_2 t$};
			\fill[lightgray] (2,0) -- (0,1.2) -- (0,2.4) --(2,2.4)--(3,0) -- cycle; 
			\fill[lightgray] (0,1.2) -- (-2,2.4) -- (0,2.4) -- cycle; 
			% % 中间区域斜线填充
			% \begin{scope}
				%   \clip (2,0) -- (0,2.5) -- (0,5) -- (3,0) -- cycle;
				%   % \foreach \x in {1.4,1.5,...,2.9} {
					%   %   \draw[thin] (\x,0) -- (0,{2.5 - (\x-1.3)*1.3});
					%   % }
				%   \fill[blue!20] (2,0) -- (0,2.5) -- (0,5) --(2,5)--(3,0) -- cycle; 
				% \end{scope}
			
			% 坐标轴
			\draw[-{Stealth},line width=1pt] (-2,0) -- (6,0) node[right] {$x$};
			\draw[-{Stealth},line width=1pt] (0,-1) -- (0,2.7) node[above] {$t$};
			% 文字说明
			%\node at (-1.5,0.5) {\text{fast}};
			%\node at (4,0.5) {\text{slow}};
			
		\end{tikzpicture}
		\caption{The gray region represents the cone $C( x_{1},x_{2},v_{1},v_{2})$.  }\label{F7}
	\end{figure}
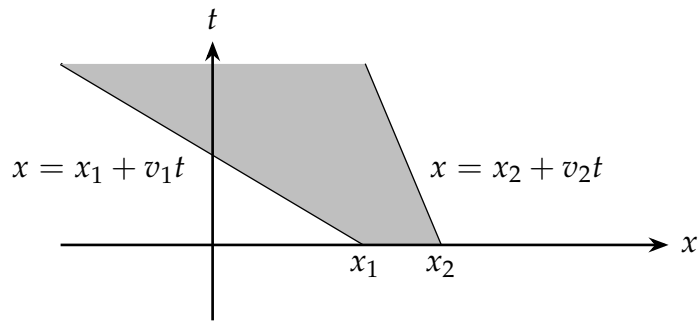

	Given $x _{1}\leq x _{2}$ and $-1/\mu ^{2} <v_{1}\leq v_{2} <0$, we analyze the long-time behavior of the solutions $E(x,t; \sigma_{0})$, $s(x,t; \sigma_{0})$, $u(x,t; \sigma_{0})$, and $r(x,t; \sigma_{0})$ within the cone
	\begin{align}\label{eq79}
		C( x_{1},x_{2},v_{1},v_{2}) :=\{ ( x,t) :x = x_{0}+vt,x_{0}\in [ x_{1},x_{2}] ,v\in[ v_{1},v_{2}] \},
	\end{align}
	as shown in Figure \ref{F7}.
	
	\begin{prop}\label{prop5.1}
		For given scattering date $\sigma_0=\{r(z)\equiv0,\{(z_k=i \eta_k,c_k)\}_{k=1}^N\}$. Divide the positive imaginary axis $i\mathbb{R}_+$ into three segments
		$$\mathcal J:=i\mathbb{R}_+\backslash (\mathcal J_-\cup \mathcal J_+), ~~\mathcal J_-:=i(0, {(-\mu^2-v_1^{-1})}^{1/2}/2),~~ \mathcal J_+:=i({(-\mu^2-v_2^{-1})}^{1/2}/2, \infty), $$
		and denote the index set 
		\begin{align*}
			\bigtriangleup_{\mathcal J}:=\{ k\in \{ 1,2,\ldots ,N\} :z_k \in \mathcal{Z}\cap \mathcal J, G( 0,\IIm z_k)  < 0\} .
		\end{align*}
		Then, for $(x,t)\in C(x_1,x_2,v_1,v_2)$ defined in \eqref{eq79} with $t\to \infty$, we have the asymptotics
		\begin{align*}
			M^{\bigtriangleup }( z;x,t|\sigma _{0}) =   ( I+\mathcal O( e^{2\mathring{ \kappa } t}) ) M^{\bigtriangleup_{\mathcal J} }(  z;x,t| \sigma _{0}( \mathcal J) ) ,
		\end{align*}
		uniformly for $z$ bounded away from $\mathcal{Z} \cup \overline{\mathcal{Z}}$. Here the constant $\mathring{ \kappa } <0$ is specified in \eqref{eq80},  and the restricted scattering data is given by 
		\begin{align*}
			\sigma_0( \mathcal J) =\{r(z)\equiv0,\{(z_k,c_k)\}_{k=1}^N: z_k\in \mathcal{Z}\cap \mathcal J  \}.
		\end{align*}
	\end{prop}
	
	\begin{proof}
		Assume $b_0$ satisfies the equation $G(0,b_0)=0$. Then $ib_0\in \mathcal J$ follows from the velocity constraint  $v_1\leq-(4b_0^2+\mu^2)^{-1}=v\leq v_2$ as $t\to \infty$. This implies the nonemptiness condition $\mathcal J\cap \{z\in i \mathbb R_+:  G( 0,\IIm z)  < 0\}\neq \emptyset$, ensuring the existence of the index set $\bigtriangleup_{\mathcal J}$. We now introduce the constant 
		\begin{align}\label{eq80}
			\mathring {\kappa}:= \max\{\kappa_-,\kappa_+\},
		\end{align}
		where $\kappa_\pm$ are determined by discrete spectral distributions
		\begin{align*}
			\begin{aligned}
				\kappa_{-}=\max _{z_{k}\in \mathcal Z\cap  \mathcal J_-}\left\{ -\eta_{k}( (4\eta_{k} ^{2}+\mu^{2})^{-1}+v_{1}) \right\} , ~~
				\kappa_{+}=\max _{z_{k}\in \mathcal Z\cap  \mathcal J_+}\left\{ \eta_{k}( (4\eta_{k} ^{2}+\mu^{2})^{-1}+v_{2}) \right\} . \end{aligned}
		\end{align*}
		From the definitions of the  intervals  $ \mathcal J _\pm$, it follows that the constants  $\kappa _\pm<0 $.  
		
		Given poles $z_k\in \mathcal Z\cap  \mathcal J_-$  (indexed by $k\in \bigtriangledown$) with complex conjugates $\overline z_k$ satisfying $\theta(\overline z_k)=-\theta(z_k)$ , we recall  \eqref{eq76.1} and $v_1\leq v$ to derive 
		\begin{align}\label{eq81}
			\lim _{t\rightarrow \infty }| e^{-2it\theta ( z_{k}) }| =e^{-2t\eta_{k} ((4\eta_{k}^{2}+\mu^2)^{-1}+v) }\leq \mathcal O( e^{2t\mathring \kappa}) .
		\end{align}
		Analogous decay estimates hold for $z_k\in \mathcal Z\cap  \mathcal J_+$ (indexed by $k\in \bigtriangleup$)  and  $\overline z_k$, indicating negligible contributions from the residues at $z_k \in \mathcal {Z\cap (J_-\cup J_+)}$ and $\overline z_k$ to the solution $M^{\bigtriangleup }$.
		
		Following \eqref{eq16}, we  equivalently reformulate the residue conditions of  $M^{\bigtriangleup }$ at the poles $z_k \in \mathcal {Z\cap (J_-\cup J_+)}$ and $\overline z_k$ as jump conditions for $\widehat M^\bigtriangleup(z;x,t|\sigma_0):=M^\bigtriangleup(z)N^{\bigtriangleup_{\mathcal J}}(z)$ on the composite contour $\widehat{\Sigma}$, as shown in Figure \ref{F8}, where
		\begin{align*}
			\widehat \Sigma:=\cup_{z_k \in \mathcal {Z\cap (J_-\cup J_+)}}(\partial D(z_k,\kappa)\cup\partial D(\overline z_k,\kappa)).
		\end{align*}
		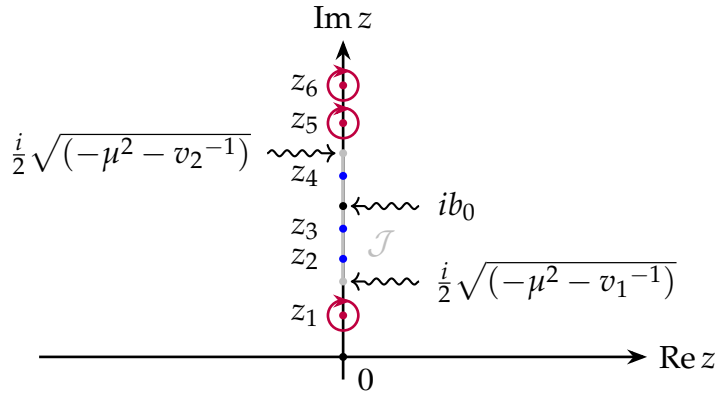
\begin{figure}
			\centering 
			\begin{tikzpicture}
				% 坐标轴
				\draw[-{Stealth},line width=1pt] (-4,0) -- (4,0) node[right] {$\RRe z$};
				\draw[-{Stealth},line width=1pt] (0,-0.3) -- (0,4.2) node[above] {$\IIm z$};
				
				\draw[lightgray,line width=1pt] (0,1) -- (0,2.7);

				\fill[lightgray] (0,2.7) circle (1.5pt);
				\node[black] at (-2.8, 2.7) {\(\frac{i}{2}\sqrt{(-\mu^2-{v_2}^{-1})}\)};
				\draw[->, thick, decorate, decoration={snake, amplitude=.4mm, segment length=2mm, post length=1mm}]  (-1,2.7) -- (-0.1,2.7);
				\fill[lightgray] (0,1) circle (1.5pt);
				\draw[->, thick, decorate, decoration={snake, amplitude=.4mm, segment length=2mm, post length=1mm}]  (1,1) -- (0.1,1);
				\node[black] at (2.8, 1) {\(\frac{i}{2}\sqrt{(-\mu^2-{v_1}^{-1})}\)};
				
				\node[lightgray] at (0.5, 1.5) {\(\mathcal J\)};
				\fill[blue] (0,2.4) circle (1.5pt);
				\node[black] at (-0.5, 2.4) {$z_4$};
				\fill[blue] (0,1.7) circle (1.5pt);
				\node[black] at (-0.5, 1.7) {$z_3$};
				\fill[blue] (0,1.3) circle (1.5pt);
				\node[black] at (-0.5, 1.3) {$z_2$};
				
				\fill[black] (0,2) circle (1.5pt);
				\draw[->, thick, decorate, decoration={snake, amplitude=.4mm, segment length=2mm, post length=1mm}]  (1,2) -- (0.1,2);
				\node[black] at (1.5, 2) {\(ib_0\)};
				
				\draw[purple, line width=1pt] (0,3.3) arc (90:-270:0.2);
				\fill[purple] (0,3.1) circle (1.5pt);
				\draw[purple, thick, -{Stealth}] (0,3.29) -- (0.1,3.29);
				\node[black] at (-0.5, 3.1) {$z_5$};
				
				\draw[purple, line width=1pt] (0,3.8) arc (90:-270:0.2);
				\fill[purple] (0,3.6) circle (1.5pt);
				\draw[purple, thick, -{Stealth}] (0,3.79) -- (0.1,3.79);
				\node[black] at (-0.5, 3.6) {$z_6$};
				
				\draw[purple, line width=1pt] (0,0.75) arc (90:-270:0.2);
				\fill[purple] (0,0.55) circle (1.5pt);
				\draw[purple, thick, -{Stealth}] (0,0.74) -- (0.1,0.74);
				\node[black] at (-0.5, 0.55) {$z_1$};
				
				% origin point O
				\fill (0, 0) circle (1.5pt);
				\node[black] at (0.3, -0.3) {0};
				
			\end{tikzpicture}
			\caption{In this example, among the purely imaginary eigenvalues $z_k$ ($k=1,\ldots,6$), $z_2$, $z_3$, and $z_4$ fall within the interval $\mathcal{J}$, and the corresponding solitons can be observed inside the cone $C(x_{1}, x_{2}, v_{1}, v_{2})$. The contour  $\widehat{\Sigma} := \cup_{k\in\{1,5,6\}} (\partial D(z_k,\kappa)\cup \partial D(\overline{z}_k,\kappa))$.}\label{F8}
		\end{figure}
		Specifically, we define
		\begin{align*}
			N^{\bigtriangleup_{\mathcal J}}(z):=
			\begin{cases}
				\begin{pmatrix}
					1 & -\frac{1}{z-z_k}c_{k}^{-1}[ ( F_{\bigtriangleup }^{-1}) '( z_{k})] ^{-2}e^{2it\theta( z_{k}) } \\
					0 & 1
				\end{pmatrix} &z\in D(z_k,\kappa),z_k\in  \mathcal Z\cap  \mathcal J_+\\
				\begin{pmatrix}
					1 & 0 \\
					-\frac{1}{z-z_k}c_{k} F_{\bigtriangleup }^{-2}( z_{k})e^{-2it\theta( z_{k}) } & 1
				\end{pmatrix}&z\in D(z_k,\kappa),z_k\in  \mathcal Z\cap  \mathcal J_-\\
				\begin{pmatrix}
					1 & 0 \\
					\frac{1}{z-\overline z_k}\overline c_{k}^{-1}( F_{\bigtriangleup }'(\overline{z}_{k})) ^{-2}e^{ -2it\theta (\overline{z}_{k}) } & 1
				\end{pmatrix}&z\in D(\overline z_k,\kappa),z_k\in  \mathcal Z\cap  \mathcal J_+\\
				\begin{pmatrix}
					1 & 	\frac{1}{z-\overline z_k}\overline c_{k}  F_{\bigtriangleup }^{2}( \overline{z}_{k})e^{ 2it\theta (\overline{z}_{k}) } \\
					0 & 1
				\end{pmatrix}&z\in D(\overline z_k,\kappa),z_k\in  \mathcal Z\cap  \mathcal J_-\\
				I &elsewise.
			\end{cases}
		\end{align*}
		Form the series of estimations similar to \eqref{eq81}, we obtain 
		\begin{align}\label{eq83}
			\|N^{\bigtriangleup_{\mathcal J}}(z)-I\|_{L^\infty (\widehat \Sigma)}\lesssim \mathcal O (e^{2t\mathring{ \kappa }}).
		\end{align}
		Thus, $\widehat M^\bigtriangleup(z)$ has continuous boundary values $\widehat M_\pm^\bigtriangleup(z)$ on $z\in \widehat \Sigma$ satisfying
		\begin{align}\label{eq82}
			\widehat M_+^\bigtriangleup(z)=\widehat M_-^\bigtriangleup(z)(N^{\bigtriangleup_{\mathcal J}}(z))^{-1}.
		\end{align}
		Additionally, $\widehat M^\bigtriangleup(z)$ retains the residue conditions of $M^{\bigtriangleup }$  at $z_k\in \mathcal{Z\cap J}$ and $\overline z_k$. 
		
		The solution $M^{\bigtriangleup_{\mathcal J} }(  z;x,t| \sigma _{0}( \mathcal J) )$ corresponds to $\widehat M^\bigtriangleup(z;x,t|\sigma_0)$ without considering the jump conditions. Their relationship can be given by
		\begin{align*}
			\widehat M^\bigtriangleup(z;x,t|\sigma_0)=\widehat{E}(z)M^{\bigtriangleup_{\mathcal J} }(  z;x,t| \sigma _{0}( \mathcal J) ).
		\end{align*}
		On $z\in \widehat \Sigma$, \eqref{eq82} implies  $\widehat{E}(z)$  satisfies  the jump condition
		\begin{align*}
			\widehat{E}_+(z) =\widehat{E}_-(z) M^{\bigtriangleup_{\mathcal J} }(  z; \sigma _{0}( \mathcal J) )(N^{\bigtriangleup_{\mathcal J}}(z))^{-1}(M^{\bigtriangleup_{\mathcal J} }(  z; \sigma _{0}( \mathcal J) ))^{-1}.
		\end{align*}
		By \eqref{eq83}, the jump  matrix of $\widehat{E}(z)$ decays uniformly  to the identity matrix as $t\to \infty$. Thus, the existence and asymptotic behavior of $\widehat E(z)$, specifically $\widehat E(z) =I +\mathcal  O (e^{2t\mathring{ \kappa }}) $, are guaranteed by the theory of small-norm Riemann–Hilbert problems.
	\end{proof}
	
	\begin{Corollary}\label{Corollary1}
		For $(x,t)\in C(x_1,x_2,v_1,v_2)$ with $t\to \infty$, the $N$-soliton and  $N$-kink solutions asymptotically reduce to $N(\mathcal J)$-soliton and $N(\mathcal J)$-kink solutions, respectively,
		\begin{align*}
			\begin{aligned} 
				&E(x,t;\sigma_{0}) =E(x,t;\sigma_{0}(\mathcal J)) + \mathcal O (e^{2t\mathring{ \kappa }}),\quad&s(x,t;\sigma_{0}) =s(x,t;\sigma_{0}(\mathcal J)) + \mathcal O (e^{2t\mathring{ \kappa }}),\\
				&u(x,t;\sigma_{0}) =u(x,t;\sigma_{0}(\mathcal J)) + \mathcal O (e^{2t\mathring{ \kappa }}),\quad&r(x,t;\sigma_{0}) =r(x,t;\sigma_{0}(\mathcal J)) + \mathcal O (e^{2t\mathring{ \kappa }}), 
			\end{aligned}
		\end{align*}
		where the scattering date $\sigma_{0}$ and $\sigma_{0}(\mathcal J)$ are given in Proposition \ref{prop5.1} and $N(\mathcal J)=|\mathcal{Z\cap J}|$.
	\end{Corollary}
	
	\subsection{Local models near the stationary phase points $\pm\zeta_0$}
	
	Reapplying the conclusions of Lemma \ref{Lemma4}, the bound \eqref{eq74} provides a point-wise (non-uniform) estimate for the decay of the jump matrix $V^{(2)}$ toward the identity matrix within the disks $D(\pm \zeta_0,\kappa)$. To precisely evaluate the contribution of the jump conditions to the asymptotic expansion of $ M^{(2)}_{\mathrm{RHP}}(z) $, we introduce localized RHP \ref{RHP7} and RHP \ref{RHP8}. 
	
	Using the formula \eqref{eq110}, we construct the local model $M_{\pm \zeta_0}$ based on the solutions to RHP \ref{RHP7} and RHP \ref{RHP8}, which exactly match the jump conditions of $ M^{(2)}_{\mathrm{RHP}}(z) $ on $\Sigma^A\cup \Sigma^B$. Here, the contours $\Sigma^A$ and $ \Sigma^B$ are defined as 
	$$\Sigma^A:=(\cup_{k=1}^4 \Sigma_{k}^A)\cap D(\zeta_0,\kappa),~~~~\Sigma^B:=(\cup_{k=1}^4 \Sigma_{k}^B)\cap D(-\zeta_0,\kappa),$$ 
	as shown in Figure \ref{F6}.
	\begin{figure}[h] 
		\centering 
			\begin{tikzpicture}
			\draw[black, line width=1pt] (-5.84,0.34) -- (-5.16,-0.34);
			\draw[black, line width=1pt] (-5.84,-0.34) -- (-5.16,0.34);
			\draw[black, -{Stealth},line width=1pt] (-5.84,0.34) -- (-5.6,0.1);
			\draw[black, -{Stealth},line width=1pt] (-5.84,0.34) -- (-5.2,-0.3);
			\draw[black, -{Stealth},line width=1pt] (-5.84,-0.34) -- (-5.6,-0.1);
			\draw[black, -{Stealth},line width=1pt] (-5.84,-0.34) -- (-5.2,0.3);
			\node[black] at (-4.5, 0) {\(\Sigma^B\)};
			%left%zeta0
			
			\draw[black, line width=1pt] (5.84,0.34) -- (5.16,-0.34);
			\draw[black, line width=1pt] (5.84,-0.34) -- (5.16,0.34);
			\draw[black, -{Stealth},line width=1pt] (5.16,0.34) -- (5.4,0.1);
			\draw[black, -{Stealth},line width=1pt] (5.16,0.34) --(5.8,-0.3);
			\draw[black, -{Stealth},line width=1pt] (5.16,-0.34) --  (5.4,-0.1);
			\draw[black, -{Stealth},line width=1pt] (5.16,-0.34) -- (5.8,0.3);
			\node[black] at (4.5, 0) {\(\Sigma^A\)};
			\fill (-5.5,0) circle (1.5pt);
			\node[black] at (-5.5, -0.5) {\(-\zeta_0\)};
			
			%right
			%zeta0
			\fill (5.5,0) circle (1.5pt);
			\node[black] at (5.5, -0.5) {\(\zeta_0\)};

			\fill (0, 0) circle (1.5pt);
			\node[black] at (0.3, -0.3) {0};			
			
		\end{tikzpicture}
		\caption{The contours $\Sigma ^B$ (left) and $\Sigma ^A$ (right).}\label{F6}
	\end{figure}
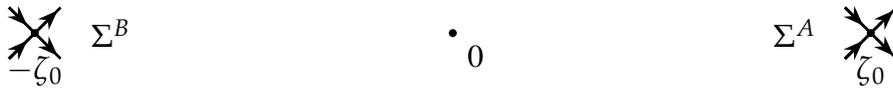
	\begin{rhp}\label{RHP7}
		Seek for an analytic matrix-valued function $M_A(z)$ defined on  $z\in \mathbb{C}\backslash \Sigma^A$ that satisfies the following conditions.
		\begin{itemize}
			\item The function $M_A(z)$ has continuous boundary values $M_{\text{A},\pm}(z)$ on $z\in \Sigma^A$, which satisfy the jump conditions
			\begin{align*}%\label{eq85}
				\begin{aligned}
					M_{A,+}(z)=M_{A,-}(z)V_A(z),
				\end{aligned}
			\end{align*}
			where the jump matrix $V_A(z):=V^{(2)}\upharpoonright_{\Sigma^A}$.
			\item The matrix  $M_A(z)$ possesses the asymptotic behavior $M_A(z)\to I$, as $|z|\to \infty$.
		\end{itemize}
	\end{rhp}
	\begin{rhp}\label{RHP8}
		Seek for an analytic matrix-valued function $M_B(z)$ defined on  $z\in \mathbb{C}\backslash \Sigma^B$ that satisfies the following conditions.
		\begin{itemize}
			\item The function $M_B(z)$ has continuous boundary values $M_{\text{B},\pm}(z)$ on $z\in \Sigma^B$, which satisfy the jump conditions
			\begin{align*}%\label{eq86}
				\begin{aligned}
					M_{B,+}(z)=M_{B,-}(z)V_B(z),
				\end{aligned}
			\end{align*}
			where the jump matrix $V_B(z):=V^{(2)}\upharpoonright_{\Sigma^B}$.
			\item The matrix  $M_B(z)$ possesses the asymptotic behavior $M_B(z)\to I$, as $|z|\to \infty$.
		\end{itemize}
	\end{rhp}
	
	We begin by analyzing RHP \ref{RHP7}. We first express the jump matrix $V_A(z)$ in terms of rescaled coordinates on the  $ \xi $-plane. Then, a solvable model, whose jump conditions correspond to the parabolic cylinder problem, can be used to effectively approximate $M_A(z)$ as $t\to \infty$. This methodology is frequently employed to derive long-time asymptotic behavior in the inverse scattering framework.
	
	Since $\zeta_0$ a stationary phase point of $\theta(z)$, it is a root of the equation
	$$(4z^2 - \mu^2)^2 {x}/{t} + 4z^2 + \mu^2 = 0.$$ 	
	Consequently, the phase function $\theta(z)$ admits an alternative representation
	\begin{align*}%\label{eq87}
		\begin{aligned}
			\theta(z) =\frac{8\zeta_0^{3}}{( 4\zeta _{0}^{2}-\mu ^{2})^{2}}+\frac{4 }{\beta}( z-\zeta_0) ^{2}+\left(\frac{8}{\beta(\mu -2z)}-\frac{4\alpha_1}{\mu^2-4z^2}\right)( z-\zeta_0) ^{3},
		\end{aligned}
	\end{align*}
	where the positive constants
	\begin{align}\label{eq87.1}
		\alpha_1:=\frac{\mu }{( 2\zeta _{0}+\mu) ^{3}},~~\beta :=\frac{( 4\zeta_0^{2}-\mu ^{2}) ^{3}}{\zeta _{0}( 4\zeta_{0}^{2}+3\mu ^{2}) }.
	\end{align}
	Let $\xi =\xi(z)$ and define the scaling operator $N_A$ as
	\begin{align*}%\label{eq88}
		\begin{aligned}
			f(z)\mapsto (N_Af)(\xi)=f\left(\sqrt{\beta/(16t) }\xi +\zeta_0\right),
		\end{aligned}
	\end{align*}
	which maps functions of $z$ to functions of the rescaled variable $\xi$. For  $z\in \Sigma^A$, we calculate
	\begin{align*}%\label{eq89}
		\begin{aligned}
			(N_AF_A^2e^{2it\theta})(\xi)=\delta_0^A\delta_1^A,
		\end{aligned}
	\end{align*}
	where 
	\begin{align}\label{eq90}
		\begin{aligned}
			\delta _{0}^{A}&:=\left( \frac{\sqrt{\beta }}{8\zeta_0\sqrt{t }}\right) ^{-2i\nu(\zeta_0) }\prod _{k\in \Delta }\left( \frac{\zeta_0-\overline{z}_{k}}{\zeta_0-z_{k}}\right) ^{2}e^{2i\beta ( \zeta_0,\zeta_0) +\frac{16it\zeta_0^{3}}{( 4\zeta_0^{2}-\mu ^{2}) ^{2}}},\\
			\delta_{1}^{A}&:=( -\xi )^{ -2i\nu(\zeta_0) }\left( \frac{\xi\sqrt{\beta  }}{8\zeta_0\sqrt{t }}+1\right) ^{2i\nu( \zeta _{0}) }e^{\frac{i}{2}\xi^2 +2it\left( \frac{8}{\beta ( \mu - \sqrt{\beta /( 4t) }\xi -2\zeta_0)}-\frac{4\alpha_1}{\mu^{2}-4( \sqrt{\beta /( 16t) }\xi +\zeta_0)^{2}} \right)( \sqrt{\beta /( 16t) }\xi) ^{3}}.
		\end{aligned}
	\end{align}
	Note that $\delta _{0}^{A}$ is independent of $\xi$ and $|\delta _{0}^{A}|=1, \overline{\delta _{0}^{A}}=1/\delta _{0}^{A}$.
	
	Now, we define the contours
	\begin{align*}
		\begin{aligned}
			\Sigma ':=\{ \xi =4\gamma e^{\pm i\pi /4}\sqrt{t/\beta },-\kappa  <\gamma  < \kappa \} ,~~\Sigma ^{\infty }:=\{ \xi =4\gamma e^{\pm i\pi /4}\sqrt{t/\beta },-\infty  <\gamma  < \infty \} 
		\end{aligned}
	\end{align*}
	oriented consistently with $\Sigma^A$, as shown in Figure \ref{F9}.
	\begin{figure}[h]
		\centering
		\begin{tikzpicture}
			\draw[black,dashed, thick] (-3,0) -- (3,0) node[right] {\text{Re} $\xi$};
			\draw[white,very thick] (-1.7,-1.7) -- (0,0);
			\draw[white,very thick] (1.7,1.7) -- (0,0);
			\draw[white,very thick] (1.7,-1.7) -- (0,0);
			\draw[white,very thick] (-1.7,1.7) -- (0,0);
			\draw[black, line width=1pt]  (-0.6,-0.6)-- (0,0);
			\draw[black,line width=1pt] (0.6,0.6) -- (0,0);
			\draw[black,line width=1pt] (0.6,-0.6) -- (0,0);
			\draw[black,line width=1pt] (-0.6,0.6) -- (0,0);
			\draw[black, -{Stealth}, line width=1pt] (0,0) -- (0.5,-0.5);
			\draw[black, -{Stealth}, line width=1pt] (0,0) -- (0.5,0.5);
			\draw[black, -{Stealth}, line width=1pt] (-0.6,0.6) -- (-0.25,0.25);
			\draw[black, -{Stealth}, line width=1pt] (-0.6,-0.6) -- (-0.25,-0.25);
			% 标记原点
			\draw[black, line width=1.5pt] (0,0) -- (0,0);
			\node[black] at (0,-0.3) {\text{0}};
			% 标记分支
			\node[black] at (-1.2,-0.4) {\(\Sigma_3'\)};
			\node[black] at (-1.2,0.4) {\(\Sigma_2'\)};
			\node[black] at (1.2,-0.4) {\(\Sigma_4'\)};
			\node[black] at (1.2,0.4) {\(\Sigma_1'\)};
		\end{tikzpicture}
		\begin{tikzpicture}
			\draw[black,dashed, thick] (-3,0) -- (3,0) node[right] {\text{Re} $\xi$};
			% 斜线（分支切割）
			\draw[black, line width=1pt] (-1.35,-1.35) -- (0,0);
			\draw[black, line width=1pt] (1.35,1.35) -- (0,0);
			\draw[black, line width=1pt] (1.35,-1.35) -- (0,0);
			\draw[black, line width=1pt] (-1.35,1.35) -- (0,0);
			\draw[black, -{Stealth}, line width=1pt] (0,0) -- (0.5,-0.5);
			\draw[black, -{Stealth}, line width=1pt] (0,0) -- (0.5,0.5);
			\draw[black, -{Stealth}, line width=1pt] (-0.6,0.6) -- (-0.25,0.25);
			\draw[black, -{Stealth}, line width=1pt] (-0.6,-0.6) -- (-0.25,-0.25);
			% 标记原点
			\draw[black, line width=1.5pt] (0,0) -- (0,0);
			\node[black] at (0,-0.3) {\text{0}};
			% 标记外部路径
			\node[black] at (-1.7,-1.4) {\(\Sigma_3^\infty\)};
			\node[black] at (-1.7,1.4) {\(\Sigma_2^\infty\)};
			\node[black] at (1.7,-1.4) {\(\Sigma_4^\infty\)};
			\node[black] at (1.7,1.4) {\(\Sigma_1^\infty\)};
		\end{tikzpicture}
		\caption{The contours $\Sigma '$ (left) and $\Sigma ^\infty$ (right), each composed of segments $\Sigma '_k$ and $\Sigma ^\infty_k$ (for $k=1,\ldots ,4$), respectively.}\label{F9}
	\end{figure}
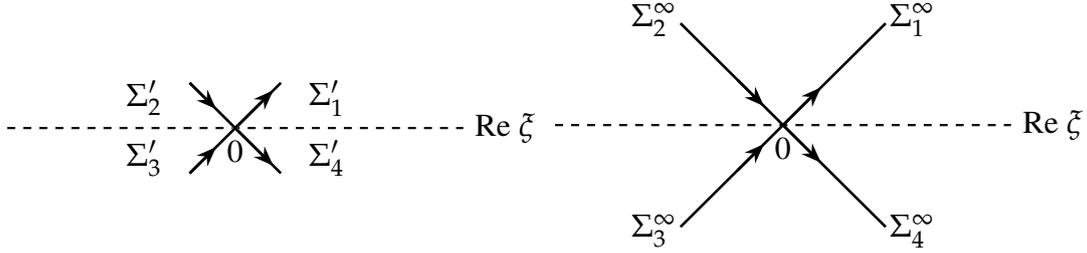
	\begin{Lemma}\label{Lemma3}
		Let $0 < \lambda < 1/2 $. As $ t \to \infty $, there exists a small constant $\kappa>0$ such that the following asymptotic estimate holds:  
		\begin{align}\label{eq91}
			|(\delta_{1}^{A})^{\pm1}-( -\xi )^{ \mp2i\nu(\zeta_0) }e^{\pm{i}\xi^2/2}|\leq c|e^{\pm i\lambda \xi^2/2}|t^{-1/2}, \xi \in \Sigma'.
		\end{align}
		As a consequence, we obtain the norm estimate
		\begin{align}\label{eq91.1}
			\|(\delta_{1}^{A})^{\pm1}-( -\xi )^{ \mp2i\nu(\zeta_0) }e^{\pm{i}\xi^2/2}\|_{L^1(\Sigma')\cap L^2(\Sigma')\cap L^\infty(\Sigma')}\leq ct^{-1/2},
		\end{align}
		where the  signs $\pm$ correspond to $\xi\in \Sigma_{1}'\cup \Sigma_{3}'$ and $ \xi\in \Sigma_{2}'\cup \Sigma_{4}'$, respectively.  
	\end{Lemma}
	
	\begin{proof} Write the left-hand side of \eqref{eq91} as
		\begin{align*}
			| (-\xi)^{-2i\nu(\zeta_0) }| \left| e^{i\lambda \xi ^{2}} \left\{ \left[ \left( {\xi\sqrt{\beta}}/(8\zeta_0\sqrt{t})+1\right) ^{2i\nu( \zeta _{0}) }-1\right]g(\xi,1) +g(\xi,1) -g(\xi,0) \right\}\right|
		\end{align*}
		where
		\begin{align*}
			g(\xi,s):=e^{i\xi^2(1-2\lambda) f(\xi,s)/2},
		\end{align*}
		and 
		\begin{align*}
			f(\xi,s):=1+\frac{s}{1-2\lambda}\frac{\xi\sqrt{\beta}}{4\sqrt{t}}\left( \frac{2}{\mu-\xi\sqrt{\beta}/(2\sqrt{t})-2\zeta_0}-\frac{\alpha_1\beta}{\mu^2-4(\xi\sqrt{\beta}/(4\sqrt{t})+\zeta_0)^2}\right).
		\end{align*}
		By parameterizing $z\in \Sigma_{1}^A\cap D(\zeta_0,\kappa)$ as $$z:=\zeta_0 +\gamma e^{i\pi/4},~~0\leq \gamma <\kappa,$$  we obtain the corresponding rescaled variable $\xi(\gamma) = 4 \gamma e^{i\pi/4}\sqrt{t/\beta } $. Recalling the definition of function $\nu(z)$ in \eqref{eq38}, we see that $|( -\xi )^{-2i\nu(\zeta_0)}|$ is  bounded, since
		\begin{align*}
			\begin{aligned}
				|( -\xi )^{-2i\nu(\zeta_0)}|=e^{2\nu(\zeta_0)\arg (-\xi)}=(1+|r(\zeta_0)|^2)^{3/4}.
			\end{aligned}
		\end{align*}
		Moreover, the terms $|e^{i\lambda \xi^2/2}|$ and $|e^{i(1-2\lambda) \xi^2/2}|$ are trivially bounded. Applying the mean value theorem, the term
		\begin{align*}
			\begin{aligned}
				&\left|\left( {\xi\sqrt{\beta}}/(8\zeta_0\sqrt{t})+1\right) ^{2i\nu( \zeta _{0})}-1\right| =\left|\int_1^{1+\xi\sqrt{\beta}/(8\zeta_0\sqrt{t})}ds^{2i\nu(\zeta_0)}\right|\\
				&\leq \sup\left\{|2i\nu(\zeta_0) s ^{2i\nu(\zeta_0) -1}|:s=1+\tau \xi\sqrt{\beta}/(8\zeta_0\sqrt{t}), 0\leq\tau \leq1 \right\} |\xi|\sqrt{\beta}/(8\zeta_0\sqrt{t})\\
			\end{aligned}
		\end{align*}
		is bounded. Since for $s=1+\tau\gamma e^{i\pi/4}/(2\zeta_0),$ we have
		\begin{align*}
			|2i\nu(\zeta_0) s ^{2i\nu(\zeta_0) -1}|=2|\nu(\zeta_0)||s|^{-1} e^{-2\nu(\zeta_0) \arg s}< \infty.
		\end{align*}
		Furthermore, we have 
		\begin{align*}
			\begin{aligned}
				\left|e^{i\lambda \xi^2/2}\left(\left( {\xi\sqrt{\beta}}/(8\zeta_0\sqrt{t})+1\right) ^{2i\nu( \zeta _{0})}-1\right)\right| \leq c\left|e^{i\lambda \xi^2/2} \xi/\sqrt{t}\right|\leq ct^{-1/2}.
			\end{aligned}
		\end{align*}
		The exponential term 
		\begin{align*}
			\left| g(\xi(\gamma),s)\right|=e^{-8(1-2\lambda)t\gamma^2 \RRe[ f(\xi(\gamma),s)]/\beta} 
		\end{align*}
		is uniformly bounded, where
		\begin{align*}
			f(\xi(\gamma),s):=1+\frac{s\gamma e^{i\pi/4}}{1-2\lambda}\left( \frac{2}{\mu-2\gamma e^{i\pi/4}-2\zeta_0}-\frac{\alpha_1\beta}{\mu^2-4(\gamma e^{i\pi/4}+\zeta_0)^2}\right).
		\end{align*}
		By the continuity of 	$f(\xi(\gamma),s)$  with respect to $\gamma$ and the fact $f(\xi(0),s)=1$, there exists a small $\kappa>0$ such that  $$\RRe (f(\xi(\gamma),s))>0 ,~~\text{for}~ 0\leq\gamma <\kappa. $$ Furthermore, Remark \ref{Remark1} ensures that the terms inside the parentheses in $f(\xi(\gamma),s)$ remain bounded. Similarly, we derive
		\begin{align*}
			\left|e^{i\lambda \xi^2/2}\left( g(\xi,1) -g(\xi,0)\right) \right|=\left|e^{i\lambda \xi^2/2} \right|\left|\int_0^1 \frac{\partial}{\partial s}g(\xi,s)ds \right|
			\leq c\left|e^{i\lambda \xi^2/2} \xi^3/\sqrt t\right|\leq ct^{-1/2}.
		\end{align*}
	\end{proof}
	
	In context of the Beals–Coifman theory, on the contour $\Sigma^\infty$, we define a pair of matrix-valued functions 
	\begin{align*}
		\begin{aligned}
			w_{A}^{+}(\xi) :=\begin{cases}
				\begin{pmatrix}
					0 & \frac{\overline{r}(\zeta_0) }{1+|r(\zeta_0)|^2}\delta _{1}^{A} \\
					0 & 0
				\end{pmatrix}  &\xi \in \Sigma'_1\\
				\begin{pmatrix}
					0 & \overline{r}(\zeta_0) \delta_1^A \\
					0 & 0
				\end{pmatrix}&\xi \in \Sigma _{3}'
			\end{cases},
			w_{A}^{-}(\xi):=\begin{cases}\begin{pmatrix}
					0 & 0 \\
					r(\zeta_0)/\delta _{1}^{A} & 0
				\end{pmatrix}&\xi \in \Sigma_{2}'\\
				\begin{pmatrix}
					0 & 0 \\
					\frac{ r(\zeta_0) }{1+|r(\zeta_0)|^{2}}{1}/{\delta ^{A}_1} & 0
				\end{pmatrix}&\xi \in \Sigma'_4\end{cases},
		\end{aligned}
	\end{align*}
	where $w_{A}^{\pm}(\xi) $ vanish identically on the remaining parts of $\Sigma^\infty$. Furthermore, introducing the Cauchy projection operators defined on $\Sigma^\infty$
	\begin{align*}
		\left( C_\pm f\right) (\xi)=\lim _{\varepsilon \to 0^{+}}\frac{1}{2\pi i }\int _{\Sigma^\infty }\frac{f(s)}{s-( \xi\pm i\varepsilon ) }ds,~~\xi\in \Sigma^\infty, 
	\end{align*}
	we define the operator
	\begin{align*}
		A:=C_+( \cdot w_A^-)+C_-(\cdot w_A^+).
	\end{align*}
	Let $w_A(\xi):=w_A^+(\xi)+w_A^-(\xi)$. The matrix function
	\begin{align}\label{eq100}
		M_{A^t}(\xi) :=I+\frac{1}{2\pi i}\int _{\Sigma ^{\infty }}\frac{( ( 1-A) ^{-1}I) (s) w_{A}(s)}{s-\xi }ds,~~\xi \in \mathbb{C} \backslash \Sigma ^{\infty }
	\end{align}
	solves the RHP
	\begin{align}\label{eq101.11}
		\begin{cases}
			M _{A^t,+}(\xi) =M _{A^t,-}(\xi) V_{A^{t}}(\xi) &\xi\in \Sigma ^{\infty }\\
			M_{A^t}(\xi) \rightarrow I &|\xi|\rightarrow \infty ,\end{cases}
	\end{align}
	where the jump matrix is given by $V_{A^t}(\xi)=I+w_A(\xi).$ It is easy to verify that 
	\begin{align*}
		V_{A^t}(\xi)=\begin{cases}
			(\delta_0^A)^{-1/2\sigma_3}(N_AV_A)(\xi)(\delta_0^A)^{1/2\sigma_3}&\xi\in {\Sigma'}\\
			I&\xi \in \Sigma^\infty \backslash\Sigma',
		\end{cases}
	\end{align*}
	which formulates the relationship between $M_A(z)$ and $M_{A^t}(\xi)$. 
	
	Building on Lemma \ref{Lemma3}, we further apply the Beals–Coifman theory to relate the RHP \eqref{eq101.11} to the solvable model associated with the parabolic cylinder problem described in \eqref{eq101.12}. To this end, on $\Sigma^\infty$, we introduce another pair of matrix-valued functions 
	\begin{align}\label{eq105}
		\begin{aligned}
			w_{A^\infty}^{+}(\xi) &:=\begin{cases}
				\begin{pmatrix}
					0 & \frac{\overline{r}(\zeta_0) }{1+|r(\zeta_0)|^2}( -\xi )^{ -2i\nu(\zeta_0) }e^{\frac{i}{2}\xi^2}  \\
					0 & 0
				\end{pmatrix}  &\xi \in \Sigma^\infty_1\\
				\begin{pmatrix}
					0 & \overline{r}(\zeta_0) ( -\xi )^{ -2i\nu(\zeta_0) }e^{\frac{i}{2}\xi^2}  \\
					0 & 0
				\end{pmatrix}&\xi \in \Sigma _{3}^\infty
			\end{cases}\\
			w_{A^\infty}^{-}(\xi)&:=\begin{cases}\begin{pmatrix}
					0 & 0 \\
					r(\zeta_0)( -\xi )^{ 2i\nu(\zeta_0) }e^{-\frac{i}{2}\xi^2}  & 0
				\end{pmatrix}&\xi \in \Sigma_{2}^\infty\\
				\begin{pmatrix}
					0 & 0 \\
					\frac{ r(\zeta_0) }{1+|r(\zeta_0)|^{2}}(- \xi )^{ 2i\nu(\zeta_0) }e^{-\frac{i}{2}\xi^2} & 0
				\end{pmatrix}&\xi \in \Sigma^\infty_4\end{cases}
		\end{aligned}
	\end{align}
	where $w_{A^\infty}^{\pm}(\xi)$ vanish identically on the remaining parts of $ \Sigma^\infty$. Similarly, define the corresponding operator
	\begin{align*}
		A^\infty:=C_+( \cdot w_{A^\infty}^-)+C_-(\cdot w_{A^\infty}^+),
	\end{align*}
	and set $w_{A^\infty}(\xi):=w_{A^\infty}^+(\xi)+w_{A^\infty}^-(\xi)$, with jump matrix $V_{A^\infty}(\xi):=I+w_{A^\infty}(\xi)$. The RHP  
	\begin{align}\label{eq101.12}
		\begin{cases}
			M _{A^\infty,+}(\xi) =M _{A^\infty,-}(\xi) V_{A^{\infty}}(\xi) &\xi\in \Sigma ^{\infty }\\
			M_{A^\infty}(\xi) \rightarrow I &|\xi|\rightarrow \infty ,\end{cases}
	\end{align}
	admits a unique solution given by
	\begin{align}\label{eq101}
		M_{A^\infty}(\xi) :=I+\frac{1}{2\pi i}\int _{\Sigma ^{\infty }}\frac{( ( 1-A^\infty) ^{-1}I) (s) w_{A^\infty}(s)}{s-\xi }ds,~~\xi \in \mathbb{C} \backslash \Sigma ^{\infty }.
	\end{align}
	
	For $0<\lambda<1/2$, we have  the estimate
	\begin{align}\label{eq101.1}
		\begin{aligned}| ( -\xi) ^{\mp2i\nu( \zeta _{0}) }e^{\pm i\xi^{2}/2}| 
			\leq e^{\mp 2\nu(\zeta_0)\arg ( -\xi) }| e^{\pm i\xi^{2}\lambda/2 }| e^{-8\kappa^{2}( 1-\lambda ) {t}/{\beta }}
			\lesssim t^{-1/2}| e^{\pm i\xi^{2}\lambda/2 }|, \end{aligned}
	\end{align}
	where  the signs $\mp$ correspond to $\xi\in (\Sigma_{1}^\infty\cup \Sigma_{3}^\infty)\backslash\Sigma'$ and $ \xi\in  (\Sigma_{2}^\infty\cup \Sigma_{4}^\infty)\backslash\Sigma'$, respectively.  
	From Lemma \ref{Lemma3} and \eqref{eq101.1}, there exists a constant $c>0$ such that  
	\begin{align*}
		\|A-A^\infty\|_{L^2(\Sigma^\infty)}\leq ct^{-1/2}.
	\end{align*}
	As a consequence, through rigorous analysis following the methodology in \cite{MR1207209}, we derive conclude that  as $t\to \infty$
	\begin{align}\label{eq102}
		\|(1-A^\infty)^{-1}\|_{L^2(\Sigma^\infty)},~~\|(1-A)^{-1}\|_{L^2(\Sigma^\infty)}<\infty.
	\end{align}
	According to \eqref{eq100} and \eqref{eq101}, for large $|\xi|$, the asymptotic expansions of $M_{A^t}(\xi)$ and $M_{A^\infty}(\xi) $ can be written as
	\begin{align*}
		M_{A^t}(\xi) &:=I-\frac{M_{A^t,1}}{\xi}+\mathcal{O}(\xi^{-2})=I-\frac{1}{2\pi i \xi}\int _{\Sigma ^{\infty }}( ( 1-A) ^{-1}I) (s) w_{A}(s)ds+\mathcal{O}(\xi^{-2}),\\
		M_{A^\infty}(\xi) &:=I-\frac{M_{A^\infty,1}}{\xi}+\mathcal{O}(\xi^{-2})=I-\frac{1}{2\pi i \xi}\int _{\Sigma ^{\infty }}( ( 1-A^\infty) ^{-1}I) (s) w_{A^\infty}(s)ds+\mathcal{O}(\xi^{-2}).
	\end{align*}
	\begin{prop} As $t\to \infty$, the following holds:
		\begin{align*}
			\frac{1}{2\pi i}\int _{\Sigma ^{\infty }}( ( 1-A) ^{-1}I) (s) w_{A}(s)ds=\frac{1}{2\pi i }\int _{\Sigma ^{\infty }}( ( 1-A^\infty) ^{-1}I) (s) w_{A^\infty}(s)ds+\mathcal{O}(t^{-1/2}).
		\end{align*}
	\end{prop}
	\begin{proof}
		Set $\hat w(\xi):=w_{A^{\infty}}(\xi)-w_A(\xi)$ and define the operator $$C_{\hat w}:=C_+(\cdot \hat w^- )+C_-(\cdot \hat w^+ ),$$
		where $\hat w^\pm:=w_{A^{\infty}}^\pm-w_A^\pm $. Then, the integral associated with $A^\infty$ decomposes as
		\begin{align*}
			\begin{aligned} 
				\int_{\Sigma^{\infty }}( (1-A^{\infty })^{-1}I) (s)  w_{A^{\infty }}(s) ds
				&=\int _{\Sigma ^{\infty }}( ( 1-A) ^{-1}I) (s) w_{A}(s) ds+\int_{\Sigma ^{\infty }}\hat{w}(s) ds\\
				&+\int_{ \Sigma ^{\infty }}( ( 1-A) ^{-1}C_{\hat{w}}I) (s)w_{A^{\infty }}(s) ds\\
				&+\int _{\Sigma ^{\infty }}( ( 1-A) ^{-1}AI) (s) \hat{w}(s) ds\\
				&+\int_{ \Sigma ^{\infty }}( ( 1-A) ^{-1}C_{\hat{w}}( 1-A^{\infty }) ^{-1}A^{\infty }I)(s)w_{A^{\infty }}(s) ds\\
				&:=\mathrm{I+II+III+IV+V}.\end{aligned}
		\end{align*}
		Applying the bounds from Lemma \ref{Lemma3}, \eqref{eq101.1}, and \eqref{eq102}, the terms $\mathrm{II-V}$ satisfy
		\begin{align*}
			\begin{aligned} 
				|\mathrm{II}|&\leq ct^{-1/2},\\
				|\mathrm{III}|&\leq c\|\hat{w}\|_{L^2(\Sigma^\infty)}\|w_{A^{\infty }}\|_{L^2(\Sigma^\infty)}\leq ct^{-1/2},\\
				|\mathrm{IV}|&\leq c\|w_A\|_{L^2(\Sigma^\infty)}\| \hat{w}\|_{L^2(\Sigma^\infty)}\leq ct^{-1/2},\\
				|\mathrm{V}|&\leq c\|{\hat{w}}\|_{L^2(\Sigma^\infty)}\| w_{A^{\infty }}\|^2_{L^2(\Sigma^\infty)}\leq ct^{-1/2}.
			\end{aligned}
		\end{align*}
	\end{proof}
	Through the standard computations, as detailed in \cite{MR1207209,MR3795020}, for solving the RHP associated with $M_{A^\infty}(\xi)$, we derive the explicit expression for the leading-order term $ M_{A^\infty,1} $ in the asymptotic expansion of $ M_{A^\infty}(\xi) $ 
	\begin{align}\label{eq104}
		M_{A^\infty,1}=\begin{pmatrix}
			0 &-i\beta_{12} \\
			i\beta_{21}& 0
		\end{pmatrix},
	\end{align}
	where the coefficients $\beta_{12} $ and $\beta_{21} $ are given by 
	\begin{align*}
		\beta_{12}= \frac{\sqrt{2\pi}e^{-\pi\nu(\zeta_0)/2}e^{-i\pi/4}}{r(\zeta_0)\Gamma(i\nu(\zeta_0))},~~\beta_{21}= -\frac{\sqrt{2\pi}e^{-\pi\nu(\zeta_0)/2}e^{i\pi/4}}{\overline r(\zeta_0)\Gamma(-i\nu(\zeta_0))},
	\end{align*}
	and satisfy the symmetry relation $\beta_{12}=-\overline \beta_{21}.$ 
	
	We then consider the RHP \ref{RHP8}. Utilizing the representation
	\begin{align*}%\label{eq92}
		\begin{aligned}
			\theta(z) =- \frac{8\zeta_0^3}{(  4\zeta_0^2-\mu^2)^2} - \frac{4}{\beta}(z+\zeta_0)^2 - \left( \frac{8}{\beta(\mu - 2z)} - \frac{4\alpha_2}{\mu^2 - 4z^2} \right) (z+\zeta_0)^3.
		\end{aligned}
	\end{align*}
	where the positive constant
	$$\alpha_2 := \frac{\mu}{(2\zeta_0 - \mu)^3}.$$
	Define the scaling operator $N_B$ as
	\begin{align*}%\label{eq93}
		\begin{aligned}
			f(z)\mapsto (N_Bf)(\xi)=f\left(\sqrt{\beta/(16t) }\xi -\zeta_0\right).
		\end{aligned}
	\end{align*}
	For $z\in \Sigma^B$, let's calculate
	\begin{align*}
		\begin{aligned}
			(N_BF_B^2e^{2it\theta})(\xi)=\delta_0^B\delta_1^B,
		\end{aligned}
	\end{align*}
	where 
	\begin{align*}
		\begin{aligned}
			\delta _{0}^{B}&:=\left( \frac{\sqrt{\beta }}{8\zeta_0\sqrt{t }}\right) ^{2i\nu(\zeta_0) }\prod _{k\in \Delta }\left( \frac{\zeta_0+\overline{z}_{k}}{\zeta_0+ z_{k}}\right) ^{2}e^{-2i\beta ( \zeta_0,\zeta_0) -\frac{16it\zeta_0^{3}}{( 4\zeta_0^{2}-\mu ^{2}) ^{2}}},\\
			\delta_{1}^{B}&:= \xi ^{ 2i\nu(\zeta_0) }\left( -\frac{\xi\sqrt{\beta  }}{8\zeta_0\sqrt{t }}+1\right) ^{-2i\nu( \zeta _{0}) }e^{-\frac{i}{2}\xi^2 -2it\left( \frac{8}{\beta ( \mu - \sqrt{\beta /( 4t) }\xi -2\zeta_0)}-\frac{4\alpha_2}{\mu^{2}-4( \sqrt{\beta /( 16t) }\xi +\zeta_0)^{2}} \right)( \sqrt{\beta /( 16t) }\xi) ^{3}}.
		\end{aligned}
	\end{align*}
	Comparing with $\delta _{0}^{A}$ defined in \eqref{eq90}, we find $\delta _{0}^{B}=\overline {\delta _{0}^{A}}$.  
	
	As $ t \to \infty $, with $0 < \lambda < 1/2 $, analogous to estimates \eqref{eq91} and \eqref{eq91.1}, we yield
	\begin{align*}%\label{eq95}
		|(\delta_{1}^{B})^{\pm1}-\xi ^{ \pm2i\nu(\zeta_0) }e^{\mp{i}\xi^2/2}|\leq c|e^{\mp i\lambda \xi^2/2}|t^{-1/2},~~\xi\in \Sigma',
	\end{align*}
	and
	\begin{align*}%\label{eq96}
		\|(\delta_{1}^{B})^{\pm1}-\xi ^{ \pm2i\nu(\zeta_0) }e^{\mp{i}\xi^2/2}\|_{L^1( \Sigma')\cap L^2( \Sigma')\cap L^\infty( \Sigma')}\leq ct^{-1/2},
	\end{align*}
	where the $\pm$ signs correspond to  $ \xi \in \Sigma_{2}'\cup \Sigma_{4}'$ and $\xi\in \Sigma_{1}'\cup \Sigma_{3}'$, respectively.  
	
	Thus, skipping the step of defining $M_{B^t}(\xi)$, we directly give the matrix-valued functions on $\Sigma^\infty$
	\begin{align}\label{eq106}
		\begin{aligned}
			w_{B^\infty}^{+}(\xi) &:=\begin{cases}
				\begin{pmatrix}
					0 & \frac{{r}(\zeta_0) }{1+|r(\zeta_0)|^2}\xi ^{ 2i\nu(\zeta_0) }e^{-\frac{i}{2}\xi^2}  \\
					0 & 0
				\end{pmatrix}  &\xi \in \Sigma^\infty_2\\
				\begin{pmatrix}
					0 & {r}(\zeta_0) \xi ^{ 2i\nu(\zeta_0) }e^{-\frac{i}{2}\xi^2}  \\
					0 & 0
				\end{pmatrix}&\xi \in \Sigma _{4}^\infty
			\end{cases}\\
			w_{B^\infty}^{-}(\xi)&:=\begin{cases}\begin{pmatrix}
					0 & 0 \\
					\overline r(\zeta_0)\xi ^{ -2i\nu(\zeta_0) }e^{\frac{i}{2}\xi^2}  & 0
				\end{pmatrix}&\xi \in \Sigma_{1}^\infty\\
				\begin{pmatrix}
					0 & 0 \\
					\frac{ \overline r(\zeta_0) }{1+|r(\zeta_0)|^{2}}\xi ^{ -2i\nu(\zeta_0) }e^{\frac{i}{2}\xi^2} & 0
				\end{pmatrix}&\xi \in \Sigma^\infty_3,\end{cases}
		\end{aligned}
	\end{align}
	where $w_{B^\infty}^{\pm}(\xi)$ also vanish identically on the remaining parts of $ \Sigma^\infty$. Similarly,
	\begin{align*}
		B^\infty:=C_+( \cdot w_{B^\infty}^-)+C_-(\cdot w_{B^\infty}^+).
	\end{align*}
	Let $w_{B^\infty}(\xi):=w_{B^\infty}^+(\xi)+w_{B^\infty}^-(\xi)$, and define $V_{B^\infty}(\xi):=I+w_{B^\infty}(\xi)$. The RHP  
	\begin{align*}
		\begin{cases}
			M _{B^\infty,+}(\xi) =M _{B^\infty,-}(\xi) V_{B^{\infty}}(\xi) &\xi\in \Sigma ^{\infty }\\
			M_{B^\infty}(\xi) \rightarrow I &|\xi|\rightarrow \infty .\end{cases}
	\end{align*}
	admits a unique solution
	\begin{align*}%\label{eq107}
		M_{B^\infty}(\xi) :=I+\frac{1}{2\pi i}\int _{\Sigma ^{\infty }}\frac{( ( 1-B^\infty) ^{-1}I) (s) w_{B^\infty}(s)}{s-\xi }ds,~~\xi \in \mathbb{C} \backslash \Sigma ^{\infty }.
	\end{align*}
	By comparing the definitions in \eqref{eq105} and \eqref{eq106}, we derive the symmetry relation
	\begin{align*}
		V_{A^\infty}(\xi)=\overline{V_{B^\infty}}(-\overline \xi),
	\end{align*}
	which, by the uniqueness of the solution to the solvable RHP, implies
	\begin{align}\label{eq107.1}
		M_{B^\infty}(\xi)=\overline{M_{A^\infty}}(-\overline \xi).
	\end{align}
	Assuming the asymptotic expansion of
	\begin{align*}
		M_{B^\infty}(\xi)=I-M_{B^\infty,1}/\xi+\mathcal O(\xi^{-2}).
	\end{align*}
	and substituting the corresponding expansion of $M_{A^\infty}$ into \eqref{eq107.1}, we obtain the leading-order term 
	\begin{align*}
		M_{B^\infty,1}=-\overline{M_{A^\infty,1}}.
	\end{align*}
	
	\begin{prop}
		Let $\xi = 4(z-\zeta_0)\sqrt{t/\beta}$, the solution to RHP \ref{RHP7} has the asymptotic expansion
		\begin{align}\label{eq108}
			M_{A}(z(\xi);x,t)=I+\frac{1}{\xi}\begin{pmatrix}
				0 &i\beta_{12}\delta_0^A \\
				i\overline{\beta_{12}\delta_0^A}& 0
			\end{pmatrix}+\mathcal{O}(t^{-1})+\mathcal{O}(\xi^{-2}).
		\end{align}
		where 
		$\delta_0^A$, $\beta_{12}$, and $\beta_{21}$ are defined in \eqref{eq90} and \eqref{eq104}.
		Similarly, let $\xi = 4(z+\zeta_0)\sqrt{t/\beta}$, the solution to RHP \ref{RHP8} satisfies 
		\begin{align}\label{eq109}
			M_{B}(z(\xi);x,t)=I+\frac{1}{\xi}\begin{pmatrix}
				0 &i\overline{\beta_{12}\delta_0^A} \\
				i\beta_{12}\delta_0^A& 0
			\end{pmatrix}+\mathcal{O}(t^{-1})+\mathcal{O}(\xi^{-2}).
		\end{align}
		Here
		\begin{align}\label{eq109.1}
			|\beta_{12}\delta_0^A|^2=|\nu(\zeta_0)|,
		\end{align}
		\begin{align}\label{eq109.2}
			\begin{aligned}\arg ( \beta _{12}\delta _{0}^{A}) &=-\pi/4-\arg r(\zeta_0)-\arg  \Gamma ( i\nu(\zeta_0)) -2\nu(\zeta_0) \ln ( \sqrt{\beta }/( 8\zeta_0\sqrt{t}) ) \\
				&+4\sum _{k\in \bigtriangleup }\arg ( \zeta_0-z_{k}) +2\int _{\zeta_0}^{\infty }\ln \left( \frac{s+\zeta_0}{s-\zeta_0}\right) d_{s}\nu (s) +\frac{16t\zeta_0^{3}}{( 4\zeta_0^{2}-\mu ^{2}) ^{2}}.\end{aligned}
		\end{align}
	\end{prop}
	
	As a result of the preceding analysis, we have independently examined the contributions of the jump matrix $ V^{(2)}(z) $ within neighborhoods of the stationary phase points $ \pm\zeta_0 $. By incorporating the local models via the piecewise construction
	\begin{align}\label{eq110}
		M_{\pm \zeta _{0}}(z) =\begin{cases} M_{\mathrm{out}}(z)M_{A}(z) &z\in D( \zeta_0,\kappa) \\
			M_{\mathrm{out}}(z)M_{B}(z) &z\in D( -\zeta_0,\kappa), \end{cases}
	\end{align}
	we unify these two contributions into $M_{\pm \zeta _{0}}(z)$,  which explicitly captures the dominant contributions from the pointwise-convergent components of the jump matrix $ V^{(2)}(z) $. Critically, this incorporation introduces only an exponentially decaying error term, as rigorously justified in \cite{MR1207209}.  
	
	\subsection{Construction and analysis of the error matrix $E_\circ(z)$}
\begin{figure}[h] 
		\centering 
			\begin{tikzpicture}
			% RRe z 
			%\draw[black, thick, dashed] (-8,0) -- (-3.5,0);
			%\draw[black, thick, dashed] (-2.5,0) -- (2.5,0);
			%\draw[black, thick, dashed] (3.5,0) -- (8,0) node[right] {\text{Re} $z$};
			%zeta1
			\fill (0,1.5) circle (1.5pt);
			\node[black] at (0, 2) {\(\zeta_1\)};
			%-zeta1
			\fill (0,-1.5) circle (1.5pt);
			\node[black] at (0, -2) {\(-\zeta_1\)};
			%left
			%sigma4c
			\draw[black, line width=1pt] (-2.5,0) -- (0,1.5);
			\draw[black, -{Stealth},line width=1pt] (-2.5,0) -- (-1,0.9);
			\node[black] at (-1.5,1.1) {\(\Sigma_4^C\)};
			%omega12
			%\node[black] at (-1.3, 0.3) {\(\Omega_{12}\)};
			
			%sigma3c
			\draw[black, line width=1pt] (-2.5,0) -- (0,-1.5);
			\draw[black,-{Stealth}, line width=1pt] (-2.5,0) -- (-1, -0.9);
			\node[black] at (-1.5,-1.1) {\(\Sigma_3^C\)};
			%omega11
			%\node[black] at (-1.3, -0.3) {\(\Omega_{11}\)};
			
			%sigma6D
			\draw[black, line width=1pt] (0,1.5) -- (0,0);
			\draw[black, -{Stealth},line width=1pt] (0,1.5) -- (0,0.5);
			\node[black] at (0.4, 0.8) {\(\Sigma_6^D\)};
			
			%sigma5D
			\draw[black, line width=1pt] (0,-1.5) -- (0,0);
			\draw[black, -{Stealth},line width=1pt] (0,-1.5) -- (0,-0.5);
			\node[black] at (0.4, -0.8) {\(\Sigma_5^D\)};
			
			%sigma1C
			\draw[black, line width=1pt] (2.5,0) -- (0,1.5);
			\draw[black, -{Stealth}, line width=1pt](0,1.5) -- (1.5,0.6);
			\node[black] at (1.5,1.1) {\(\Sigma_1^C\)};
			%omega9
			%\node[black] at (1.2, 0.3) {\(\Omega_{9}\)};
			
			%sigma2C
			\draw[black, line width=1pt] (2.5,0) -- (0,-1.5);
			\draw[black,-{Stealth}, line width=1pt] (0,-1.5) -- (1.5, -0.6);
			\node[black] at (1.5,-1.1) {\(\Sigma_2^C\)};
			%omega10
			%\node[black] at (1.3, -0.3) {\(\Omega_{10}\)};
			
			%sigma4D
			\draw[black, line width=1pt] (-3.5,2) -- (-3.5,0);
			\draw[black, -{Stealth},line width=1pt] (-3.5,2) -- (-3.5,1);
			\node[black] at (-3,1.8) {\(\Sigma_4^D\)};
			
			%sigma3D
			\draw[black, line width=1pt] (-3.5, -2) -- (-3.5,0);
			\draw[black, -{Stealth},line width=1pt] (-3.5, -2) -- (-3.5,-1);
			\node[black] at (-3, -1.8) {\(\Sigma_3^D\)};
			
			%sigmaB
			\draw[black, line width=1pt] (-8,-2.5) -- (-5.86,-0.36);
			\draw[black, line width=1pt] (-5.14,0.36) -- (-3.5,2);
			\draw[black, line width=1pt] (-8,2.5) -- (-5.86,0.36);
			\draw[black, line width=1pt] (-5.14,-0.36) -- (-3.5,-2);
			%sigma3B
			\draw[black, -{Stealth}, line width=1pt] (-8,-2.5) -- (-6.0,-0.5);
			\node[black] at (-6.5,-1.8) {\(\Sigma_3^B\)};
			%sigm1B
			\draw[black, -{Stealth}, line width=1pt] (-5.14,0.36) -- (-4.75,0.75);
			\node[black] at (-4.5,1.8) {\(\Sigma_1^B\)};
			%sigma4B
			\draw[black, -{Stealth}, line width=1pt] (-5.14,-0.36) -- (-4.75, -0.75);
			\node[black] at (-4.5,-1.8) {\(\Sigma_4^B\)};
			%sigm2B
			\draw[black, thick,black, -{Stealth}, line width=1pt] (-8,2.5) -- (-6.0,0.5);
			\node[black] at (-6.5,1.8) {\(\Sigma_2^B\)};
			
			%right
			%sigma1D
			\draw[black, line width=1pt] (3.5,2) -- (3.5,0);
			\draw[black, -{Stealth},line width=1pt] (3.5,0) -- (3.5,1);
			\node[black] at (3,1.8) {\(\Sigma_1^D\)};
			
			%sigma2D
			\draw[black, line width=1pt] (3.5, -2) -- (3.5,0);
			\draw[black, -{Stealth},line width=1pt] (3.5, 0) -- (3.5,-1);
			\node[black] at (3, -1.8) {\(\Sigma_2^D\)};
			
			%sigmaA
			\draw[black, line width=1pt] (8,-2.5) -- (5.86,-0.36);
			\draw[black, line width=1pt] (5.14,0.36) -- (3.5,2);
			\draw[black, line width=1pt] (8,2.5) -- (5.86,0.36);
			\draw[black, line width=1pt] (5.14,-0.36) -- (3.5,-2);
			
			%sigma4A
			\draw[black, -{Stealth}, line width=1pt](5.86,-0.36) -- (6.25,-0.75);
			\node[black] at (6.5,-1.8) {\(\Sigma_4^A\)};
			%sigm2A
			\draw[black, -{Stealth}, line width=1pt] (3.5,2) -- (5,0.5);
			\node[black] at (4.5,1.8) {\(\Sigma_2^A\)};
			%sigma1A
			\draw[black, -{Stealth}, line width=1pt] (5.86,0.36) -- (6.25, 0.75);
			\node[black] at (6.5,1.8) {\(\Sigma_1^A\)};
			%sigm3A
			\draw[black, thick,black, -{Stealth}, line width=1pt] (3.5,-2) -- (5,-0.5);
			\node[black] at (4.5,-1.8) {\(\Sigma_3^A\)};
			
			\draw[black,line width=1pt] (2.5,0) arc (180:0:0.5);			
			\draw[black, thick, -{Stealth}] (3,0.5) -- (3.1,0.5);			
			\node[black] at (3,0.8) {$\Sigma_1^\kappa$};			
			\draw[black,line width=1pt] (2.5,0) arc (-180:0:0.5);			
			\draw[black, thick, -{Stealth}] (3,-0.5) -- (3.1,-0.5);			
			\node[black] at (3,-0.8) {$\Sigma_2^\kappa$};
			\draw[black,line width=1pt] (-3.5,0) arc (180:0:0.5);			
			\draw[black, thick, -{Stealth}] (-3,0.5) -- (-2.9,0.5);		
			\node[black] at (-3,0.8) {$\Sigma_4^\kappa$};		
			\draw[black, line width=1pt] (-3.5,0) arc (-180:0:0.5);		
			\draw[black, thick, -{Stealth}] (-3,-0.5) -- (-2.9,-0.5);			
			\node[black] at (-3,-0.8) {$\Sigma_3^\kappa$};
			\draw[black, line width=1pt] (-6,0) arc (180:0:0.5);
			\draw[black, thick, -{Stealth}] (-5.5,0.5) -- (-5.4,0.5);		
			\draw[black,line width=1pt] (-6,0) arc (-180:0:0.5);
			\draw[black, line width=1pt] (5,0) arc (180:0:0.5);
			\draw[black, thick, -{Stealth}] (5.5,0.5) -- (5.6,0.5);		
			\draw[black,line width=1pt] (5,0) arc (-180:0:0.5);
			\node[black] at (-7, 0) {\(D( - \zeta_{0},\kappa) \)};
			\node[black] at (6.8, 0) {\(D( \zeta_{0},\kappa) \)};
			%\draw[green, line width=1pt] (-5.854,0.354) -- (-5.146,-0.354);
			%\draw[green, line width=1pt] (-5.854,-0.354) -- (-5.146,0.354);
			%\draw[green, line width=1pt] (5.854,0.354) -- (5.146,-0.354);
			%\draw[green, line width=1pt] (5.854,-0.354) -- (5.146,0.354);
			left%zeta0
			\fill (-5.5,0) circle (1.5pt);
			\node[black] at (-5.5, -0.5) {\(-\zeta_0\)};
			%omega5
			%\node[black] at (-4.8, 0.3) {\(\Omega_5\)};
			%omega8
			%\node[black] at (-4.8, -0.3) {\(\Omega_8\)};
			
			%omega6
			%\node[black] at (-6.2, 0.3) {\(\Omega_6\)};
			%omega7
			%\node[black] at (-6.2, -0.3) {\(\Omega_7\)};
			
			%-mu-kappa
			%\fill (-3.5,0) circle (1pt);
			%\draw[->, thick, decorate, decoration={snake, amplitude=.4mm, segment length=2mm, post length=1mm}] (-2.7, -0.8) -- (-3.4,-0.1);
			%\node[black] at (-2.8, -1) {\(-\frac{\mu}{2}-\kappa\)};
			%-mu+kappa
			%\fill (-2.5,0) circle (1pt);
			%\draw[->, thick, decorate, decoration={snake, amplitude=.4mm, segment length=2mm, post length=1mm}] (-2.7, 0.8) -- (-2.6,0.1);
			%\node[black] at (-2.8, 1) {\(-\frac{\mu}{2}+\kappa\)};
			
			%right
			%zeta0
			\fill (5.5,0) circle (1.5pt);
			\node[black] at (5.5, -0.5) {\(\zeta_0\)};
			
			%omega2
			%\node[black] at (4.8, 0.3) {\(\Omega_2\)};
			%omega3
			%\node[black] at (4.8, -0.3) {\(\Omega_3\)};
			
			%omega1
			%\node[black] at (6.2, 0.3) {\(\Omega_1\)};
			%omega4
			%\node[black] at (6.2, -0.3) {\(\Omega_4\)};
			
			%mu-kappa
			%\fill (3.5,0) circle (1.5pt);
			%\draw[->, thick, decorate, decoration={snake, amplitude=.4mm, segment length=2mm, post length=1mm}] (2.9, -0.8) -- (3.4,-0.1);
			%\node[black] at (2.8, -1) {\(\frac{\mu}{2}-\kappa\)};
			%mu+kappa
			%\fill (2.5,0) circle (1.5pt);
			%\draw[->, thick, decorate, decoration={snake, amplitude=.4mm, segment length=2mm, post length=1mm}] (2.7, 0.8) -- (2.6,0.1);
			%\node[black] at (2.8, 1) {\(\frac{\mu}{2}+\kappa\)};
			
			%omega13
			%\node[black] at (-1,2.3) {\(\Omega_{13}\)};
			
			% origin point O
			\fill (0, 0) circle (1.5pt);
			\node[black] at (0.3, -0.3) {\(\text{0}\)};			
		\end{tikzpicture}
		\caption{The jump contour $\Sigma^{E_{\circ}}$ for RHP \ref{RHP9}.}\label{F10}
	\end{figure}
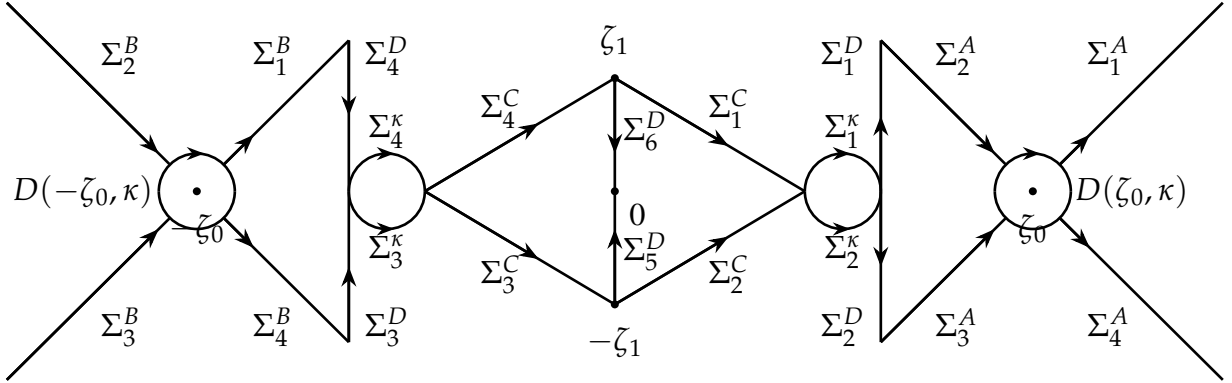

	By combining $M_{\mathrm{out}}(z)$ and $M_{\pm \zeta_0}(z)$, as given in the RHP \ref{RHP6} and \eqref{eq110} respectively, and reverting to \eqref{eq73}, we define the jump contour $\Sigma^{E_{\circ}}$ for $E_{\circ}(z;x,t)$ as
	\begin{align*}
		\Sigma ^{{E_{\circ}}}:=( \Sigma _{2}\cup \partial D( \pm \zeta_{0},\kappa) )\backslash ( \Sigma ^{A}\cup\Sigma ^{B}) ,
	\end{align*}
	as shown in Figure \ref{F10}.
	\setcounter{rhp}{8}
	\begin{rhp}\label{RHP9}
		Seek for an analytic matrix-valued function ${E_{\circ}}(z)$ defined on  $z\in \mathbb{C}\backslash\{\Sigma^{E_{\circ}}\}$ that satisfies the following conditions.
		\begin{itemize}
			\item ${E_{\circ}}(z)$ has continuous boundary values ${E_{\circ,\pm}}(z)$ on $z\in \Sigma^{E_{\circ}}$, which satisfy the jump conditions
			\begin{align*}%\label{eq111}
				\begin{aligned}
					E_{\circ,+}(z)=E_{\circ,-}(z)V^{({E_{\circ}})}(z;x,t),
				\end{aligned}
			\end{align*}
			where
			\begin{align}\label{eq112}
				V^{({E_{\circ}})}(z) :=
				\begin{cases} M_{\mathrm{out}}(z)V^{(2)}(z)M^{-1}_{\mathrm{out}}(z) &z\in \Sigma ^{{E_{\circ}}}\backslash \partial D( \pm \zeta_0,\kappa) \\
					M_{\mathrm{out}}(z)M_A(z(\xi))M^{-1}_{\mathrm{out}}(z) &z\in \partial D(  \zeta_0,\kappa) \\
					M_{\mathrm{out}}(z)M_B(z(\xi))M^{-1}_{\mathrm{out}}(z) &z\in \partial D( -\zeta_0,\kappa) .\end{cases} 
			\end{align}
			\item  The matrix  ${E_{\circ}}(z)$  exhibits the asymptotic behavior ${E_{\circ}}(z)\to I$, as $|z|\to \infty$.
		\end{itemize}
	\end{rhp}
	\begin{proof}
		For $z\in \Sigma ^{{E_{\circ}}}\backslash \partial D( \pm \zeta_0,\kappa)$, recalling the definition \eqref{eq73} and the jump condition \eqref{eq70}, we derive
		\begin{align*}
			\begin{aligned}E_{\circ,+}(z) &=M_{\mathrm{RHP},+}^{ (2) }(z)M^{-1}_{\mathrm{out}}(z)=M_{\mathrm{RHP},-}^{ (2) }(z)V^{(2)}(z)M^{-1}_{\mathrm{out}}(z)\\
				&=E_{\circ,-}(z) M_{\mathrm{out}}(z)V^{(2)}(z)M^{-1}_{\mathrm{out}}(z).\end{aligned}
		\end{align*}
		For $z\in \partial D( \zeta_0,\kappa)$, where $M_{\mathrm{RHP}}^{ (2) }(z)$ is jump-free, combining \eqref{eq110}, we compute
		\begin{align*}
			\begin{aligned}E_{\circ,+}(z) &=M_{\mathrm{RHP}}^{ (2) }(z)M^{-1}_{\mathrm{out}}(z)=E_{\circ,-}(z)M_{\mathrm{out}}(z)M_A(z(\xi))M^{-1}_{\mathrm{out}}(z) ,\end{aligned}
		\end{align*}
		Analogous computations apply for $z\in \partial D( -\zeta_0,\kappa)$, with $M_A$ replaced by $M_B$.
	\end{proof}
	
	\begin{prop}
		Let the complex variable be expressed as $z=a+ib$.  Given the constants $c_1<0$ and $c_2>0$ from Lemma \ref{Lemma4}, define $c_0:=\min \{-c_1|a\pm\zeta_0|^2,c_2\}$. Then, the following asymptotic estimate holds
		\begin{align}\label{eq113}
			\| V^{({E_{\circ}})}(z) -I\| _{L^{\infty }(\Sigma^{E_{\circ}})} =\begin{cases}\mathcal O( e^{-c_0 t }) &z\in \Sigma ^{{E_{\circ}}}\backslash \partial D( \pm \zeta_0,\kappa) \\
				\mathcal O( t^{-1/2}) &z\in \partial D( \pm \zeta_0,\kappa). \end{cases}
		\end{align}
	\end{prop}
	
	\begin{proof}
		For $z \in \Sigma^{E_{\circ}} \setminus \partial D(\pm \zeta_0, \kappa)$, we note that  $V^{(2)}\equiv I$ within $\cup_{z_k\in \mathcal{Z}\cup \overline{\mathcal{Z}}} D( z_{k},\kappa )$.  Then 
		\begin{align*}%\label{eq113}
			\begin{aligned}
				\| M_{\mathrm{out}}(z)&V^{(2)}(z)M^{-1}_{\mathrm{out}}(z)-1\| _{L^{\infty }(\Sigma^{E_{\circ}} \setminus \partial D(\pm \zeta_0, \kappa))} \\
				&=\| M_{\mathrm{out}}(z)(V^{(2)}(z)-1)M^{-1}_{\mathrm{out}}(z)\| _{L^{\infty }(\Sigma^{E_{\circ}} \setminus (\partial D(\pm \zeta_0, \kappa)\cup_{z_k\in \mathcal{Z}\cup \overline{\mathcal{Z}}} D( z_{k},\kappa )))} ,
			\end{aligned}
		\end{align*}
		and the estimation follows directly from recalling \eqref{eq74}. 
		
		Furthermore, since $M_{\mathrm{out}}(z)$ is bounded within $D(\pm \zeta_0, \kappa)$, there exists a constant $c > 0$ such that
		\begin{align*}
			\| V^{({E_{\circ}})}(z) -I\| _{L^{\infty }}\leq c\| M_A(z(\xi)) -I\| _{L^{\infty }},~~z\in \partial D( \zeta_0,\kappa),
		\end{align*} and 
		\begin{align*}
			\| V^{({E_{\circ}})}(z) -I\| _{L^{\infty }}\leq c\| M_B(z(\xi)) -I\| _{L^{\infty }},~~z\in \partial D( -\zeta_0,\kappa).
		\end{align*}
		Applying \eqref{eq108}-\eqref{eq109.1}, the stated results are obtained.
	\end{proof}
	It follows that the RHP \ref{RHP9} falls into the small-norm Riemann–Hilbert problem, ensuring the existence and uniqueness of its solution. Consequently, by the Beals-Coifman theory, we define $w(z) := V^{({E_{\circ}})}(z) - I$ and the operator
	\begin{align*}
		C_w:=C_-(\cdot w),
	\end{align*}
	which satisfies $\|C_w \|_{L^2(\Sigma^{E_{\circ}})}=\mathcal{O}(t^{-1/2})$.
	Then, ${E_{\circ}}(z)$ can be represented as 
	\begin{align}\label{eq114}
		{E_{\circ}}(z)=I+\frac{1}{2\pi i}\int_{\Sigma^{E_{\circ}}}\frac{(I+h(s))w(s)}{s-z}ds,
	\end{align}
	where $h(z)$ is the solution to the integral equation $(1-C_w)h=C_wI$. We proceed to estimate
	\begin{align}\label{eq115}
		\begin{aligned}
			\| h\| _{L^{2}( \Sigma ^{{E_{\circ}}}) }&=\| ( 1-C_{w}) ^{-1}C_wI\| _{L^{2}( \Sigma^{E_{\circ}}) }=\| \sum ^{\infty }_{k=0}( C_{w}) ^{k}C_{w} I\| _{L^{2}( \Sigma^ {E_{\circ}}) }\\
			&\leq \sum ^{\infty }_{k=0}\| C_w\| _{L^{2}(\Sigma^{E_{\circ}})}^{k}\| C_wI\| _{L^{2}(\Sigma^{E_{\circ}})} =\frac{ct^{-{1}/{2}}}{1-\| C_{w}\| _{L^{2}(\Sigma^{E_{\circ}})}}.
		\end{aligned}
	\end{align}
	
	From \eqref{eq114}, the asymptotic expansion of ${E_{\circ}}(z)$ as $|z|\to \infty$ is given by
	\begin{align*}
		{E_{\circ}}(z) :=I+\frac{E_{\circ,1}(x,t)}{z}+\mathcal O(z^{-2}) =I-\frac{1}{2\pi i z}\int _{\Sigma^{E_{\circ}}}( I+h(s)) w(s) ds+\mathcal O( z^{-2}) .
	\end{align*}
	To determine the asymptotic behavior of $E_{\circ,1}(x,t)$ as $t\to \infty$, we decompose
	\begin{align*}
		E_{\circ,1}(x,t)&=-\frac{1}{2\pi i}\left( \int _{\partial D(\pm\zeta_0,\kappa)}w(s)ds+\int_{ \Sigma ^{{E_{\circ}}}\backslash \partial D(\pm \zeta_0 ,\kappa)}w(s) ds+\int _{\Sigma^{E_{\circ}}}h(s) w(s) ds\right) \\
		&:=\mathrm{I+II+III}.
	\end{align*}
	By the exponential decay bound \eqref{eq113}, $|\mathrm{II}|$ is exponentially small as  $t \to \infty$. Applying H\"older's inequality and using estimates \eqref{eq113} and \eqref{eq115}, we deduce $|\mathrm{III}| \lesssim \mathcal{O}(t^{-1})$. Furthermore, using \eqref{eq112}, \eqref{eq108}, and \eqref{eq109}, the contour integral term $\mathrm{I}$ (oriented clockwise) is evaluated via residues, resulting in
	\begin{align}\label{eq116}
		\begin{aligned}
			E_{\circ,1} (x,t)= \frac{\sqrt{\beta}}{4\sqrt{t}} &\left( M_{\mathrm{out}}(\zeta_0) \begin{pmatrix}
				0 &i\beta_{12}\delta_0^A \\
				i\overline{\beta_{12}\delta_0^A}& 0
			\end{pmatrix} M^{-1}_{\mathrm{out}}(\zeta_0) \right. \\
			&+ \left. M_{\mathrm{out}}(-\zeta_0) \begin{pmatrix}
				0 &i\overline{\beta_{12}\delta_0^A} \\
				i\beta_{12}\delta_0^A& 0
			\end{pmatrix} M^{-1}_{\mathrm{out}}(-\zeta_0) \right) + \mathcal{O}(t^{-1}).
		\end{aligned}
	\end{align}
	Additionally, consider
	\begin{align*}
		{E_{\circ}}(\pm \mu/2) =I+ \frac{1}{2\pi i}\int_{\Sigma^{E_{\circ}}}\frac{(I+h(s))w(s)}{s\mp\mu/2}ds.
	\end{align*}
	Since $|s\pm\mu/2|\geq \kappa$ for $s\in \Sigma^{E_{\circ}}$ (avoiding poles at $ \pm\mu/2 $) and using the calculations for $E_{\circ,1}$, we obtain
	\begin{align}\label{eq117}
		\begin{aligned}
			{E_{\circ}}(\pm \mu/2) =I+\mathcal{O}(t^{-1/2}). 
		\end{aligned}
	\end{align}
	
	\subsection{Asymptotic analysis of the pure $\overline{\partial}$ problem}
	The resolution of the $\overline \partial-$problem \ref{RHP5} is equivalent to solving the integral equation
	\begin{align}\label{eq118}
		\begin{aligned}
			M ^{(3)}(z)= I-\frac{1}{\pi}\iint _{\mathbb{C} }\frac{ \overline{\partial }M^{(3) }(s) }{s-z}dA(s),
		\end{aligned}
	\end{align}
	where $dA(s)$ denotes the Lebesgue area measure on $\mathbb{C}$. Substituting the $\overline \partial-$equation \eqref{eq72} into the right-hand side, we reformulate it via the solid Cauchy operator as
	\begin{align*}%\label{eq119}
		\begin{aligned}
			(1-\mathbb{S})M ^{(3)}(z)= I,
		\end{aligned}
	\end{align*}
	where the operator $\mathbb{S}$ is defined by
	$$(\mathbb{S} f) (z):= -\frac{1}{\pi}\iint _{\mathbb{C} }\frac{ f(s)W(s) }{s-z}dA(s).$$
	
	\begin{prop}\label{prop8}
		For $W(z)$ defined in $\overline \partial-$problem \ref{RHP5} and $r(z)\in H^{1,1}(\mathbb R)$, the operator norm satisfies
		\begin{align}\label{eq119.1}
			\lim _{t\rightarrow \infty }\| \mathbb{S}\| _{L^{\infty }( \mathbb{C} ) }=0.
		\end{align}
		Consequently, the operator $1-\mathbb{S}$ is invertible for sufficiently large $t$.
	\end{prop}
	
	\begin{proof}
		Since $\overline \partial R^{(2)}(z)=0$ for $z\in \Omega_{13}$ and $ z\in \cup_{z_k\in \mathcal{Z}\cup \overline{\mathcal{Z}}} D( z_{k},\kappa )$, and the function $M^{(2)}_{\mathrm{RHP}}(z)$ is bounded in the region $$\Omega ^{\# }:= (\cup ^{12}_{k=1}\Omega _{k}) \backslash ( \cup_{z_{j}\in \mathcal Z\cup \overline {\mathcal Z}}D( z_{j},\kappa ) ), $$ we assume $\|f\|_{L^\infty(\mathbb{C})}=1$ and then yield 
		\begin{align}\label{eq120}
			\begin{aligned}\|\mathbb S\| _{L^{\infty }(\mathbb C) }&=\left| -\frac{1}{\pi}\iint _{\mathbb{C}}\frac{  f(s) M^{(2)}_{\mathrm{RHP} }(s) \overline{\partial }R ^{(2) }(s)(M^{(2)}_{\mathrm{RHP} }(s))^{-1}}{s-z}dA(s)\right| \\
				&\lesssim \frac{1}{\pi}\|M^{(2)}_{\mathrm{RHP}}\|_{L^\infty(\Omega^\#)}\|(M^{(2)}_{\mathrm{RHP}})^{-1}\|_{L^\infty(\Omega^\#)}\iint _{\Omega^ {\#}}\frac{| \overline{\partial }R^{(2)}(s) | }{| s-z| }dA(s). \end{aligned} 
		\end{align}
		Let $s:=a+ib$ and $z:=\alpha+i\beta$.  The behavior of this integral is governed by the explicit expression of $\overline{\partial }R^{(2)}$ given in \eqref{eq68} and further estimates in Lemma \ref{lemma1}.
		
		For $s \in \Omega_1$, in view of $\operatorname{Re}[i \theta(z)]$ in \eqref{eq32}, we note that the inequality $G(a, b) \leq G(b + \zeta_0, b) \leq 0$ holds, since $G(a,b)$ is monotonically decreasing with respect to $a\in (b+\zeta_0,\infty)$ for fixed $b\in (0,\infty)$. Thus, we have
		\begin{align*}
			\begin{aligned}\iint _{\Omega_1}\frac{| \overline{\partial }R^{(2)}(s) | }{| s-z| }dA(s) \leq \mathrm{I+II},\end{aligned} 
		\end{align*}
		where 
		\begin{align*}
			\mathrm{I}&:=\int _{0}^{\infty }e^{2tbG( b+\zeta_0,b) }\int _{b+\zeta_0}^{\infty }\frac{| r'( a) | }{| s-z| }dadb\\
			\mathrm{II}&:=\int _{0}^{\infty }e^{2tbG( b+\zeta_0,b) }\int _{b+\zeta_0}^{\infty }\frac{1 }{\sqrt{a-\zeta_0}| s-z| }dadb.
		\end{align*}
		For $q>1$, it is easy to see
		\begin{align*}
			\|{(s-z)}^{-1}\| _{L_{a}^{q}( b+\zeta_{0},\infty ) }\lesssim | b-\beta | ^{1/q-1}
		\end{align*}
		holds uniformly in $z$. Choosing $q=2$ and applying Hölder's inequality to term $\mathrm{I}$ yields
		\begin{align}\label{eq121}
			\mathrm{I}\lesssim\| r\| _{H^{1,1}(\mathbb{R} )}\left(\int _{0}^{\beta }e^{2tbG( b+\zeta_0,b) }( \beta -b) ^{-1/2}db+\int _{\beta}^\infty e^{2tbG( b+\zeta_0,b)}( b-\beta ) ^{-1/2}db\right),
		\end{align}
		where both integrals on the right-hand side are uniformly convergent for $t>0$.
		For the dual exponent $p>2$ (with $q\in (1,2)$), the norm
		\begin{align}\label{eq121.1}
			\|(a-\zeta_0)^{-1/2}\|_{L^p(b+\zeta_{0},\infty)}=\left( \frac{1}{({p}/{2}-1)b^{{p}/{2}-1}}\right) ^{{1}/{p}}
		\end{align}
		implies the bound for  term $\mathrm{II}$
		\begin{align*}%\label{eq122}
			\mathrm{II}\lesssim \int _{0}^{\beta }e^{2tbG( b+\zeta_0,b) }b^{1/p-1/2}( \beta -b) ^{1/q-1}db+\int _{\beta }^{\infty }e^{2tbG( b+\zeta_0,b)  }b^{1/p-1/2}( b-\beta ) ^{1/q-1}db.
		\end{align*}
		Thereby, the integrals $\mathrm{I}$ and $\mathrm{II}$ converge uniformly for $t>0$. 
		
		For $s\in \cup_{k=9}^{12}\Omega_k$, an argument analogous to that of \eqref{eq121} establishes the uniform convergence of 
		\begin{align*}%\label{eq122}
			\iint_{\Omega_k}  \frac{|\overline \partial \chi_{\mathcal{Z}}(s)| }{| s-z| }e^{\pm2tbG(a,b)}dadb,
		\end{align*}
		with respect to $t>0$.
		
		The same reasoning applies to all $s\in \cup_{k=2}^{12}\Omega_k$ ensuring the uniform convergence of the integral in \eqref{eq120}. It follows that, due to the exponential decay of $e^{\pm2tbG( b+\zeta_0,b) }$ as $t\to \infty $, the limit in \eqref{eq119.1} holds.
	\end{proof}
	
	From \eqref{eq118}, the asymptotic expansion of $M^{(3)}(z)$ as $|z|\to \infty$ is given by
	\begin{align*}
		M^{(3)}(z) :=I+\frac{M^{(3)}_1(x,t)}{z}+\mathcal O(z^{-2}) =I+\frac{1}{z\pi }\iint _{\mathbb{C}} M^{(3)}(s)W(s)dA(s)+\mathcal O( z^{-2}) .
	\end{align*}
	To determine the asymptotic behavior of $M^{(3)}_1(x,t)$ as $t\to \infty$, we estimate
	\begin{align*}%\label{eq123}
		| M_{1}^{( 3)} | \leq \iint _{\Omega ^\# }| \overline{\partial }R^{ (2) }( s) | dA(s).
	\end{align*}
	For $z\in \Omega_1$, applying Lemma \ref{lemma1} again, we have 
	\begin{align*}%\label{eq124}
		\begin{aligned}
			\iint _{\Omega _1}| \overline{\partial }R^{(2)}(s)| dA(s) \leq \mathrm{I+II},
		\end{aligned}		
	\end{align*}
	where
	\begin{align*}
		\mathrm{I}&:=\int _{0}^{\infty }\int _{b+\zeta_0}^{\infty }{| r'( a) | }e^{2tbG( a,b) }dadb\\
		\mathrm{II}&:=\int _{0}^{\infty }\int _{b+\zeta_0}^{\infty }\frac{1 }{\sqrt{a-\zeta_0} }e^{2tbG( a,b) }dadb.
	\end{align*}
	It is reasonable to assume the existence of ${\kappa}_{\circ} >0$ sufficiently small such that 
	\begin{align*}
		G(a,b)+{\kappa_{\circ}}(a-b-\zeta_0)\leq G(b+\zeta_0,b).
	\end{align*}
	Thus,  by applying H\"older's inequality, one yield
	\begin{align*}%\label{124.1}
		\begin{aligned}
			\mathrm{I} &\leq \| r\| _{H^{1,1}( \mathbb{R}) }\int _{0}^{\infty }\left\|e^{2tbG(a,b)}\right\|_{L^2_a(b+\zeta_0,\infty)}db\\
			&\leq  \| r\| _{H^{1,1}( \mathbb{R}) }\int _{0}^{\infty }\left\|e^{2tb(G(b+\zeta_0,b)-{\kappa_{\circ}}(a-b-\zeta_0) )}\right\|_{L^2_a(b+\zeta_0,\infty)}db\\
			&\lesssim t^{-1/2}\int _{0}^{\infty }e^{2tbG(b+\zeta_0,b)}b^{-1/2}db
		\end{aligned}
	\end{align*}
	For term $\mathrm{II}$, recall \eqref{eq121.1} with $p>2$, we have
	\begin{align*}
		\begin{aligned}
			\mathrm{II} &\lesssim \int _{0}^{\infty }\frac{1}{b^{1/2-1/p}}\left\|e^{2tbG(a,b)}\right\|_{L^q_a(b+\zeta_0,\infty)}db\\
			&\lesssim t^{-1/q}\int _{0}^{\infty } {b^{1/p-1/2-1/q}} e^{2tbG(b+\zeta_0,b)}  db
		\end{aligned}
	\end{align*}
	where $q\in (4/3,2)$. 
	
	Through analogous methods, estimating the contributions from the remaining regions $s\in \cup_{k=2}^{12}\Omega_k$, we establish the bound for the leading term
	\begin{align}\label{eq125}
		| M_{1}^{( 3)} | \lesssim \iint _{\Omega ^\# }| \overline{\partial }R^{ (2) }( s) | dA(s)\lesssim\mathcal {O}(t^{-\gamma_0}),
	\end{align}
	where $\gamma_0>1/2$ arises from $t^{-1/2}$ and additional decay provided by  $e^{2tbG(b+\zeta_0,b)}$. 
	
	Moreover, consider
	\begin{align*}
		M^{ (3) }( \pm \mu/{2}) =I-\frac{1}{\pi}\iint _{\mathbb{C} }\frac{\overline \partial M ^{ (3) }(s) }{s{\mp }\mu/2}dA(s),
	\end{align*}
	which, as $|s\pm\mu/2|\geq \kappa$ for $s\in \Omega^\#$, implies 
	\begin{align}\label{eq126}
		\begin{aligned}
			|M^{ (3) }( \pm \mu/{2})-I| \lesssim\iint _{\Omega^\#}{\overline \partial R ^{ (2) }(s) }dA(s)
			\lesssim\mathcal{O}({t^{-\gamma_0}}).
		\end{aligned}
	\end{align}
	
	\section{Long time asymptotics for RMB equations}\label{S5}
	In this section, we investigate the long-time asymptotic behavior of solutions to the RMB equations \eqref{eq01} within the soliton region $-1/\mu^2<x/t<0$ with $\mu \in (0,1]$ as $t\to \infty$. 
	
	Through inverse transformations \eqref{eq35}, \eqref{eq63}, \eqref{eq69}, and \eqref{eq73}, the solution to RHP \ref{RHP1} can be represented in the form
	\begin{align*}
		M( z;x,t) =M^{(3)}( z;x,t) E_\circ( z;x,t) M_{\mathrm {out}}( z;x,t) (R^{(2)}(z)) ^{-1}F^{\sigma _{3}}(z).
	\end{align*}
	Then, by applying the reconstruction formulae provided in Theorem \ref{theo2}, we recover the solutions of the RMB equations. Restricting to $z\in \Omega_{13}$, where  $R^{(2)}(z)\equiv I$, and utilizing the asymptotic expansion
	\begin{align*}
		F^{\sigma _{3}}(z)=I+\frac{1}{z}\left( \sum _{k\in \bigtriangleup }2z_{k}-i\int _{\mathbb{R} \backslash [ -\zeta_0,\zeta_0] }\nu( s) ds\right) \sigma _{3}+\mathcal O(z^{-2}),
	\end{align*}
	the electric field is given by
	\begin{align*}
		E(x,t) =-4i\left( M^{(3)}_1(x,t)  +E_{\circ,1}(x,t)  +M_{\mathrm {out},1}(x,t)  \right)_{12}.
	\end{align*}
	Combining  Proposition \ref{prop5}, \eqref{eq116}, and \eqref{eq125}, the asymptotic expansion of $E(x,t)$  as $t\to \infty$ reads
	\begin{align*}
		E(x,t) =E(x,t;\sigma_\delta)+t^{-1/2}\sqrt{\beta} (f_1(x,t)+f_2(x,t)) +\mathcal{O} (t^{-\gamma_\circ}),
	\end{align*}
	where the constant $ \beta$ is given in \eqref{eq87.1} and $\gamma_\circ>1/2$. Moreover
	\begin{align*}
		f_1(x,t)&:=	\left( M_{\mathrm{out}}(\zeta_0) \begin{pmatrix}
			0 &\beta_{12}\delta_0^A \\
			\overline{\beta_{12}\delta_0^A}& 0
		\end{pmatrix} M^{-1}_{\mathrm{out}}(\zeta_0) \right) _{12}\\
		f_2(x,t)&:=\left(M_{\mathrm{out}}(-\zeta_0) \begin{pmatrix}
			0 &\overline{\beta_{12}\delta_0^A} \\
			\beta_{12}\delta_0^A& 0
		\end{pmatrix} M^{-1}_{\mathrm{out}}(-\zeta_0)\right) _{12},
	\end{align*}
	where $M_{\mathrm {out}}(z)=M^{\bigtriangleup }( z;\sigma _{\delta }) $ is the solution to RHP \ref{RHP10} with the scattering data $\sigma_{\delta }$ as given in \eqref{eq75.1}. The modulus and argument of $\beta_{12}\delta_0^A$ are given explicitly in \eqref{eq109.1} and \eqref{eq109.2}, respectively.
	
	Additionally, noting that $\pm \mu/2\in \Omega_{13}$, we have 
	\begin{align*}
		\rho( \pm \mu /{2};x,t) =\widehat{\rho }( \pm \mu /{2};x,t) \sigma _{3}\widehat{\rho }^{-1}( \pm \mu /{2};x,t) ,
	\end{align*}
	where
	\begin{align*}
		\widehat{\rho}( \pm \mu /{2};x,t) :=M^{(3)}( \pm \mu /{2};x,t) E_{\circ}( \pm \mu /{2};x,t) M_{\mathrm {out}}( \pm \mu /{2};x,t).
	\end{align*}
	Applying Proposition \ref{prop5} together with \eqref{eq117} and \eqref{eq126},  we obtain the following long-time asymptotics for the components of the Bloch vector
	\begin{align*}
		\begin{aligned}
			&s(x,t) =s( x,t;\sigma _{\delta }) +\mathcal O( t^{-{1}/{2}}) 
			&u(x,t) =u( x,t;\sigma _{\delta })+\mathcal O( t^{-{1}/{2}})  \\
			&r(x,t) =r( x,t;\sigma _{\delta }) +\mathcal  O( t^{-{1}/{2}}) .
		\end{aligned}
	\end{align*}
	
	Within the cone $C(x_1,x_2,v_1,v_2)$ defined in \eqref{eq79}, Corollary \ref{Corollary1} establishes the reduction to observable $N(\mathcal { J}) $-soliton and $N(\mathcal { J})$-kink
	\begin{align*}
		\begin{aligned}
			&E(x,t;\sigma_\delta) \to E(x,t;\sigma_\delta(\mathcal{J}))  ~~&s(x,t;\sigma_\delta) \to s(x,t;\sigma_\delta(\mathcal{J}))\\
			&u( x,t;\sigma _{\delta })  \to u(x,t;\sigma_\delta(\mathcal{J}))~~&~r( x,t;\sigma _{\delta }) \to r(x,t;\sigma_\delta(\mathcal{J})),
		\end{aligned}
	\end{align*}
	in which the scattering data $\sigma_\delta(\mathcal J):=\{\{(z_k,c_{\delta,k}= c_k\delta^{-2}(z_k))\}_{k=1}^N:z_k\in \mathcal {Z\cap J} \}$. Correspondingly, $M_{\mathrm {out}}(z)=M^{\bigtriangleup _{\mathcal J}}( z;\sigma _{\delta }(\mathcal J))$  is taken as the solution to RHP \ref{RHP10} with index set $\bigtriangleup=\bigtriangleup _{\mathcal J}$ and the scattering data $\sigma_\delta(\mathcal J)$.
	
\section*{Acknowledgements}
	The research is supported by the National Natural Science Foundation of China (Grant No. 12471239), and the Guangdong Basic and Applied Basic Research Foundation (Grant No. 2024A1515013106).

	\bibliographystyle{abbrv} % 样式文件，可选择plain, unsrt, alpha, IEEEtran等abbrv
	\bibliography{references}

\begin{thebibliography}{10}

\bibitem{Aiyer1983}
R.~N. Aiyer.
\newblock Hamiltonian and recursion operator for the reduced {M}axwell-{B}loch
  equations.
\newblock {\em J. Phys. A: Math. Gen.}, 16:1809--1811, 1983.

\bibitem{beals}
R.~Beals and R.~R. Coifman.
\newblock Scattering and inverse scattering for first order systems.
\newblock {\em Comm. Pure Appl. Math.}, 37:39--90, 1984.

\bibitem{MR954382}
R.~Beals, P.~Deift, and C.~Tomei.
\newblock {\em Direct and inverse scattering on the line}.
\newblock American Mathematical Society: Rhode Island, 1988.

\bibitem{Biondini2021}
G.~Biondini, S.~Li, and D.~Mantzavinos.
\newblock Long-time asymptotics for the focusing nonlinear {S}chr\"odinger
  equation with nonzero boundary conditions in the presence of a discrete
  spectrum.
\newblock {\em Comm. Math. Phys.}, 382:1495--1577, 2021.

\bibitem{MR3795020}
M.~Borghese, R.~Jenkins, and K.~D. T.-R. McLaughlin.
\newblock Long time asymptotic behavior of the focusing nonlinear
  {S}chr\"odinger equation.
\newblock {\em Ann. Inst. H. Poincar\'e Anal. Non Lin\'eaire}, 35:887--920,
  2018.

\bibitem{Boutet2025}
A.~Boutet~de Monvel, J.~Lenells, and D.~Shepelsky.
\newblock The focusing {NLS} equation with step-like oscillating background:
  {A}symptotics in a transition zone.
\newblock {\em J. Differential Equations}, 429:747--801, 2025.

\bibitem{Bullough1979}
R.~K. Bullough, P.~M. Jack, P.~W. Kitchenside, and R.~Saunders.
\newblock Solitons in laser physics.
\newblock {\em Phys. Scripta}, 20:364--381, 1979.

\bibitem{Caudrey1974}
P.~J. Caudrey, J.~C. Eilbeck, and J.~D. Gibbon.
\newblock Exact multisoliton solution of the reduced {M}axwell-{B}loch
  equations of non-linear optics.
\newblock {\em J. Inst. Math. Appl.}, 14:375--386, 1974.

\bibitem{Caudrey1973}
P.~J. Caudrey, J.~C. Eilbeck, J.~D. Gibbon, and R.~K. Bullough.
\newblock Exact multisoliton solutions of inhomogeneously broadened
  self-induced transparency equations.
\newblock {\em J. Phys. A: Math., Nucl. Gen.}, 6:L53--L56, 1973.

\bibitem{Charlier2023}
C.~Charlier, J.~Lenells, and D.-S. Wang.
\newblock The ``good'' {B}oussinesq equation: long-time asymptotics.
\newblock {\em Anal. PDE}, 16:1351--1388, 2023.

\bibitem{Chen2021}
G.~Chen and J.~Liu.
\newblock Soliton resolution for the focusing modified {K}d{V} equation.
\newblock {\em Ann. Inst. H. Poincar\'e - Anal. Non Lin\'eaire}, 38:2005--2071,
  2021.

\bibitem{Cuccagna2016}
S.~Cuccagna and R.~Jenkins.
\newblock On the asymptotic stability of {$N$}-soliton solutions of the
  defocusing nonlinear {S}chr\"odinger equation.
\newblock {\em Comm. Math. Phys.}, 343:921--969, 2016.

\bibitem{MR1207209}
P.~Deift and X.~Zhou.
\newblock A steepest descent method for oscillatory {R}iemann-{H}ilbert
  problems. {A}symptotics for the {MK}d{V} equation.
\newblock {\em Ann. of Math.}, 137:295--368, 1993.

\bibitem{MR1989226}
P.~Deift and X.~Zhou.
\newblock Long-time asymptotics for solutions of the {NLS} equation with
  initial data in a weighted {S}obolev space.
\newblock {\em Comm. Pure Appl. Math.}, 56:1029--1077, 2003.

\bibitem{Dieng2019}
M.~Dieng, K.~D.-R. McLaughlin, and P.~D. Miller.
\newblock Dispersive asymptotics for linear and integrable equations by the
  {$\overline{\partial}$} steepest descent method.
\newblock In {\em Nonlinear dispersive partial differential equations and
  inverse scattering}, pages 253--291. Springer: New York, 2019.

\bibitem{Dodd}
R.~K. Dodd, J.~C. Eilbeck, J.~D. Gibbon, and H.~C. Morris.
\newblock {\em Solitons and Nonlinear Wave Equations}.
\newblock Academic Press: London, 1982.

\bibitem{Eilbeck1972a}
J.~C. Eilbeck.
\newblock Reflection of short pulses in linear optics.
\newblock {\em J. Phys. A: Gen. Phys.}, 5:1355--1363, 1972.

\bibitem{Eilbeck1972b}
J.~C. Eilbeck and R.~K. Bullough.
\newblock The method of characteristics in the theory of resonant or
  nonresonant nonlinear optics.
\newblock {\em J. Phys. A: Gen. Phys.}, 5:820--829, 1972.

\bibitem{J.C.Eilbeck_1973}
J.~C. Eilbeck, J.~D. Gibbon, P.~J. Caudrey, and R.~K. Bullough.
\newblock Solitons in nonlinear optics {I. A} more accurate description of the
  $2\pi$ pulse in self-induced transparency.
\newblock {\em J. Phys. A: Math., Nucl. Gen.}, 6:1337--1347, 1973.

\bibitem{Gibbon1973}
J.~D. Gibbon, P.~J. Caudrey, R.~K. Bullough, and J.~C. Eilbeck.
\newblock An {$N$}-soliton solution of a nonlinear optics equation derived by a
  general inverse method.
\newblock {\em Lett. Nuovo Cimento}, 8:775--779, 1973.

\bibitem{Grauel1984}
A.~Grauel.
\newblock Soliton connection of the reduced {M}axwell-{B}loch equations.
\newblock {\em Lett. Nuovo Cimento}, 41:263--268, 1984.

\bibitem{Grauel1986}
A.~Grauel.
\newblock The {P}ainlev\'e{} test, {B}\"acklund transformation and solutions of
  the reduced {M}axwell-{B}loch equations.
\newblock {\em J. Phys. A: Math. Gen.}, 19:479--484, 1986.

\bibitem{Huang2019}
L.~L. Huang and Y.~Chen.
\newblock Localized excitations and interactional solutions for the reduced
  {M}axwell-{B}loch equations.
\newblock {\em Commun. Nonlinear Sci. Numer. Simul.}, 67:237--252, 2019.

\bibitem{Jenkins2018}
R.~Jenkins, J.~Liu, P.~Perry, and C.~Sulem.
\newblock Soliton resolution for the derivative nonlinear {S}chr\"odinger
  equation.
\newblock {\em Comm. Math. Phys.}, 363:1003--1049, 2018.

\bibitem{MR475310}
G.~L. Lamb, Jr.
\newblock Analytical descriptions of ultrashort optical pulse propagation in a
  resonant medium.
\newblock {\em Rev. Modern Phys.}, 43:99--124, 1971.

\bibitem{Li2024}
S.~Li and P.~D. Miller.
\newblock On the {M}axwell-{B}loch system in the sharp-line limit without
  solitons.
\newblock {\em Comm. Pure Appl. Math.}, 77:457--542, 2024.

\bibitem{Maimistov2001}
A.~I. Maimistov.
\newblock Completely integrable models of nonlinear optics.
\newblock {\em Pramana - J. Phys.}, 57:953--968, 2001.

\bibitem{Maimistov}
A.~I. Maimistov and A.~M. Basharov.
\newblock {\em Nonlinear optical waves}.
\newblock Kluwer Academic Publishers: Dordrecht, 1999.

\bibitem{PhysRev.183.457}
S.~L. McCall and E.~L. Hahn.
\newblock Self-induced transparency.
\newblock {\em Phys. Rev.}, 183:457--485, 1969.

\bibitem{McLaughlin2006}
K.~T.-R. McLaughlin and P.~D. Miller.
\newblock The {$\overline{\partial}$} steepest descent method and the
  asymptotic behavior of polynomials orthogonal on the unit circle with fixed
  and exponentially varying nonanalytic weights.
\newblock {\em Int. Math. Res. Pap.}, 2006:48673, 1--77, 2006.

\bibitem{McLaughlin2008}
K.~T.-R. McLaughlin and P.~D. Miller.
\newblock The {$\overline{\partial}$} steepest descent method for orthogonal
  polynomials on the real line with varying weights.
\newblock {\em Int. Math. Res. Not.}, 2008:rnn075, 1--66, 2008.

\bibitem{Pelinovsky2017}
D.~E. Pelinovsky and Y.~Shimabukuro.
\newblock Existence of global solutions to the derivative {NLS} equation with
  the inverse scattering transform method.
\newblock {\em Int. Math. Res. Not.}, 2018:5663--5728, 2018.

\bibitem{Wadati19731}
M.~Wadati.
\newblock The modified {K}orteweg-de {V}ries equation.
\newblock {\em J. Phys. Soc. Japan}, 34:1289--1296, 1973.

\bibitem{Wadati1972}
M.~Wadati and M.~Toda.
\newblock The exact {$N$}-soliton solution of the {Korteweg-de V}ries equation.
\newblock {\em J. Phys. Soc. Japan}, 32:1403--1411, 1972.

\bibitem{Wei2018}
J.~Wei, X.~Wang, and X.~Geng.
\newblock Periodic and rational solutions of the reduced {M}axwell-{B}loch
  equations.
\newblock {\em Commun. Nonlinear Sci. Numer. Simul.}, 59:1--14, 2018.

\bibitem{Xu2013}
S.~Xu, K.~Porsezian, J.~He, and Y.~Cheng.
\newblock Circularly polarized few cycle optical rogue waves:{ Rotating reduced
  Maxwell-Bloch equations}.
\newblock {\em Phys. Rev. E}, 88:062925, 2013.

\bibitem{Yang2022}
Y.~Yang and E.~Fan.
\newblock On the long-time asymptotics of the modified {C}amassa-{H}olm
  equation in space-time solitonic regions.
\newblock {\em Adv. Math.}, 402:108340, 1--78, 2022.

\bibitem{Zakharov}
V.~E. Zakharov and A.~B. Shabat.
\newblock Exact theory of two-dimensional self-focusing and one-dimensional
  self-modulation of waves in nonlinear media.
\newblock {\em Soviet Physics JETP}, 34:62--69, 1972.

\bibitem{MR1000732}
X.~Zhou.
\newblock The {R}iemann-{H}ilbert problem and inverse scattering.
\newblock {\em SIAM J. Math. Anal.}, 20:966--986, 1989.

\bibitem{zhou}
X.~Zhou.
\newblock {$L^2$}-{S}obolev space bijectivity of the scattering and inverse
  scattering transforms.
\newblock {\em Comm. Pure Appl. Math.}, 51:697--731, 1998.

\end{thebibliography}
\end{document}